\documentclass[11pt]{amsart}

\usepackage[OT2, T1]{fontenc}
\usepackage{url}
\usepackage{amsmath}
\usepackage{stackengine}
\usepackage{array}

\usepackage{graphicx}
\usepackage{amsfonts}
\usepackage{amssymb}
\usepackage{amstext}
\usepackage{amsthm}
\usepackage{enumitem}
\usepackage{bm}
\usepackage{hyperref}
\usepackage{colonequals}
\usepackage{enumitem}
\usepackage[alphabetic,lite]{amsrefs}
\usepackage{cleveref}
\usepackage[all,cmtip]{xy}
\usepackage{fullpage}
\usepackage{appendix}
\usepackage{scalerel}
\usepackage{mathrsfs}
\usepackage{comment}
\usepackage{tikz-cd}
\usepackage{tikz}
\usepackage{galois}
\usepackage{dynkin-diagrams}
\usepackage{amssymb} 
\def\acts{\curvearrowright}

\newcommand{\subscript}[2]{$#1 _ #2$}
\numberwithin{equation}{subsection}

\newtheorem{theorem}{Theorem}
\newtheorem{lemma}[theorem]{Lemma}
\newtheorem{proposition}[theorem]{Proposition}
\newtheorem{corollary}[theorem]{Corollary}
\newtheorem{setup}[theorem]{Setup}

\newtheorem{conj}[theorem]{Conjecture}

\theoremstyle{definition}
\newtheorem{construction}[theorem]{Construction}
\newtheorem{defn}[theorem]{Definition}

\theoremstyle{remark}
\newtheorem{remark}[theorem]{Remark}

\newtheorem{claim}[theorem]{Claim}
\newtheorem*{assumption}{Assumption}
\newtheorem{example}[theorem]{Example}

\numberwithin{theorem}{section}

\newcommand{\bA}{\mathbb{A}}
\newcommand{\bB}{\mathbb{B}}
\newcommand{\bC}{\mathbb{C}}
\newcommand{\bD}{\mathbb{D}}

\newcommand{\bF}{\mathbb{F}}
\newcommand{\bG}{\mathbb{G}}
\newcommand{\bH}{\mathbb{H}}

\newcommand{\bK}{\mathbb{K}}
\newcommand{\bL}{\mathbb{L}}
\newcommand{\bM}{\mathbb{M}}

\newcommand{\bQ}{\mathbb{Q}}
\newcommand{\bR}{\mathbb{R}}

\newcommand{\bV}{\mathbb{V}}
\newcommand{\bW}{\mathbb{W}}
\newcommand{\bZ}{\mathbb{Z}}

\newcommand{\bbA}{\mathbf{A}}

\newcommand{\bbD}{\mathbf{D}}
\newcommand{\bbI}{\mathbf{I}}
\newcommand{\bbJ}{\mathbf{J}}

\newcommand{\bbH}{\mathbf{H}}
\newcommand{\bbK}{\mathbf{K}}

\newcommand{\bbQ}{\mathbf{Q}}

\newcommand{\bbT}{\mathbf{T}}

\newcommand{\cA}{\mathcal{A}}

\newcommand{\cC}{\mathcal{C}}
\newcommand{\cD}{\mathcal{D}}

\newcommand{\cM}{\mathcal{M}}
\newcommand{\cN}{\mathcal{N}}

\newcommand{\cQ}{\mathcal{Q}}

\newcommand{\cW}{\mathcal{W}}
\newcommand{\cX}{\mathcal{X}}

\newcommand{\Sh}{\mathcal{S}h}

\newcommand{\et}{{\text{\'et}}}
\newcommand{\dR}{{\mathrm{dR}}}
\newcommand{\cris}{{\mathrm{cris}}}

\newcommand{\ord}{{\mathrm{ord}}}
\newcommand{\can}{{\mathrm{can}}}

\newcommand{\inj}{\hookrightarrow}

\newcommand{\la}{\langle}
\newcommand{\ra}{\rangle}

\DeclareMathOperator{\FIsoc}{\textbf{F-Isoc}}
\DeclareMathOperator{\LS}{\textbf{LS}}
\newcommand{\Fpbar}{\bF}
\DeclareMathOperator{\GL}{GL}

\DeclareMathOperator{\GSp}{GSp}

\DeclareMathOperator{\SO}{SO}
\DeclareMathOperator{\GSpin}{GSpin}

\DeclareMathOperator{\Spf}{Spf}

\DeclareMathOperator{\Hg}{Lg}
\DeclareMathOperator{\SSU}{SU}

\DeclareMathOperator{\Lv}{Lv}
\DeclareMathOperator{\Gal}{Gal}
\DeclareMathOperator{\End}{End}
\DeclareMathOperator{\Hom}{Hom}
\DeclareMathOperator{\Aut}{Aut}
\DeclareMathOperator{\Lie}{Lie}
\DeclareMathOperator{\Res}{Res}
\DeclareMathOperator{\Spec}{Spec}
\DeclareMathOperator{\Ext}{Ext}
\DeclareMathOperator{\ad}{ad}
\DeclareMathOperator{\Span}{Span}
\DeclareMathOperator{\der}{der}
\DeclareMathOperator{\Fil}{Fil}

\DeclareMathOperator{\diag}{diag}

\DeclareMathOperator{\Id}{Id}

\DeclareMathOperator{\Tr}{Tr}
\DeclareMathOperator{\coker}{coker}

\DeclareMathOperator{\rk}{rk}

\DeclareMathOperator{\gr}{gr}
\DeclareMathOperator{\im}{im}

\DeclareMathOperator{\red}{red}

\DeclareMathOperator{\Vect}{Vect}

\DeclareMathOperator{\loc}{loc}

\DeclareMathOperator{\Cl}{Cl}

\DeclareMathOperator{\KS}{KS}
\DeclareMathOperator{\SU}{U}
\DeclareMathOperator{\SUU}{SU}

\DeclareMathOperator{\Def}{Def}

\DeclareMathOperator{\Spin}{Spin}
\DeclareMathOperator{\UU}{U}

\DeclareMathOperator{\Orb}{Orb}

\DeclareMathOperator{\MTT}{MT}

\DeclareMathOperator{\formalBr}{Br}

\begin{document}
\title{$p$-adic monodromy and mod $p$ unlikely intersections, I
}
\author{Ruofan Jiang}
\address{Dept.\ of Mathematics, University of California, Berkeley}
\email{ruofanjiang@berkeley.edu}
\begin{abstract}
We formulate characteristic $p$ analogues of the Mumford--Tate and the André--Oort conjectures for ordinary mod $p$ Shimura varieties of Hodge type, and set up general frameworks for studying them. We prove the two conjectures for (subvarieties of) arbitrary products of GSpin Shimura varieties, by reducing them, via a notion of linearity for mod $p$ Shimura varieties, to a third conjecture of Ax--Schanuel type. Along the way, we solve Chai's Tate-linear conjecture for products of GSpin Shimura varieties, and reveal an intimate relation among the four conjectures mentioned above. Our proof uses Crew's parabolicity conjecture which is recently proven by D'Addezio. 
\end{abstract}
\maketitle
\tableofcontents
\section{Introduction}
For a smooth projective variety $Y$ over a number field, the Mumford--Tate conjecture states that the base change to $\bQ_l$ of the Mumford--Tate group has the same neutral component with the $l$-adic étale monodromy group of $Y$. When $Y$ is an abelian variety, there are an abundance of results establishing many special cases, see \cites{Serre1,Chi,Ta1,Ta2,Pre,V08} (and there are more results not mentioned). Even though, the conjecture remains wide open for general abelian varieties, even for abelian fourfolds. 

The André--Oort conjecture arises in the relatively new field of unlikely intersections. It is related to the distribution of special subvarieties in a Shimura variety. It states that a subvariety of a Shimura variety that contains a Zariski dense collection of special subvarieties must itself be special. The
conjecture for Shimura varieties of Hodge type was proved by Tsimmerman (\cite{T18}). A proof for the conjecture in full generality has been announced recently, see \cite{P22}.

The main results of this paper concern characteristic $p$ analogues of the two conjectures. The reason that we treat two seemingly different conjectures in one paper is that they turn out to be much more intimately related to each other (and to a third conjecture of Ax--Schanuel type, see Conjecture~\ref{conj:AxSchanuel}) in the characteristic $p$ setting than their classical counterparts.

Before anything can be said, there are apparent difficulties in formulating the analogues. Let $Y$ be a smooth projective variety over a finitely generated field of characteristic $p$. Any attempt to formulating the Mumford--Tate conjecture for $Y$ will be obstructed by a lack of the notion of Hodge structure or  Mumford--Tate group. Of course, one can ask a weaker question (cf. \cite{LP95}):\\

\textit{Do the $l$-adic étale monodromy groups of $Y$ admit rational models over $\bQ$? }\\

If $Y$ is an ordinary principally polarized abelian variety, one do can  define a characteristic $p$ analogue of Mumford--Tate group, and formulate a characteristic $p$  analogue of the Mumford--Tate conjecture. Let $\mathcal{A}_g$ and $\mathscr{A}_g$ be the Siegel modular variety and its canonical integral model, with suitable level structures. The following conjecture is already implicit in an early work of Chai ({\cite[Conjecture 7.6, Conjecture 7.7]{Ch03}}):
\begin{conj}[Characteristic $p$ analogue of the Mumford--Tate conjecture for ordinary Shimura varieties of Hodge type]\label{conj:MTforAg}
Suppose $X_0$ is a geometric connected smooth 
 variety over $\mathbb{F}_q$ with a morphism $f_0$ into $\mathscr{A}_{g,\mathbb{F}_q}$, whose image lies in the ordinary locus. Let $A$ be the pullback abelian scheme over $X_0$. Define $\MTT(A)$ to be the generic Mumford--Tate group 
of a minimal special subvariety of $\mathscr{A}_{g,\overline{\bQ}}$ whose Zariski closure in $\mathscr{A}_{g,\overline{\bZ}}$ contains the image of $f_{0,\overline{\Fpbar}_p}$\footnote{Minimal special subvarieties having this property share the same generic Mumford--Tate group. Indeed, they are Hecke translates of each other. See \S\ref{subsub:smallest}.} Then the neutral component of the $l$-adic étale (\textit{resp}. crystalline) monodromy group of $A$ equals $\MTT(A)_{\bQ_l}$ (\textit{resp}. $\MTT(A)_{\bQ_p}$). In particular, $\MTT(A)$ is a rational model over $\bQ$ of the $l$-adic étale and crystalline monodromy groups of $A$. 
\end{conj}

Now we consider the characteristic $p$ analogue of the André--Oort conjecture. We call a subvariety of $\mathscr{A}_{g,\overline{\mathbb{F}}_p}$ \textbf{special}, if it is open dense in an irreducible component of the mod $p$ reduction of a classical special subvariety in characteristic 0. Naïvely, one can formulate the conjecture as follows: if a subvariety $X$ of the mod $p$ reduction of a Shimura variety contains a Zariski dense collection of special subvarieties, then $X$ is special. Unfortunately, this is not true: since every $\overline{\Fpbar}_p$-point in the mod $p$ reduction of a Shimura variety is already special, any positive dimensional subvariety contains a Zariski dense set of special points. 

To make it possibly true, one need to pose extra conditions on the collection of special subvarieties. A natural condition is that the special subvarieties in the collection are positive dimensional. However, even this is not enough to guarantee that $X$ is special, see Example~\ref{ex:ntla}. Nevertheless, with this condition, one can at least expect that $X$ ``admits a special factor''. In the following, let $\mathscr{A}_{g,\overline{\Fpbar}_p}^{\ord}$ 
be the ordinary stratum of $\mathscr{A}_{g,\overline{\Fpbar}_p}$.
\begin{conj}[Characteristic $p$ analogue of the André--Oort conjecture for ordinary Shimura varieties of Hodge type]\label{Conj:AOforSV}
Suppose $X\subseteq \mathscr{A}_{g,\overline{\Fpbar}_p}^{\ord}$ is a closed subvariety which contains a Zariski dense collection $\Xi$ of positive dimensional special subvarieties, then there is a special subvariety $Y\subseteq \mathscr{A}_{g,\overline{\Fpbar}_p}^{\ord}$ which splits as an almost product of two special subvarieties $Y_1$ and $Y_2$ with $\dim Y_1>0$, such that $X$ splits as an almost product of $Y_1$ with a possibly non-special subvariety $X'\subseteq Y_2$ \footnote{This will be termed as \textit{quasi-weakly special}, see Definition~\ref{def:quasi-weakly}.}. In particular, $X$ is an almost product of a positive dimensional special subvariety with a possibly non-special subvariety. 
\end{conj}
In some cases, one can show that $X$ is already special. For example, one of our main results (Theorem~\ref{Mainthm:AOSingle}) shows that if $X$ lies on a GSpin Shimura variety, then it is already special. However, if $X$ lies on a product of two or more Shimura varieties, it is in general not special, as the readers may have already seen in Example~\ref{ex:ntla}.

Suppose that Conjecture~\ref{Conj:AOforSV} holds, then the size of the possibly non-special factor $X'$ can also be well controlled:  Projecting $\Xi$ to $X'$ gives a collection $\Xi'$ of special subvarieties on $X'$. If $\Xi'$ contains a Zariski dense collection of positive dimensional special subvarieties, we can apply the conjecture to $X'$ to get a further splitting. Now use $X'$ to denote the new non-special factor thus appear, and iterate. Eventually we will arrive at a point where $X'$ is small enough so that $\Xi'$ does not contain a Zariski dense collection of positive dimensional special subvarieties. This optimizes the size of the non-special factor.

Besides ours, there are other attempts of formulating characteristic $p$ analogues of the André--Oort conjecture. In \cite[\S 7]{CO06}, Chai and Oort formulate an analogue using hypersymmetric points. In \cite{BaR20}, the authors define a ``mod $p$ special point'' to be the mod $p$ reduction of the whole Galois orbit of a CM point (see also \cite{Ro23}), and then prove a mod $p$ variant of André--Oort conjecture for a product of two modular curves, on condition of GRH. In general, it is relatively easy to formulate an analogue, but very hard to verify it in any special case. Our formulation of the conjecture is very different from the aforementioned works, and can be verified for products of arbitrary many GSpin Shimura varieties. Some aspects of our method are valid for general Shimura varieties, indicating that Conjecture~\ref{Conj:AOforSV} is very likely to hold in general.

\subsection{Main results I: the Mumford--Tate and the André--Oort conjectures}\label{sub:Mainresults}
We prove Conjecture~\ref{conj:MTforAg} and \ref{Conj:AOforSV} for products of GSpin Shimura varieties. In the following, let $\Fpbar=\overline{\Fpbar}_p$. A (locally closed) irreducible subvariety of $\mathscr{A}_{g,\Fpbar}$ is said to be \textbf{special}, if its Zariski closure is an irreducible component of
the reduction of the Zariski closure in $\mathscr{A}_g$ of a \textbf{special subvariety of Hodge type} in the sense of \cite[\S 1.1]{M98}. Loosely speaking, a subvariety is special if it is the reduction of a Shimura subvariety.


Let's briefly review what a GSpin Shimura datum is. Let $(L_{\bQ},Q)$ be a quadratic space with signature $(2,b)$. Associated to it is a GSpin Shimura datum $(\GSpin(L_{\bQ}),\mathcal{D}_L)$, with reflex field $\bQ$. It is known that $(\GSpin(L_{\bQ}),\mathcal{D}_L)$ admits a so called \textbf{Kuga--Satake embedding} into a Siegel Shimura datum. In particular, GSpin Shimura data are of Hodge type.


Now let $\cA_g$ be a Siegel modular variety with integral canonical model $\mathscr{A}_g$. Let $\bbI$ be a finite index set. For each $i\in \bbI$, let $\mathcal{S}_i$ be a GSpin Shimura variety.  Let $\mathcal{S}_\bbI$ be the product of various $\mathcal{S}_i$'s. Suppose that there is a Hodge embedding $\mathcal{S}_\bbI\hookrightarrow \mathcal{A}_g$. The \textbf{naïve integral model} of $\mathcal{S}_\bbI$ in $\mathscr{A}_g$ is defined to be the Zariski closure of $\mathcal{S}_\bbI$ in $\mathscr{A}_g$, and is denoted by $\mathscr{S}_\bbI$. In good cases, naïve integral models coincide with integral canonical models, by the work of Kisin (\cite{KM09}, see \S\ref{subsub:canonicalintegralmodelofHodge}). Let the ordinary locus $\mathscr{S}_{\bbI,\Fpbar}^{\ord}$ be $\mathscr{S}_{\bbI,\Fpbar
}\cap\mathscr{A}_{g,\Fpbar}^{\ord}$. In this paper, we will always assume that $\mathscr{S}_{\bbI,\Fpbar}^{\ord}\neq \emptyset$.

\begin{theorem}[$=$Theorem \ref{thm:MTTllinear2}. Characteristic $p$  Mumford--Tate conjecture for ordinary strata of products of GSpin Shimura varieties]\label{Mainthm:MTGSpin}
Conjecture~\ref{conj:MTforAg} is true if $f_0$ factors through the naïve integral model $\mathscr{S}_{\bbI,\bF_q}$. 
\end{theorem}

In particular, the theorem implies that the characteristic $p$ analogue of the Mumford--Tate conjecture holds for an arbitrary product of principally polarized ordinary abelian surfaces and K3 surfaces over a finitely generated field over $\mathbb{F}_p$. \\


The Mumford--Tate result~\ref{Mainthm:MTGSpin} is also one of the main ingredients for establishing the characteristic $p$ analogue of the André--Oort conjecture for products of GSpin Shimura varieties. 

Let's first state two special cases of our result on the characteristic $p$ André--Oort conjecture, which are less technical. The first is the single GSpin case (i.e., $\#\bbI=1$), which is the most akin to the classical story:  
\begin{theorem}[Characteristic $p$ André--Oort conjecture for ordinary GSpin Shimura varieties]\label{Mainthm:AOSingle} Suppose that $X$ is a closed subvariety of $\mathscr{A}_{g,\Fpbar}^{\ord}$ that contains a Zariski dense collection of positive dimensional special subvarieties. If $X$ is contained in the naïve integral model  $\mathscr{S}_{\Fpbar}$ of a GSpin Shimura subvariety, then $X$ is special. In particular, Conjecture~\ref{Conj:AOforSV} is true in this case.  
\end{theorem}

The product case is much more subtle, even for products of modular curves. As noted before, a subvariety containing a Zariski dense collection of positive special subvarieties may fail to be special. The author learned the following example from Chai:
\begin{example}\label{ex:ntla}Let $\bbI=\{1,2,3\}$. For each $i\in\bbI$, let $\mathcal{S}_i$ be a modular curve with integral canonical model $\mathscr{S}_i$. Consider the subvariety $X:=C\times \mathscr{S}_{3,\Fpbar}^{\ord}$ in the triple product space $ \mathscr{S}_{\bbI,\Fpbar}^{\ord}$, where $C$ is a non-special curve in $\mathscr{S}_{1,\Fpbar}^{\ord}\times \mathscr{S}_{2,\Fpbar}^{\ord}$. Then $X$ is not a special subvariety of $\mathscr{S}_{\bbI,\Fpbar}^{\ord}$. However, $X$ contains a Zariski dense collection of special curves $\{x\times \mathscr{S}_{3,\Fpbar}^{\ord}|x\in C(\Fpbar)\}$. 
\end{example}
However, Example~\ref{ex:ntla} is the only obstruction towards having a special subvariety. The following is another special case that is orthogonal to Theorem~\ref{Mainthm:AOSingle}:
\begin{theorem}[Characteristic $p$ André--Oort conjecture for products of ordinary modular curves]\label{Mainthm:AOmodularcurve}
  Let $\bbI$ be a finite index set. For each $i\in\bbI$, let $\mathcal{S}_i$ be a modular curve with integral canonical model $\mathscr{S}_i$. Let $X$ be a closed subvariety of $\mathscr{S}_{\bbI,\Fpbar}^{\ord}$ that contains a Zariski dense collection of positive dimensional special subvarieties. Let $\mathbf{I}_{S}\subseteq \bbI$ be the set of indices $i$ such that $X$ contains a Zariski dense collection of special subvarieties whose projections to $\mathscr{S}_{i,\Fpbar}$ are 
positive dimensional. Then $X$ is the product of a special subvariety of $\mathscr{S}_{\mathbf{I}_{S},\Fpbar}$ and a subvariety of $\mathscr{S}_{\mathbf{I}-\mathbf{I}_S,\Fpbar}$. In particular, Conjecture~\ref{Conj:AOforSV} is true. 
\end{theorem}
Now we are ready state the most general result: 

\begin{theorem}[$=$ Theorem~\ref{thm:AOprod}. Characteristic $p$ André--Oort conjecture for ordinary strata of products of GSpin Shimura varieties]\label{Mainthm:AOGSpin}
Suppose that $X$ is a closed subvariety of $\mathscr{A}_{g,\Fpbar}^{\ord}$ that contains a Zariski dense collection of positive dimensional special subvarieties. If $X$ is contained in the naïve integral model  $\mathscr{S}_{\bbI,\Fpbar}$ of a product of GSpin Shimura subvarieties, then Conjecture~\ref{Conj:AOforSV} holds.
\end{theorem}
It is possible to state a quantitative version of Theorem~\ref{Mainthm:AOGSpin} similar to Theorem~\ref{Mainthm:AOmodularcurve}. We choose not to do it here since it is a bit technical, but see \S\ref{subsub:AOsepfact} for more details. The deduction of Theorem~\ref{Mainthm:AOSingle} and \ref{Mainthm:AOmodularcurve} from Theorem~\ref{Mainthm:AOGSpin} will be given at the end of \S\ref{subsec:AOGSpin}. 

The reason why $X$ is already special in the single GSpin case~\ref{Mainthm:AOSingle} is that, a GSpin Shimura variety does not contain products of Shimura subvarieties of too big size. More precisley, if a GSpin Shimura variety admits a special subvariety that splits as an almost product of two smaller special subvarieties, then the two smaller special subvarieties are necessarily special curves. 
 
\subsection{Main results II: the Tate-linear and the Ax--Schanuel conjectures} Linearity is a fundamental concept which characterizes special subvarieties of a Shimura variety. It will play a crucial role in our treatment of the conjectures. Several different notions of linearity exist. We will give a brief review of them. For simplicity, in the following we consider linearity for subvarieties of $\cA_g$. This is already enough for dealing with Shimura subvarieties of Hodge type. 
\subsubsection{Linearity in char 0} Consider the uniformization map $\pi: \bH_g\rightarrow \cA_{g,\bC}$ with deck group an arithmetic subgroup of $\text{GSp}_{2g}(\bZ)$. One can make sense of algebraic subvarieties of $\bH_g$, cf. \cite[\S 3]{EY11}. A subvariety $V\subseteq \cA_{g,\bC}$ is called bi-algebraic if components of $\pi^{-1}(V)$ are algebraic. Since the morphism $\pi$ is highly transcendental, bi-algebraic subvarieties are rare, and have very special properties. In fact, being bi-algebraic puts a strong linear condition on $V$ (or $\pi^{-1}(V)$). This linearity can be understood better from more classical settings. Indeed, consider an abelian variety $A'$ over $\bC$ with a uniformization map $e:\bC^n\rightarrow A'$. It is classically known that an irreducible subvariety $V\subseteq A'$ is bi-algebraic if and only if $V$ is a translate of an abelian subvariety, and in this case, each component of $e^{-1}(V)$ is a linear subspace of $\bC^n$. In the case of Shimura varieties, similar phenomena happen. As a consequence of \cite{EY11}, a subvariety $V\subseteq \cA_{g,\bC}$ is bi-algebraic if and only if $V$ is a weakly special subvariety. In this case, each component of $\pi^{-1}(V)$ is a ``linear subspace'' of $\bH_g$ in the sense that it is totally geodesic. Indeed, a subvariety of a Euclidean space is linear if and only if it is totally geodesic. Therefore the property of being totally geodesic is a natural generalization of linearity in the classical sense.

The interpretation of linearity in terms of totally geodesic property can be found in much earlier works of Moonen. In \cite{M97}, Moonen showed that a subvariety of a Shimura variety is weakly special if and only if it is totally geodesic. In this and a subsequent paper \cite{M98}, Moonen also investigated the notion of ``formal linearity''. To state it, we recall that Serre--Tate theory implies that the completion of $\mathscr{A}_{g}$ at an ordinary mod $p$ point admits the structure of a formal torus, cf. \cite{K81}. It is usually called \textbf{Serre--Tate formal torus} or \textbf{Serre--Tate formal coordinates} or \textbf{canonical coordinates}. A subvariety of $\cA_g$ is called \textbf{formally linear} if the completion of its Zariski closure in $\mathscr{A}_{g}$ at an ordinary mod $p$ point is a union of torsion translates of subtori of the Serre--Tate formal torus. By an earlier result of Noot (\cite{N96}), special subvarieties of $\cA_g$ are formally linear. In \cite{M98}, Moonen showed the converse, i.e., formal linearity characterizes special subvarieties. 
\subsubsection{Linearity in char $p$} Linearity in terms of bi-algebraicity or totally geodesic property doesn't generalize easily to mod $p$ Shimura varieties. However, formal linearity does generalize to the ordinary Newton stratum: a subvariety $V\subseteq\mathscr{A}_{g,\Fpbar}^{\ord}$ is called formally linear at an $\Fpbar$-point $x$, if every irreducible component of $V^{/x}$ is a subtorus of the (mod $p$) Serre--Tate formal torus. If $V^{/x}$ is furthermore irreducible, then we adopt Chai's terminology, and call $V$ \textbf{Tate-linear} at $x$, cf. \cite[Definition 5.1]{Ch03}. A subvariety of $\mathscr{A}_{g,\Fpbar}^{\ord}$ is called \textbf{Tate-linear}, if it is {Tate-linear} at a point (hence at every smooth point, see \cite[Proposition 5.3, Remark 5.3.1]{Ch03}). 

Tate-linearity appears naturally in Chai's work on Hecke orbit conjecture, since the Zariski closure of an ordinary Hecke orbit is Tate-linear, see \cite{Chai06}. Note that Noot's result (\cite{N96}) implies that a special subvariety  of $\mathscr{A}_{g,\Fpbar}^{\ord}$ is Tate-linear at an ordinary point. It is natural to ask the converse: does Tate-linearity characterize special subvarieties?  This is very much Chai's Tate-linear conjecture: 
\begin{conj}[Tate-linear conjecture, see {\cite[Conjecture 7.2, Remark 7.2.1, Proposition 5.3]{Ch03}}]\label{conj:TTl} If an irreducible subvariety of $\mathscr{A}_{g,\Fpbar}^{\ord}$ is Tate-linear at a point, then it is special.
\end{conj}
The conjecture is still open. In this paper, we will also consider the following generalization:
\begin{conj}[Characteristic $p$ analogue of the Ax--Schanuel conjecture for ordinary Shimura varieties of Hodge type]\label{conj:AxSchanuel} Let $(X,x)$ be a pointed smooth connected variety over $\Fpbar$ that admits a morphism
$f$ into $\mathscr{A}_{g,\Fpbar}^{\ord}$. Let $\mathscr{T}_{f,x}\subseteq \mathscr{A}_{g,\Fpbar}^{/x}$ be the smallest formal subtorus through which $f^{/x}$ factors. Then $f$ factors through a special subvariety whose formal germ admits $\mathscr{T}_{f,x}$ as an irreducible component. 
\end{conj}
Conjecture~\ref{conj:TTl} is then a special case when $f$ is an immersion.

Conjecture~\ref{conj:AxSchanuel} can be thought of as an analogue of the Ax--Schanuel conjecture for Shimura varieties (\cite[Theorem 1.1]{T19}), in the following sense: Consider the uniformization map $\pi:\bH_g\rightarrow \mathcal{A}_{g,\bC}$ in the classical setting. A special case of the Ax--Schanuel conjecture asserts that if $V\subseteq \mathcal{A}_{g,\bC}$ is a subvariety, such that a component of $\pi^{-1}(V)$ is contained in a totally geodesic subvariety, then $V$ is contained in a weakly special subvariety. Conjecture~\ref{conj:AxSchanuel} is exactly a char $p$ analogue of this: \begin{itemize}
    \item taking completion is an analogue of taking uniformization, 
    \item the Serre--Tate formal torus $\mathscr{A}_{g,\Fpbar}^{/x}$ is an analogue of $\bH_g$, 
    \item a formal subtorus $\mathscr{T}_{f,x}\subseteq\mathscr{A}_{g,\Fpbar}^{/x}$ is an analogue of a totally geodesic subvariety of $\bH_g$,
    \item the condition ``$f^{/x}$ factoring through $\mathscr{T}_{f,x}$'' is an analogue of ``a component of $\pi^{-1}(V)$ being contained in a totally geodesic subvariety''.
\end{itemize}

One of our main results is that Conjecture~\ref{conj:AxSchanuel} is the consequence of the characteristic $p$ Mumford--Tate conjecture:
\begin{theorem}[$=$Proposition~\ref{prop:MTimpliesTl}]\label{thm:implications}
   Conjecture~\ref{conj:MTforAg} $\Rightarrow$ Conjecture~\ref{conj:AxSchanuel} $\Rightarrow$ Conjecture~\ref{conj:TTl}.
\end{theorem} 
It follows from Theorem~\ref{Mainthm:MTGSpin} that
\begin{theorem}[$=$ Theorem~\ref{thm:MTTllinear}]\label{MainThm:genTlGSpin} Conjecture~\ref{conj:TTl} and Conjecture~\ref{conj:AxSchanuel} are true if $f$ factors through the naïve integral model a product of GSpin Shimura varieties. 
\end{theorem}
However, in the paper we take the opposite direction: we first establish Theorem~\ref{MainThm:genTlGSpin} and then use it as an ingredient to establish Theorem~\ref{Mainthm:MTGSpin}.

\subsection{Proof strategies for Theorem~\ref{Mainthm:MTGSpin} and Theorem~\ref{MainThm:genTlGSpin}}
For simplicity, we will only discuss the proof strategies for single GSpin Shimura varieties. We briefly summarize the ingredients: \begin{itemize}
    \item  Ingredients for the mod $p$ Ax--Schanuel conjecture (Theorem~\ref{MainThm:genTlGSpin}):\begin{enumerate}
        \item The parabolicity conjecture (\cite{MD20}),
        \item Chai's local monodromy theorem (\cite[\S 3-4]{Ch03}, see also Theorem~\ref{thm:etale-glob}).
    \end{enumerate} 
    \item Ingredients for the mod $p$ Mumford--Tate conjecture (Theorem~\ref{Mainthm:MTGSpin}):
    \begin{enumerate}
        \item The mod $p$ Ax--Schanuel conjecture (Theorem~\ref{MainThm:genTlGSpin}). 
        \item Monodromy of coefficient objects in compatible families (\cite[Theorem 1.2.1]{DA20}).
    \end{enumerate}
\end{itemize}
Here is a sketch of our methods:  

Let $(X,x)$ be a pointed smooth connected variety over $\Fpbar$ with a morphism $f$ into $\mathscr{A}_{g,\Fpbar}^{\ord}$ that factors through the naïve integral model $\mathscr{S}$ of a GSpin Shimura subvariety $\mathcal{S}\hookrightarrow \mathcal{A}_g$. First,  we reduce to the case where $\mathcal{S}\hookrightarrow \mathcal{A}_g$ is a Kuga--Satake embedding, and the naïve integral model $\mathscr{S}$ coincides with the integral canonical model. General strategies for doing such reductions will be worked out in \S\ref{sub:inpmM} and \S\ref{subsub:prodsetup}. 

The scheme $\mathscr{S}$ is equipped with a so-called Kuga--Satake abelian scheme, which is nothing other than the pullback of the universal family from $\mathscr{A}_g$.  Let $\tilde{x}_\bC$ be the canonical lift of $x$, base changed to $\bC$. There is a smallest special subvariety $\mathcal{X}_f\subseteq \mathcal{S}$ containing $\tilde{x}_\bC$, whose Zariski closure $\mathscr{X}_f\subseteq\mathscr{S}_{\mathbb{F}}$ contains the image of $f$. We let $\MTT(f)$ be the generic Mumford--Tate group of $\mathcal{X}_f$ with respect to the Kuga--Satake abelian scheme. Let $\mathscr{T}_{f,x}$ be the smallest formal subtorus of $\mathscr{S}_{\Fpbar}^{/x}$ through which $f^{/x}$ factors. Let $f_0:X_0\rightarrow \mathscr{S}_{\Fpbar_q}$ be a finite field model of $f$. 

To prove Theorem~\ref{MainThm:genTlGSpin}, it suffices to show that $\mathscr{X}_{f,\Fpbar}^{/x}$ contains $\mathscr{T}_{f,x}$ as an irreducible component. It is classically known that  $\mathscr{X}_{f,\Fpbar}^{/x}$ is a union of formal tori (\cite{N96}), so it already contains $\mathscr{T}_{f,x}$. Therefore it suffices to show that $\mathscr{X}_{f,\Fpbar}^{/x}$ has small dimension, i.e., $\mathcal{X}_{f}$ is cut out by enough Hodge cycles. To construct Hodge cycles, our strategy is to first understand the monodromy group $G_p^-$ of the pullback $p$-divisible group $A_{X_0}[p^\infty]$. The parabolicity conjecture (\cite{MD20}) predicts that $G_p^-$ is the parabolic subgroup of the crystalline monodromy group $G_p$ of $A_{X_0}$ that fixes the slope filtration. On the other hand, Chai's local monodromy theorem (\cite[\S 3-4]{Ch03}) implies that the unipotent radical of $G_p^-$ is determined by $\mathscr{T}_{f,x}$. As a result, we can largely understand the structure of $G^-_p$, which enables us to construct enough endomorphisms of $A$ over an étale cover of $X_0$ (using the crystalline Tate conjecture (\cite{DJ98})). We then canonical lift these endomorphisms to $\tilde{x}_\bC$, and we use deformation theory to show that these endomorphisms cut out $\mathcal{X}_f$. This will imply that $\mathscr{T}_{f,x}$ is an irreducible component of $\mathscr{X}_{f,\Fpbar}^{/x}$. Note that we only need to construct endomorphisms of $A$, but not higher Hodge cycles. This is a special property of (products of) GSpin Shimura varieties: special subvarities are cut out by endomorphisms.  

Theorem~\ref{MainThm:genTlGSpin} then implies that $G^\circ_p$ and
$\MTT(f)_{\bQ_p}$ admit a common factor (see \S\ref{subsec:MTgpclass}). We will build several compatible systems from $\MTT(f)$, and use Chebotarev's density theorem and \cite[Theorem 1.2.1]{DA20} to conclude that $\MTT(f)_{\bQ_p}\simeq G_p^\circ$. This establishes  Theorem~\ref{Mainthm:MTGSpin} in the single GSpin case. 

\subsubsection{Comparison with the Mumford--Tate conjecture in characteristic 0} A very effective tool for studying the Mumford--Tate conjecture in the classical number field setting is the classification of Mumford--Tate pairs introduced by Serre (\cite{Serre1}). The method was developed in various works (\cites{Chi,Ta1,Ta2,Pre,V08}). For example, suppose that there is a reductive group $G$ with a faithful representation $\rho:G\rightarrow \GL(V)$, such that $G$ is generated by conjugates of cocharacters of weights lying in $\{0,1\}$ (this is called a \textbf{weak Mumford--Tate pair of weights $\{0,1\}$}). Then the possibility of $(G,\rho)$ is very restricted. In the Mumford--Tate setting, one takes $G$ to be either the Mumford--Tate group or the $l$-adic monodromy group, and takes the cocharacters to be the Hodge cocharacters of the reductions of an abelian variety over different places. The method can be briefly summarized as ``\textbf{studying the Mumford--Tate problem by cocharacters}".   

It is possible to study the characteristic $p$ Mumford--Tate conjecture by cocharacters, at least over the ordinary stratum. Indeed, after we set up the general framework in \S\ref{sec:2}, one sees that, as long as the image of $f$ lies in the ordinary stratum, the $l$-adic monodromy still form a weak Mumford--Tate pair of weights $\{0,1\}$. 

However, our method is orthogonal to the classical one. It does not rely on the classification of Mumford--Tate pairs. Instead, it crucially uses the Tate-linear structure in characteristic $p$, and can be summarized as ``\textbf{studying the Mumford--Tate problem by unipotents}''. In the characteristic $p$ world, it has at least two advantages: \begin{enumerate}
    \item 
It admits generalization to higher Newton strata, where the $l$-adic monodromy groups may no longer form Mumford--Tate pairs of weights $\{0,1\}$. 
\item It has the potential to go beyond the group theoretic limit of the method of Mumford--Tate pairs. 
Recall that our proof of the mod $p$ Mumford--Tate Conjecture~\ref{conj:MTforAg} relies on the mod $p$ Ax--Schanuel Conjecture~\ref{conj:AxSchanuel}, which is of an algebraization/unlikely intersection flavor. In general, one may hope that the mod $p$ Mumford--Tate problem can be approached by techniques from the areas of algebraization and unlikely intersections. 
\end{enumerate}
These direction are being explored in an ongoing work, which will come out soon after this paper.
\subsection{Proof strategies for Theorem~\ref{Mainthm:AOGSpin}}
We summarize the ingredients in the proof of the mod $p$ André--Oort conjecture (Theorem~\ref{Mainthm:AOGSpin}): \begin{enumerate}
    \item The mod $p$ Mumford--Tate conjecture (Theorem~\ref{Mainthm:MTGSpin}),
    \item The mod $p$ Ax--Schanuel conjecture (Theorem~\ref{MainThm:genTlGSpin}),
    \item  A generalization of Chai's rigidity theorem (\cite{Chai08}) to geometric bases (Theorem~\ref{lm:rigidity}), 
    \item Global canonical coordinates (\cite[\S2]{Ch03}).
\end{enumerate}
We begin by a baby example ({which the the author learned from Ananth Shankar and Yunqing Tang}): 
\begin{example}\label{ex:babyex}
   Let $\mathscr{S}$ be the Siegel modular variety $\mathscr{A}_2$ (which is also a GSpin Shimura variety). This example shows that if a divisor $X\subseteq \mathscr{S}_{\Fpbar}$ contains Zariski dense positive dimensional special subvarieties that pass through the same ordinary $\Fpbar$-point $x$ of $X$, then $X$ is special. In fact, each special subvariety gives rise to a (union of) formal subtorus of $\mathscr{S}_{\Fpbar}^{/x}$ which is also contained in $X^{/x}$. One can show that these formal subtori are Zariski dense in $X^{/x}$. Since each subtorus is invariant under scaling by $\bZ_p^*$, $X^{/x}$ is also invariant under scaling by $\bZ_p^*$. By Chai's rigidity theorem, $X^{/x}$ is a formal subtorus of $\mathscr{S}_{\Fpbar}^{/x}$. Theorem~\ref{MainThm:genTlGSpin} then implies that $X$ is special.
\end{example}
In general, the Zariski dense collection of positive dimensional special subvarieties of $X$ may not pass through the same point. Even if they do, the collection of formal subtori thus arise may not be Zariski dense in $X^{/x}$. Hence 
the strategy in Example~\ref{ex:babyex} won't work. Instead, we use the Zariski dense collection of special subvarieties to construct certain arithmetic $\bZ_p$-lisse sheaves on (an étale open subset of) $X$. The construction can be summarized as follows: 
\begin{construction}\label{constr:2}
    Consider $(X\times \mathscr{S}_{\Fpbar})^{/\Delta}$, where $\Delta$ is the graph of the immersion $X\subseteq \mathscr{S}_{\Fpbar}$. For each special subvariety $Z\subseteq X$, we consider $(Z\times Z)^{/\Delta}$, where $\Delta$ stands for the diagonal. Let $\mathfrak{Z}$ be the smallest closed formal subscheme of $(X\times \mathscr{S}_{\Fpbar})^{/\Delta}$ that contains all $(Z\times Z)^{/\Delta}$. The theory of global canonical coordinates implies that $(X\times \mathscr{S}_{\Fpbar})^{/\Delta}$ is a family of formal tori over $X$, and $(Z\times Z)^{/\Delta}$ is a family of formal tori over $Z$. The scaling-by-$\bZ^*_p$ map on $(X\times \mathscr{S}_{\Fpbar})^{/\Delta}$ preserves each $(Z\times Z)^{/\Delta}$, hence $\mathfrak{Z}$. Let $\eta$ be the generic point of $X$, then Chai's rigidity theorem implies that $\mathfrak{Z}\widehat{\times}_{X} \overline{\eta}$ is a union of formal tori. Each formal torus then gives rise to a  $\bZ_p$-lisse sheaf $\mathscr{F}$ over an étale open subset of $X$. It is \textit{essentially} the lisse sheaf that we want\footnote{We are hiding away some technical difficulty. The actual construction of $\mathscr{F}$ is much harder than what is described here, and requires a generalization of Chai's rigidity theorem to geometric bases, see Theorem~\ref{lm:rigidity} and Appendix~\ref{subsub:Starshape}} .
\end{construction}  
For simplicity, let's assume that $\mathscr{F}$ is defined over $X$. Note that $\mathfrak{Z}\subseteq (X\times X)^{/\Delta}$. One of the main features of our construction is that, for Zariski dense $x\in X(\Fpbar)$, the deformation space $X^{/x}$ is sandwiched between two formal tori:\begin{equation}\label{eq:swswsw}
    \mathscr{F}_x\otimes_{\bZ_p} \bG_m^{\wedge}\subseteq X^{/x}\subseteq \mathscr{T}_x,
\end{equation}
where $\mathscr{T}_{x}$ is the smallest formal torus containing $X^{/x}$. Note that the lisse sheaf $\mathscr{F}$ can be understood via representation theory. In the single GSpin case, it follows from the Mumford--Tate result (Theorem~\ref{Mainthm:MTGSpin}) that all inclusions in (\ref{eq:swswsw}) are equalities. This implies that $X^{/x}$ is a formal torus. Theorem~\ref{Mainthm:AOSingle} then follows from Theorem~\ref{MainThm:genTlGSpin}. 

We remark that Construction~\ref{constr:2} works not only for GSpin Shimura varieties, but general Shimura varieties of Hodge type as well. The construction always produces a sandwich of form (\ref{eq:swswsw}), reducing the question of understanding the structure of $X^{/x}$ to a question of representation theoretical nature, to which one can apply the Mumford--Tate type results. We refer the readers to the \textbf{geometric squeeze theorem}~\ref{lm:thekeylemma} for more details.  

For products of GSpin Shimura varieties, the representation theory is different, so one cannot deduce that (\ref{eq:swswsw}) is an equality as in the single GSpin case. Nevertheless, the representation theory will guarantee that $X^{/x}$ has ``a positive dimensional toric factor'', which corresponds to the fact that $X$ has a positive dimensional special factor. The following is a concrete example for the construction of the $p$-adic lisse sheaf for a triple product of modular curves:
\begin{example}
Let $X$ be as in Example~\ref{ex:ntla}. We replace $X$ by its ordinary stratum. Consider the projection $\pi_3: X\rightarrow\mathscr{S}_{3,\Fpbar}^{\ord}$. Then $\mathfrak{Z}$ in Construction~\ref{constr:2} is nothing other than the pullback via $\pi_3$ of  $(\mathscr{S}^{\ord}_{3,\Fpbar}\times \mathscr{S}_{3,\Fpbar}^{\ord})^{/\Delta}$. It is a family of rank 1 formal tori over $X$. The lisse sheaf $\mathscr{F}$ thus arise is the pullback via $\pi_3$ of the $\bZ_p$-coefficient relative étale cohomology of the universal family over $\mathscr{S}^{\ord}_{3,\Fpbar}$.
\end{example}

\subsection{Organization of the paper} In \S\ref{sec:1} we review the notion of GSpin Shimura varieties and the theory of canonical coordinates. We also introduce the notion of quasi-weakly special subvarieties. In \S\ref{sec:monodromy} we review previous works on monodromy of $F$-isocrystals.  In \S\ref{sec:2} we construct $\MTT(f)$ and $\mathcal{X}_f$, and study their properties in a general setting. In \S\ref{Sec:ttlinearsingle}  we prove Theorem~\ref{Mainthm:MTGSpin} and Theorem~\ref{MainThm:genTlGSpin}. 
In \S\ref{sec:AO} we prove Theorem~\ref{Mainthm:AOGSpin}. 

\subsection{Notations and conventions}\label{notation}
We use $p$ to denote an odd prime and $q$ to denote a positive power of $p$. We fix an algebraic closure of $\bF_p$, and denote it by $\bF$. We denote by $W$ the ring of Witt vectors of $\bF$ and by $K$ its fraction field. Let $\overline{K}$ be an algebraic closure of $K$ and let $\overline{W}$ be the integral closure of $W$ in $\overline{K}$.  
We also make the following conventions: 
\begin{itemize}
    \item (Algebraic closures)
We fix once and for all an identification of ${\bC}_p:=\widehat{\overline{\bQ}}_p$ with $\bC$, as well as an embedding of $\overline{K}$ into $\bC_p$. As a result, we have fixed an embedding for any finite extension of $W$ into $\bC$.
\item(Tori and formal tori) We use $\bG_m$ \textit{resp}. $\bG_m^{\wedge}$ to denote the rank 1 torus \textit{resp}. formal torus over $W$, $\bZ_p$ or $\Fpbar$,  depending on the context. Sometimes we will also write $\bG_{m,W}$, $\bG_{m,\bZ_p}$ and $\bG_{m,\Fpbar}$ to emphasize the base scheme.
\item(Formal completions) Suppose $\mathscr{Y}$ is a $W$-scheme and $y\in \mathscr{Y}(\Fpbar)$. We write
${\mathscr{Y}}^{/y}$ for the completion of $\mathscr{Y}$ along $y$. If $Y$ is a $\Fpbar$-scheme and $Z$ is a closed subscheme, we write $Y^{/Z}$ for the completion of $Y$ along $Z$. 
If $f:X\rightarrow Y$ is a morphism of varieties over $\Fpbar$ and $x\in X(\Fpbar)$, we write $f^{/x}:X^{/x}\rightarrow Y^{/x}$ for the completion of $f$ at $x$. Note that $Y^{/x}$ actually stands for $Y^{/f(x)}$.
\item(Irreducibility) An irreducible component of a scheme is also required to be reduced.
\end{itemize}

\subsection*{Acknowledgements} The author thanks Ananth Shankar for encouraging him to think about the possible applications of the Tate-linear conjecture to characteristic $p$ analogues of the Mumford--Tate and André--Oort conjectures, together with all of his help and enlightening conversations when the author is writing the paper. The author thanks Ching-Li Chai for his previous works and theoretical buildups, without which the results of this paper won't be possible. The author also thanks Dima Arinkin, Ching-Li Chai, Asvin G, Qiao He, Jiaqi Hou, Brian Lawrence, Yu Luo, Keerthi Madapusi Pera, Devesh Maulik, Yunqing Tang, Yifan Wei, and Ziquan Yang for valuable discussions. Finally, the author thanks Marco D'Addezio for helpful comments as well as kindly pointing out some mistakes in an earlier draft of the paper. The author is partially supported by the NSF grant DMS-2100436. 
\section{Preliminaries}\label{sec:1}
In this section, we review the notion of special subvarieties, GSpin Shimura varieties and their arithmetic deformation theory. In \S\ref{sub:Quasi-weakly}, we also study the ordinary integral models of special subvarieties and introduce the notion of quasi-weakly special subvarieties.  
\subsection{Shimura varieties of Hodge type and special subvarieties}\label{Sec:GSpinDef}
Let $H$ be a self-dual symplectic $\bZ_{(p)}$-lattice. One can form a Siegel Shimura datum $(\GSp(H_\bQ),\mathfrak{H}^{\pm})$ with reflex field $\bQ$. Let $\bK=\bK_p{\bK}^p$ be a hyperspecial level structure such that $\bK_p=\GSp(H)(\bZ_p)$, ${\bK}^p\subseteq \GSp(H_{\bQ})(\bA_f^p)$ compact open, and $\bK$ leaves $H_{\widehat{\bZ}}$ stable. For $\bK'$ sufficiently small, we get a smooth variety $\cA_{g,\bK}:=\Sh_{\bK}(\GSp(H_\bQ),\mathfrak{H}^{\pm})$ over $\bQ$, called the Siegel modular variety, which admits an integral canonical model $\mathscr{A}_{g,\bK}$ over $\bZ_{(p)}$ (\cite[Theorem 2.10]{Mil92}). Integral models of this sort will be the ambient spaces that we consider in this paper.

The proof of the existence of the integral canonical model relies on the interpretation of $\cA_{g,\bK}$ as a moduli space of abelian varieties. In fact, when $\bK$ is sufficiently small, then $H$ gives rise to a universal abelian scheme over $\cA_{g,\bK}$ (up to prime-to-$p$ isogeny), whose first $\bZ_{(p)}$-Betti cohomology is the local system induced by $H$. The abelian scheme can be extended to $\mathscr{A}_{g,\bK}$, and we denote by $A\rightarrow\mathscr{A}_{g,\bK}$ the universal abelian scheme over  $\mathscr{A}_{g,\bK}$. Let $\bH_\text{B},\bH_{\text{dR}}, \bH_{l,\text{ét}}$ be the rational Betti, de Rham, $l$-adic étale ($l\neq p$) relative first cohomology of $A\rightarrow \mathscr{A}_{g,\bK}$, and let $\bH_{\text{cris}}$ be the first rational crystalline cohomology of $A_{\mathbb{F}_p}\rightarrow \mathscr{A}_{g,\bK,\Fpbar_p}$. 
\subsubsection{Shimura varieties of Hodge type and integral canonical models}\label{subsub:canonicalintegralmodelofHodge} 
A Shimura datum $(G,\mathcal{D})$ is of \textbf{Hodge type}, if there is a Siegel Shimura datum $(\GSp(H_\bQ),\mathfrak{H}^{\pm})$, and an embedding of Shimura data $\iota:(G,\mathcal{D})\hookrightarrow (\GSp(H_\bQ),\mathfrak{H}^{\pm})$. Let $E$ be the reflex field of $(G,\mathcal{D})$. For a suitable level structure $\bbK\subseteq G(\bA_f)$ such that $\iota(\bbK)\subseteq \bK$, we get a morphism of Shimura varieties. 
\begin{equation}\label{eq:shimuraHodgemorphi}
\Sh_{\iota,\bbK,\bK}:\Sh_{\bbK}(G,\mathcal{D})\rightarrow \cA_{g,\bK,E}
\end{equation}
which is an embedding for suitable choice of $\bbK$ and $\bK$. 

Suppose that there is a reductive group $G_{\bZ_{(p)}}$ over $\bZ_{(p)}$ with generic fiber $G$, and let $\bbK_p=G_{\bZ_{(p)}}(\bZ_p)$ be the associated hyperspecial subgroup. Following \cite[\S1.3.3]{MKmodp}, up to embedding in a larger Siegel Shimura datum (using Zarhin's trick), one can assume that $\iota$ is induced by an embedding $G_{\bZ_{(p)}}\hookrightarrow \GL(H)$, where $H$ is a self-dual symplectic $\bZ_{(p)}$-lattice. Let $\bK_p=\GSp(H)(\bZ_p)$, and let $\mathfrak{p}$ be a prime of $O_E$ lying above $p$.

From \cite[\S 2.3.2]{KM09}, for sufficiently small $\bbK^p\subseteq G(\bA_f^p)$, there is a sufficiently small $\bK^p\subseteq \GSp(H_\bQ)(\bA_f^p)$, so that there is an embedding $\mathcal{S}_\bbK:=\Sh_\bbK(G,\mathcal{D})\hookrightarrow \cA_{g,\bK,E}$, such that the normalization of the Zariski closure of $\mathcal{S}_\bbK$ in $\mathscr{A}_{g,\bK,O_{E,(\mathfrak{p})}}$ is the integral canonical model of $\mathcal{S}_\bbK$. We will denote it by $\mathscr{S}_\bbK$. In \cite{XU20}, it is shown that the normalization is redundant. 

We can pullback the universal abelian scheme and local systems from $\mathscr{A}_{g,\bK,O_{E,(\mathfrak{p})}}$ to $\mathscr{S}_\bbK$. These will again be called the universal family and local systems over $\mathscr{S}_\bbK$ with respect to the embedding. 

\subsubsection{Special subvarieties and naïve integral models}\label{subsub:naivespecial} Let's continue assuming that $\bK$ is a sufficiently small hyperspecial level structure. Roughly speaking, a \textbf{special subvariety of $\cA_{g,\bK}$ of Hodge type} is (a component of) the Hecke translate of the image of a map of the form  (\ref{eq:shimuraHodgemorphi}). We refer the readers to \cite[Definition 2.5]{M97}, \cite[\S1]{M98} for precise definitions. 

Let $k$ be either a number field or $\overline{\bQ}$, and let $\mathfrak{p}$ be a prime of $O_k$ that lies above $p$. Let $\mathcal{X}\subseteq \cA_{g,\bK,k}$ be a special subvariety. The \textbf{naïve integral model} of $\mathcal{X}$ is the Zariski closure of $\mathcal{X}$ inside $\mathscr{A}_{g,\bK,O_{k,(\mathfrak{p})}}$, denoted by $\mathscr{X}$. To eliminate the ambiguity, we will always consider $k$ as a subfield of $\bC$, and require the prime $\mathfrak{p}$ to be induced from the fixed isomorphism $\bC\simeq \bC_p$ as in \S\ref{notation}.

Let $\widetilde{\mathscr{X}}$ be the normalization of $\mathscr{X}$. In general $\mathscr{X}$ (even $\widetilde{\mathscr{X}}$) can be very singular, and depends on the embedding. Unfortunately, since we are dealing with random special subvarieties of $\cA_{g,\bK}$ in this paper, we have no choice but working with naïve integral models.  
 
There is a notion of universal abelian scheme and various local systems over $\mathscr{X}$ with respect to the embedding. They are just the pullback of the corresponding objects from $\mathscr{A}_{g,\bK}$. 

\subsection{GSpin Shimura varieties}\label{subsub:GSpindef} GSpin Shimura varieties are typical examples of Shimura varieties of Hodge type. We review the definitions and basic properties of GSpin Shimura varieties, following \cite{KR00},\cite{MP15} and \cite{AGHMP17}. 

\subsubsection{Self-dual case}\label{subsub:Selfdual} For an integer $b\geq 1$, let $(L,Q)$ be a self-dual quadratic $\mathbb{Z}_{(p)}$-lattice of rank $b+2$ and signature $(2,b)$. Let $\GSpin(L,Q)\subseteq \Cl(L)$ be the group of spinor similitude, where $\Cl(\cdot)$ is the Clifford algebra. The group $\GSpin(L_\bR)$ acts on the symmetric space 
$\mathcal{D}_L=\{z\in\mathbb{P}(L_\mathbb{C})|(z,z)=0,(z,\overline{z})<0\}$
via $c:\GSpin(L_\bQ)\rightarrow \SO(L_\mathbb{Q})$. This gives rise to a Shimura datum $(\GSpin(L_\bQ),\mathcal{D}_L)$ with reflex field $\mathbb{Q}$. In the following, let $\bbK=\bbK_p\bbK^p$ 
be a sufficiently small hyperspecial level struture, with $\bbK_p= \GSpin(L)(\bZ_p)$, $\bbK^p\subseteq \GSpin(L_\bQ)(\bA_f^p)$ compact open, and $\bbK\subseteq \Cl(L_{\widehat{\mathbb{Z}}})^\times$. We then get a smooth variety $\mathcal{S}_\bbK:=\mathcal{S}h_\bbK(\GSpin(L_\bQ),\mathcal{D}_L)$ over $\mathbb{Q}$, called the GSpin Shimura variety, which admits a canonical smooth integral model $\mathscr{S}_\bbK$ over $\mathbb{Z}_{(p)}$ (\cite[Theorem 2.3.8]{KM09}). When the level structure is clear from the context, we will simply drop it from the subscript.

Let $H=\Cl(L)$ with the action of itself on the right. Equip $\Cl(L)$ with the action of $\GSpin(L)$ on the left. Then there exist (1) a self-dual symplectic structure on $H$, which gives rise to an embedding of Shimura data $(\GSpin(L_\bQ),\mathcal{D}_L)\rightarrow (\GSp(H_\mathbb{Q}),\mathfrak{H}^{\pm})$, called the \textbf{Kuga--Satake embedding}, and (2) a sufficiently small hyperspecial level structure $\bK\subseteq \GSp(H_\mathbb{Q})(\bA_f)$ fixing $H_{\widehat{\bZ}}$, with $\bK_p= \GSp(H_{\mathbb{Z}_p})(\bZ_p)$, so that there is an embedding of integral canonical models $\mathscr{S}_\bbK\hookrightarrow \mathscr{A}_{g,\bK}$. In other words, $\mathscr{S}_\bK$ is the naïve integral model of the special subvariety $\mathcal{S}_\bbK$ with respect to its embedding $\mathcal{S}_\bbK\hookrightarrow\mathscr{A}_{g,\bK}$.

The universal abelian scheme $A\rightarrow\mathscr{A}_{g,\bK}$ induced by $H$ admits a left $\Cl(L)$-action. The pullback of $A$ to $\mathscr{S}_\bbK$ is called the Kuga--Satake abelian scheme. By abuse of notation, we will again use $A$ to denote it. We will again use  $\bH_\bullet$ to denote the pullback sheaves on $\mathscr{S}_\bbK$.

The natural action of $L$ on $H$ produces a $\GSpin(L)$ invariant embedding $L\hookrightarrow \End_{\Cl(L)}(H)$. Correspondingly,  for each $ \bullet=$  $\text{B},\text{dR},\{l,\text{ét}\}$ and cris, there is a local system $\mathbb{L}_{\bullet}$. The local system is equipped with a natural quadratic form $\mathbf{Q}$ such that $f\comp f=\mathbf{Q}(f)\Id$ for a section $f$ of $\mathbb{L}_\bullet$. It also admits an embedding $\mathbb{L}_\bullet\hookrightarrow \mathcal{E}nd_{\Cl(L)}(\mathbb{H}_{\bullet})$, which is compatible with various $p$-adic Hodge theoretic comparisons, see \cite[\S 4.3]{AGHMP17} and \cite[Proposition 3.11, 3.12, 4.7]{MP15}. 

\begin{remark}
 Low dimensional GSpin Shimura varieties recover some PEL type Shimura varieties. For example,  modular curve, Hilbert modular surfaces and the Siegel modular variety $\mathcal{A}_2$ are special cases of GSpin Shimura varieties when $b=1,2,3$, respectively. 
\end{remark}
\begin{remark}\label{rmk:embwedge}
 A natural isomorphism exists between the $\bQ$-vector spaces $\wedge L_\bQ$ and $\Cl(L_\bQ)$. This can be realized by selecting an orthogonal basis $\{e_1,...,e_{b+2}\}$ of $L_{\bQ}$, and identifying $e_{i_1}\wedge e_{i_2}\wedge...\wedge e_{i_k}$ with $e_{i_1}e_{i_2} ... e_{i_k}$ (it is independent of the choice of the orthogonal basis). The natural action of $\wedge L_{\bQ}$ on $H_{\bQ}$ results in a $\GSpin(L_{\bQ})$ invariant embedding $\wedge L_{\bQ}\hookrightarrow \End_{\Cl(L_{\bQ})}(H_{\bQ})$. As a consequence, there is a natural embedding of local systems $\wedge \mathbb{L}_\bullet\hookrightarrow \mathcal{E}nd_{\Cl(L)}(\mathbb{H}_\bullet)$. 
\end{remark}
\subsubsection{Non-self-dual case}\label{subsub:Nonselfdual} Suppose that $(L,Q)$ is not self dual. We can still define a Shimura datum $(\GSpin(L_\bQ),\mathcal{D}_{L})$ and define a Kuga--Satake embedding into a Siegel Shimura datum, exactly in the same manner as described above. However, these constructions may not behave well integrally. 
   
It is explained in \cite[Lemma 6.8]{MP15} that one can always find an isometric embedding of $(L,Q)$ into a self-dual $\bZ_{(p)}$-lattice $(\widetilde{L},\widetilde{Q})$, such that $\widetilde{L}=L\oplus L^{\perp}$, where ${L}^{\perp}$ is negative definite. We can then realize $\Sh_{\bbK}(\GSpin(L_\bQ),\mathcal{D}_{L})$ as a special cycle of $\Sh_{\widetilde{\bbK}}(\GSpin(\widetilde{L}_\bQ),\mathcal{D}_{\widetilde{L}})$. One can construct a naïve integral model of $\mathcal{S}_\bbK$ by taking Zariski closure inside the integral canonical model of $\Sh_{\widetilde{\bbK}}(\GSpin(\widetilde{L}_\bQ),\mathcal{D}_{\widetilde{L}})$.

\subsection{Canonical coordinates}\label{subsub:Cancoorpdiv} In this section, we briefly review several theories of canonical coordinates. The most classical theory concerns the deformation of an ordinary $p$-divisible group over $\Fpbar$ (see  \S\ref{sub:ST}). Later on Chai generalized it to a family of ordinary $p$-divisible group over more general bases (see \S\ref{sub:globalST}). These results can be used to describe the formal structures of naïve integral models of Shimura subvarieties at ordinary points (see \S\ref{subsub:canspecialsub} and \S\ref{subsub:pHodge}).

\subsubsection{The canonical Hodge cocharacter}\label{subsub:canHodge}
Let $\mathscr{G}$ be an ordinary $p$-divisible group over $\Fpbar$, and let $\mathscr{G}^{\et}$, $\mathscr{G}^{\loc}$ be its étale and local part. Let $F$ be the Frobenius on the Dieudonné module $\bD(\mathscr{G})$. Define a $\bZ_p$-module  
\begin{equation}
 \label{eq:fiberfunomega}
 \omega(\mathscr{G})=\{v\in \bD(\mathscr{G})|Fv=v \text{ or }pv\}.
 \end{equation}
There are canonical identifications \begin{equation}\label{eq:XTp}
    \omega(\mathscr{G}^{\loc})\simeq X^*(\mathscr{G}^{\loc}),\,\,\, \omega(\mathscr{G}^{\et})\simeq T_p(\mathscr{G}^{\et})^{\vee},
\end{equation}
where the symbols $X^*$ and $T_p$ stand for the character lattice and the $p$-adic Tate module. Note that $$\omega(\mathscr{G})=\omega(\mathscr{G}^{\loc})\oplus \omega(\mathscr{G}^{\et}).$$ The $\bZ_p$-module $\omega(\mathscr{G})$ is a \textbf{canonical} $\bZ_p$-structure of $\bD(\mathscr{G})$, in the sense that  $\bD(\mathscr{G})=\omega(\mathscr{G})\otimes W$ and 
\begin{equation}\label{eq:recoveer }F=\begin{bmatrix}
  \Id_{\omega(\mathscr{G}^{\et})} &0\\
  0&p\cdot\Id_{\omega(\mathscr{G}^{\loc})}  
\end{bmatrix}\sigma. 
\end{equation}
Define the \textbf{canonical Hodge cocharacter} as
\begin{equation}\label{eq:Newtoncochar}
    \mu:\bG_{m,\bZ_p}\rightarrow \GL(\omega(\mathscr{G})),\,\,\,t\rightarrow \begin{bmatrix}
  \Id_{\omega(\mathscr{G}^{\et})} &0\\
  0&t\cdot\Id_{\omega(\mathscr{G}^{\loc})}  
\end{bmatrix}.
\end{equation}
The inverse cocharacter $\mu^{-1}$ induces the slope filtration on $\bD(\mathscr{G})$.

\subsubsection{Canonical coordinates for ordinary $p$-divisible groups}\label{sub:ST}
Let the notation be as in \S\ref{subsub:canHodge}. The theory of canonical coordinates (\cite{K79}) asserts that the formal deformation space $\Def(\mathscr{G}/W)$ admits the structure of a formal torus:
\begin{align}\label{eq:serretate}
         \Def(\mathscr{G}/W)\simeq &\Lie U_{\GL,\mu^{-1}}\otimes_{\mathbb{Z}_p}\mathbb{G}_{m,W}^\wedge\\
         \label{eq:serretatevar}\simeq &\Hom_{\bZ_p}\left(\Lie U_{\GL,\mu},\mathbb{G}_{m,W}^\wedge\right).
\end{align}
Here the $\bZ_p$-algebraic groups $U_{\GL,\mu},U_{\GL,\mu^{-1}}$ are the unipotent and opposite unipotent of $\mu$ in $\GL(\omega(\mathscr{G}))$, and $\mathbb{G}_{m,W}^\wedge$ is the formal torus representing the constant multiplicative $p$-divisible group $\mu_{p^\infty,W}$. Note that $\Lie U_{\GL,\mu^{-1}}$ and $\Lie U_{\GL,\mu}$ can be canonically identified with the cocharacter and character lattices of $\Def(\mathscr{G}/W)$. In the same spirit as (\ref{eq:XTp}), we can also canonically identify 
 $\Lie U_{\GL,\mu^{-1}}$ \textit{resp}. $\Lie U_{\GL,\mu}$ with the $\bZ_p$-linear space $X_*(\mathscr{G}^{\loc})\otimes_{\mathbb{Z}_p}  T_p(\mathscr{G}^{\et})^{\vee}$ \textit{resp}.  $X^*(\mathscr{G}^{\loc})\otimes_{\mathbb{Z}_p}  T_p(\mathscr{G}^{\et})$, where $X_*$ stands for the cocharacter lattice.

The identifications (\ref{eq:serretate})$\sim$ (\ref{eq:serretatevar}) will be called the \textbf{canonical coordinates} (or \textbf{Serre--Tate coordinates}) over $\Def(\mathscr{G}/W)$. The group law of the formal torus comes from Baer sums of the extensions. The unique element in $\Def(\mathscr{G}/W)$ corresponding to identity will be called the \textbf{canonical lift}. For a formal $W$-algebra $R$ and a $p$-divisible group ${\widehat{\mathscr{G}}}$ over $\Spf R$ deforming $\mathscr{G}$, (\ref{eq:serretatevar}) yields a $\mathbb{Z}_p$-linear map \begin{equation}\label{eq:serretatepair}
    q_{\widehat{\mathscr{G}}}: \Lie U_{\GL,\mu}\simeq X^*(\mathscr{G}^{\loc})\otimes_{\mathbb{Z}_p}  T_p(\mathscr{G}^{\et})\rightarrow \mathbb{G}^{\wedge}_{m,W}(R),
\end{equation}
called the \textbf{canonical pairing}.

\begin{remark}[Deformation of endomorphisms]\label{subsub:canonicaldefendo} If $\mathscr{G}$ is equipped with a family $V$ of endomorphisms such that each $s\in V$ decomposes as 
$s^{\loc}\times s^{\et}$, then $V$ deforms  endomorphisms of ${\widehat{\mathscr{G}}}$ if and only if \begin{equation}\label{eq:deformendo}    q_{\widehat{\mathscr{G}}}(s^{\loc}x\otimes y)=q_{\widehat{\mathscr{G}}}(x\otimes s^{\et}y),\text{ for all }x\in X^*(\mathscr{G}^{\loc}), y\in T_p(\mathscr{G}^{\et}), s\in V.
\end{equation}
Let $\Def(\mathscr{G},V/W)\subseteq \Def(\mathscr{G}/W)$ be the deformation subspace where $V$ also deforms. Then  there is a $\bZ_p$-sublattice $\Lambda_S\subseteq \Lie U_{\GL,\mu}$ such that $$\Def(\mathscr{G},V/W)=\Hom_{\bZ_p}\left(\Lie U_{\GL,\mu}/\Lambda_S,\mathbb{G}_{m,W}^\wedge\right).$$ Let $\overline{\Lambda_V}\subseteq \Lie U_{\GL,\mu}$ be the saturation of $\Lambda_V$. Then $\Lie U_{\GL,\mu}/\Lambda_V$ decomposes as a direct sum of the free $\bZ_p$-module $\Lie U_{\GL,\mu}/\overline{\Lambda_V}$ and the torsion $\bZ_p$-module 
$\overline{\Lambda_V}/\Lambda_V$. It follows that
\begin{equation}\label{eq:decopmose2222}\Def(\mathscr{G},V/W)=\Hom_{\bZ_p}(\Lie U_{\GL,\mu}/\overline{\Lambda_V}, \bG_m^{\wedge})\times \Hom_{\bZ_p}(\overline{\Lambda_V}/\Lambda_V, \bG_m^{\wedge}).
\end{equation}
Here $\Hom_{\bZ_p}(\Lie U_{\GL,\mu}/\overline{\Lambda_V}, \bG_m^{\wedge})$ is a formal subtorus of $\Def(\mathscr{G}/W)$, while $\Hom_{\bZ_p}(\overline{\Lambda_V}/\Lambda_V, \bG_m^{\wedge})$ is a finite flat subgroup scheme over $W$. In particular, the induced reduced structure $\Def(\mathscr{G},V/W)_{\red}$ is a formal subtorus of $\Def(\mathscr{G}/W)$.

If $\mathscr{G}=A[p^\infty]$ for an abelian variety $A/\Fpbar$. Let $V$ be a collection of endomorphisms as above, and let $\lambda$ be a polarization. Since the deformation of a polarization is again governed by the canonical pairing (\cite[Theorem 2.1]{K81}), a similar analysis shows that $\Def((\mathscr{G},\lambda),V/W)$ is the product of a formal subtorus of $\Def(\mathscr{G}/W)$ with a finite flat group subscheme over $W$.
\end{remark}
\subsubsection{Global canonical pairings}\label{sub:globalST} Let $S$ be a $\bZ_{(p)}$-scheme where $p$ is locally nilpotent. Let $\mathscr{G}$ be an ordinary $p$-divisible group over $S$, so it sits an exact sequence \begin{equation*}
    1\rightarrow \mathscr{G}^{\loc}\rightarrow \mathscr{G}\rightarrow \mathscr{G}^{\et}\rightarrow 1.
\end{equation*}
In this case, $T_p(\mathscr{G}^{\et})$ and $X^*(\mathscr{G}^{\loc})$ are lisse $\bZ_p$-sheaves over $S_{\et}$. Let $\nu_{p^\infty,S}$ be the sheaf over $S_{\et}$ defined by 
$$\varprojlim_n \coker([p^n]:\bG_{m,S}\rightarrow \bG_{m,S}).$$ Then $\mathscr{G}$ admits a \textbf{canonical pairing} (cf. \cite[Definition 2.5]{Ch03}): 
\begin{equation}\label{eq:globalpairing}
q_{\mathscr{G}}\in \Hom_{S_{\et}}\left(X^*(\mathscr{G}^{\loc})\otimes_{\mathbb{Z}_p}  T_p(\mathscr{G}^{\et}),\, \nu_{p^\infty,S}\right)
\end{equation}
that arises from the extension data $[\mathscr{G}]\in \Ext(\mathscr{G}^{\et},\mathscr{G}^{\loc})$. The pairing (\ref{eq:globalpairing}) is functorial in $S$, and admits (\ref{eq:serretatepair}) as a special case. See \cite[\S2]{Ch03} for more details.

\subsubsection{Canonical coordinates for special subvarieties}\label{subsub:canspecialsub} Let $x\in\mathscr{A}^{\ord}_{g,\bK}(\Fpbar)$. It is classically known that $\mathscr{A}^{/x}_{g,\bK,W}\subseteq \Def(A_x[p^{\infty}]/W)$ is a formal subtorus of dimension $g(g+1)/2$ (\cite{K79}). Results of this type also hold for special subvarieties of $\mathscr{A}_{g,\bK}$.

Let the setup be as in \S\ref{subsub:naivespecial}: Let $\mathcal{X}\subseteq \cA_{g,\bK,k}$ be a special subvariety. Let the place of $\mathfrak{p}$ of $k$ be induced by $\bC\simeq \bC_p$. We will always assume that the reduction $\mathscr{X}_{\Fpbar}:=\mathscr{X}\times_{O_{k,(\mathfrak{p})}}\Fpbar$ contains an ordinary point $x$. 
%
A crucial result is the following:

\begin{theorem}[{\cite[Theorem 3.7, Corollary 3.8]{N96}}]\label{thm:noottheorem} The completion $\mathscr{X}^{/x}_{\overline{W}}$ is a finite union of torsion point translates of formal subtori of $\Def(A_x[p^{\infty}]/W)_{\overline{W}}$, and each lift of $x$ to $\mathscr{X}_{\overline{W}}^{/x}$ is contained in a unique irreducible component. Furthermore, $\mathscr{X}_{\overline{W}}$ is flat at $x$, and the normalization of $\mathscr{X}_{\overline{W}}$ is smooth at (points over) $x$.  
\end{theorem} 

When $\mathscr{X}$ is the integral canonical model of a Shimura variety of Hodge type with hyperspecial level structure (\S\ref{subsub:canonicalintegralmodelofHodge}), then the above theorem can be strengthened. Consider a Hodge type Shimura datum $(G,\mathcal{D})$ with a sufficiently small hyperspecial level structure $\bbK=\bbK_p\bbK^p$, with an embedding of Shimura varieties $\Sh_\bbK(G,\mathcal{D})\hookrightarrow \cA_{g,\bK,E}$. Let
$\mathscr{S}_{\bbK}\subseteq \mathscr{A}_{g,\bK,O_{E,(\mathfrak{p})}}$ be the embedding of integral canonical models, such that $\mathscr{S}_{\bbK,\Fpbar}$ contains an ordinary point $x\in \mathscr{A}_{g,\bK}(\Fpbar)$.  
\begin{theorem}[{\cite[Theorem 1.1]{SZ21}}]\label{thm:torusHodge} 
Notation as above, $\mathscr{S}^{/x}_{\bbK,W}\subseteq  \Def(A_x[p^{\infty}]/W)$ is a formal subtorus. In particular, the canonical lift of $x$ is contained in $ \mathscr{S}_{\bbK}(W)$.
\end{theorem}
It is possible to explicitly write down the coordinates of the formal subtori that arise in Theorem~\ref{thm:noottheorem} and  Theorem~\ref{thm:torusHodge} in terms of the Shimura datum. This amounts to identifying Hodge tensors and crystalline tensors in a certain $\bQ_p$-space, which we will do in \S\ref{subsub:pHodge}.

\subsubsection{$p$-adic Hodge comparisons}\label{subsub:pHodge}
Let $x\in \mathscr{A}^{\ord}_{g}(\Fpbar)$. Let $\tilde{x}\in \mathscr{A}^{\ord}_{g}(W)$ be the canonical lift of $x$ and let $\tilde{x}_\bC$ be its base change to $\bC$ along a fixed map $W\hookrightarrow \bC$ (\S\ref{notation}). Identify $H$ with $\bbH_{\mathrm{B},\tilde{x}_\bC}$ as $\bZ_{(p)}$-lattices. Then there are canonical identifications of $\bZ_p$-modules:
\begin{equation}\label{eq:comparison111}
  H_{\bZ_p}\simeq  T_p(A_{\tilde{x}}[p^\infty])^{\vee}\simeq \omega(A_x[p^\infty]).
\end{equation}
The first isomorphism is the Betti-étale comparison, and the second isomorphism follows from the reduction of torsion points of $A_{\tilde{x}}[p^\infty]$ and (\ref{eq:XTp}), which plays a role similar to that of \cite[Corollary 1.4.3(3)]{KM09}. More precisely, via $p$-adic Hodge theory, the Hodge tensors $\{c_{\mathrm{B},\gamma,\tilde{x}_\bC}\}\subseteq H^\otimes$ coming from complex multiplications are sent to Galois invariant tensors $\{c_{p,\et,\gamma,\tilde{x}}\}\subseteq(T_p(A_{\tilde{x}}[p^\infty]))^{\otimes}$, and then to $F$-invariant tensors $\{c_{\cris,\gamma,x}\}\subseteq \Fil^0\bD(A_{x}[p^\infty])^\otimes$ via Fontaine's $D_{\cris}$ functor. Here the filtration on $\bD(A_{x}[p^\infty])$ is nothing other than the one induced by $\mu\otimes W$, which is the filtration corresponding to  $\tilde{x}$ via Grothendieck--Messing theory. From the formation (\ref{eq:fiberfunomega}), a tensor $c_{\cris,\gamma,x}$ canonically gives rise to a tensor $\overline{c}_{\cris,\gamma,x}\subseteq \omega(A_x[p^\infty])^{\otimes}$. The isomorphisms (\ref{eq:comparison111}) have the special property that they carry $c_{\mathrm{B},\gamma,\tilde{x}_\bC}$ to $c_{p,\et,\gamma,\tilde{x}}$, then to $\overline{c}_{\cris,\gamma,x}$. 

The Hodge filtration on $H_{\bC}\simeq H_{\bC_p}$ descends to $H_{\bQ_p}$, and is taken to the filtration on $ \omega(A_x[p^\infty])_{\bQ_p}$ induced by $\mu$ under (\ref{eq:comparison111}). We will abuse the terminology, and say that the Hodge cocharacter $\mu_{\tilde{x}_{\bC}}$ is identified with $\mu$ under the comparison isomorphism.\\

 Let $\{s_\alpha\}\subseteq H^{\otimes}_{\bQ}$ be a collection of tensors fixed by the Mumford--Tate group of $\tilde{x}_{\bC}$, and let $\mathcal{G}\subseteq \GL(H_{\bQ})$ be the stabilizer. By Tannakian formalism, the $p$-adic Hodge comparison maps take each $s_\alpha$ to a Galois invariant tensor $s_{\alpha,p,\et,\tilde{x}}\subseteq (T_p(A_{\tilde{x}}[p^\infty])_{\bQ_p})^{\otimes}$, then to an $F$-invariant tensor $s_{\alpha,\cris,x}\in \Fil^0\bD(A_{x}[p^\infty])_{\bQ_p}^\otimes$. Each $s_{\alpha,\cris,x}$ canonically gives rise to a $\overline{s}_{\alpha,\cris,x}\subseteq (\omega(A_x[p^{\infty}])_{\bQ_p})^\otimes$. Then (\ref{eq:comparison111}) (tensor to $\bQ_p$) takes $s_\alpha$ to $ s_{\alpha,p,\et,\tilde{x}}$ then to $\overline{s}_{\alpha,\cris,x}$. \textbf{Therefore we can canonically view $\mathcal{G}_{\bQ_p}$ as the subgroup of $\GL(\omega(A_x[p^\infty])_{\bQ_p})$ fixing $\{\overline{s}_{\alpha,\cris,x}\}$.} \\

Now let $\mathcal{X}$ and $\mathscr{X}$ be as in the context of Theorem~\ref{thm:noottheorem}. Without loss of generality, we suppose that $\mathcal{X}$ is geometrically irreducible. Let $\tilde{x}'$ be a quasi-canonical lift of $x$ in $\mathscr{X}({\overline{W}})$, whose existence is guaranteed by the toric structure of $\mathscr{X}^{/x}_{\overline{W}}$. We identify the rational Hodge structure of $\tilde{x}'_{\bC}$ with that of the canonical lift $\tilde{x}_{\bC}$, thereby identifying the generic Mumford--Tate group  $\MTT(\mathcal{X}_{\bC})$ as a subgroup of $\GL(H_{\bQ})$. The identification does not depend on the $\tilde{x}'$ chosen.

Take $\{s_\alpha\}\subseteq H_\bQ^\otimes$ to be the tensors cutting out $\MTT(\mathcal{X}_{\bC})$ as a subgroup of $\GL(H_{\bQ})$, which is easily seen fixed by $\MTT(\tilde{x}_{\bC})$. By the process discussed in the previous paragraphs, we see that $\MTT(\mathcal{X}_{\bC})_{\bQ_p}$ is canonically identified as a subgroup of  $\GL(\omega(A_x[p^\infty])_{\bQ_p})$ fixing $\{\overline{s}_{\alpha,\cris,x}\}$. 

Following \S\ref{subsub:Cancoorpdiv}, let $U_{\GL,\mu^{-1}}$ be the opposite unipotent of $\mu$ in $\GL(\omega(A_x[p^{\infty}]))$, and identify the cocharacter lattice of $\Def(A_x[p^{\infty}]/W)$ with $\Lie U_{\GL,\mu^{-1}}$. Let $\Lie U_{\MTT_{\bQ_p},\mu^{-1}}$ be the unipotent of $\MTT(\cX_\bC)_{\bQ_p}\subseteq \GL(\omega(A_x[p^\infty])_{\bQ_p})$ with respect to $\mu^{-1}$.  
\begin{proposition}\label{prop:cocharofSTspecial}
Notation as above. Let $\mathcal{C}\subseteq\mathscr{X}_{\overline{W}}^{/x}$ be the unique irreducible component containing a quasi-canonical lift $\tilde{x}_{\overline{W}}$, as predicted by Theorem~\ref{thm:noottheorem}. Let $\mathcal{C}_0$ be the subtorus of $\Def(A_x[p^{\infty}]/W)_{\overline{W}}$ determined by the rational cocharacter lattice $$\Lie U_{\MTT_{\bQ_p},\mu^{-1}} \subseteq \Lie U_{\GL,\mu^{-1},\bQ_p}.$$
Then $\mathcal{C}$ is the translate of $\mathcal{C}_0$ by $\tilde{x}_{\overline{W}}$.  
\end{proposition}
In the context of Theorem~\ref{thm:torusHodge}, we can take $\tilde{x}_{\overline{W}}$ to be the canonical lift, which descends to $W$. Instead of taking the generic Mumford--Tate group $\MTT(\Sh_\bbK(G,D)_{\bC})$, we can take the structure group $G$ of the Shimura variety. Since $\Lie_{G_{\bQ_p},\mu^{-1}}=\Lie_{\MTT_{\bQ_p},\mu^{-1}}$, the proposition implies that $\mathscr{S}^{/x}_{\bbK,W}$ is the subtorus of $\Def(A_x[p^{\infty}]/W)$ with rational cocharacter lattice $\Lie U_{G_{\bQ_p},\mu^{-1}}$. This is also proven in \cite[Proposition 5.4.3]{Hften21}. 

We will give a proof of Proposition~\ref{prop:cocharofSTspecial} using different techniques in \S\ref{subsub:tmIMPLICA}.

\subsection{Quasi-weakly special subvarieties}\label{sub:Quasi-weakly} 
In this section we introduce the notion of quasi-weakly special subvarieties for $\mathscr{A}_{g,\Fpbar}^{\ord}$, which is an analogue of the notion of weakly special subvarieties in char 0, and is needed for stating the char $p$ analogue of the André--Oort conjecture. The classical definition of weakly special subvarieties can be found in \cite[\S 3.1]{M97}. Roughly speaking, it is a product of a special subvariety with a possibly non-special point. 

Let the setup being the same as in \S\ref{subsub:naivespecial}: $k\subseteq\bC$ is a number field or $\overline{\bQ}$, $\mathfrak{p}$ is a prime of $O_{k}$ induced from the fixed isomorphism $\bC\simeq \bC_p$ (\S\ref{notation}), 
$\mathcal{X}\subseteq \mathscr{A}_{g,\bK,k}$ is a special subvariety, and let $\mathscr{X}/O_{k,(\mathfrak{p})}$ be the naïve integral model. 
We define the \textbf{ordinary naïve integral model} ${\mathscr{X}}^{\ord}\subseteq \mathscr{X}$ as the open subscheme by throwing out the non-ordinary loci of the special fiber.

\subsubsection{Special correspondences}\label{subsub:specialcorr}
Suppose that $\mathcal{X}\subseteq \mathcal{A}_{g_1,\bK_1,k}$ and $\mathcal{Y}\subseteq \mathcal{A}_{g_2,\bK_2,k}$ are geometrically connected special subvarieties. A \textbf{special correspondence} between $\mathcal{X}$ and $\mathcal{Y}$ is a special subvariety $\mathcal{W}\subseteq \mathcal{X}\times \mathcal{Y}$ such that the projection of $\mathcal{W}$ to $\mathcal{X}$ and $\mathcal{Y}$ are finite surjective. We denote it by $\mathcal{X}\xrightarrow{\mathcal{W}}\mathcal{Y}$. 

Let $\mathscr{W}\subseteq \mathscr{X}\times \mathscr{Y}$ be the ordinary naïve integral model of $\mathcal{W}$ in $\mathscr{A}_{g_1,\bK_1,k}\times \mathscr{A}_{g_2,\bK_2,k}$. Intuitively, $\mathscr{W}$ may be viewed as 
an integral correspondence between $\mathscr{X}$ and $\mathscr{Y}$, though there are more technical issues. For example, it is not clear if the projection of $\mathscr{W}$ to $\mathscr{X}$ or $\mathscr{Y}$ is still finite surjective. 

However, if we restrict ourselves to ordinary strata, things behave much better, see Lemma~\ref{lm:extensiontonormal}. 
Let $X\subseteq \mathscr{X}_{\Fpbar}^{\ord}$ and $Y\subseteq \mathscr{Y}_{\Fpbar}^{\ord}$ be irreducible subvarieties. We say $X$ and $Y$ \textbf{specially correspond with respect to} $\mathcal{W}$ if there is an irreducible subvariety $Z\subseteq \mathscr{W}_{\Fpbar}^{\ord}$, such that the projection of $Z$ to $\mathscr{X}_{\Fpbar}^{\ord}$ admits the same Zariski closure as $X$, and the projection of $Z$ to $\mathscr{Y}_{\Fpbar}^{\ord}$ admits the same Zariski closure as $Y$. We say $X$ and $Y$ \textbf{specially correspond} if they specially correspond with respect to some $\mathcal{W}$.

\begin{lemma}\label{lm:extensiontonormal}
Let $\mathcal{X},\mathscr{X}$, $\mathcal{Y},\mathscr{Y}$, and $\mathcal{W},\mathscr{W}$ be as in \S\ref{subsub:specialcorr}, and let $\pi_1:\mathscr{W}^{\ord}\rightarrow \mathscr{X}^{\ord}$ and $\pi_2:\mathscr{W}^{\ord}\rightarrow \mathscr{Y}^{\ord}$. The following are true:
\begin{enumerate}
    \item\label{lm:extnaiveit00}
  $\pi_1,\pi_2$ are quasi-finite,
  \item\label{lm:extnaiveit01}
If $\mathscr{X}^{\ord}_{\Fpbar}$ and $\mathscr{Y}^{\ord}_{\Fpbar}$ are nonempty, then $\pi_1,\pi_2$ are surjective.
\end{enumerate}

\end{lemma}
\proof\begin{enumerate}
    \item  It suffices to show that $\pi_{1,\Fpbar}:\mathscr{W}^{\ord}_{\Fpbar}\rightarrow \mathscr{X}^{\ord}_{\Fpbar}$ is quasi-finite. Using the property of Jacobson schemes, one reduces to showing that an $\Fpbar$-point of $\mathscr{X}^{\ord}_{\Fpbar}$ admits finite fiber. 

If $\mathscr{W}^{\ord}(\Fpbar)=\emptyset$, then the claim is true. So we will assume that $\mathscr{W}^{\ord}(\Fpbar)\neq \emptyset$. Let $z\in \mathscr{W}^{\ord}(\Fpbar)$. By abuse of notation, we will also denote by $z$ the point $\pi_{1}(z)$. Since $\mathcal{X}$ and $\mathcal{W}$ are special, by Theorem~\ref{thm:noottheorem}, $\mathscr{X}^{/z}_{\overline{W}}\subseteq \mathscr{A}_{g_1,\bK_1,\overline{W}}^{/z}$ and 
$\mathscr{W}^{/z}_{\overline{W}}\subseteq \mathscr{A}_{g_1,\bK_1,\overline{W}}^{/z}\times \mathscr{A}_{g_2,\bK_2,\overline{W}}^{/z}$ are finite union of torsion translates of formal subtori. Since $\mathscr{A}_{g_1,\bK_1,\overline{W}}^{/z}\times \mathscr{A}_{g_2,\bK_2,\overline{W}}^{/z}\rightarrow \mathscr{A}_{g_1,\bK_1,\overline{W}}^{/z}$ is a morphism of formal tori, the projection $\pi_{1}^{/x}:\mathscr{W}^{/z}_{\overline{W}}\rightarrow \mathscr{X}^{/z}_{\overline{W}}$ is \textbf{linear} in the following sense: If $\mathcal{C}$ is an irreducible component of the source such that $\pi_{1}^{/z}|_{\mathcal{C}}$ lands in an irreducible component $\mathcal{D}$ of the target, then $\mathcal{C}$ (\textit{resp}. $\mathcal{D}$) is a translate of a formal subtori $\mathcal{C}_0$ (\textit{resp}. $\mathcal{D}_0$) by a torsion point $c$ (\textit{resp}. $d$), and there exists a homomorphism of formal tori $\phi:\mathcal{C}_0\rightarrow \mathcal{D}_0$ such that $\pi_{1}^{/z}|_{\mathcal{C}}$ is the composition 
$\mathcal{C}\xrightarrow{-c} \mathcal{C}_0 \xrightarrow{\phi} \mathcal{D}_0 \xrightarrow{+d} \mathcal{D}$. Now that $\pi_{1}^{/z}$ is finite over generic fiber, by linearity, we deduce that $\pi_{1}^{/z}$ is finite. 

If some point of ${\mathscr{X}}^{\ord}(\Fpbar)$ admits a positive dimensional pre-image in $\mathscr{W}^{\ord}(\Fpbar)$, we pick $z\in \mathscr{W}^{\ord}(\Fpbar)$ sitting inside a positive dimensional component of its preimage. But then $\pi_{1}^{/z}$ can not be finite. It follows that $\pi_{1}$ is quasi-finite.

\item It suffices to check that $\pi_{1,\Fpbar}$ is surjective over $\Fpbar$-points. The following argument can be seen as ``verifying the valuative criterion for CM points''. Let $x\in \mathscr{X}^{\ord}(\Fpbar)$. We quasi-canonical lift it to a CM point $\tilde{x}_{\overline{K}}\in \mathcal{X}(\overline{K})$. Since $\pi_1$ is finite surjective over the generic fibers, we can lift $\tilde{x}_{\overline{K}}$ to a point $\tilde{y}_{\overline{K}}\in \mathcal{W}(\overline{K})$. \begin{claim}
    $\tilde{y}_{\overline{K}}$ specializes to a point $y$ in $\mathscr{W}^{\ord}(\Fpbar)$.
\end{claim} Once this is done, then $y$ is a preimage of $x$ under $\pi_{1,\Fpbar}$. In particular, $\pi_{1,\Fpbar}$ is surjective.

Now $\tilde{y}_{\overline{K}}$ is a CM point, so it has good reduction. As a result, it always specializes to a point $y\in \mathscr{W}(\Fpbar)$. We need to show that $y$ is ordinary. Now let's base change from $\overline{K}$ to $\bC$. Let $\MTT(\mathcal{W}_{\bC})$ be the generic Mumford--Tate group of $\mathcal{W}_{{\bC}}$ for the universal abelian scheme on $\mathcal{A}_{g_1,\bK_1,\bC}\times \mathcal{A}_{g_2,\bK_2,\bC}$, let $\MTT(\mathcal{X}_{\bC})$ (\textit{resp}. $\MTT(\mathcal{Y}_{\bC})$) be the generic Mumford--Tate group of $\mathcal{X}_{\bC}$ (\textit{resp}. $\mathcal{Y}_{\bC}$) for the universal abelian scheme on $\mathcal{A}_{g_1,\bK_1,\bC}$ (\textit{resp}. $ \mathcal{A}_{g_2,\bK_2,\bC}$). Then we have an embedding $\MTT(\mathcal{W}_\bC)\hookrightarrow \MTT(\mathcal{X}_\bC)\times \MTT(\mathcal{Y}_\bC)$ that surjects onto each factors. Furthermore, the projection of $\MTT(\mathcal{W}_\bC)$ onto each factors induces isomorphisms on adjoints. This implies that $\ker(\MTT(\mathcal{W}_\bC)\rightarrow \MTT(\mathcal{X}_\bC))$ is central in $\MTT(\mathcal{W}_\bC)$. Let $T:=\ker(\MTT(\tilde{y}_\bC)\rightarrow \MTT(\tilde{x}_\bC))\subseteq \ker(\MTT(\mathcal{W}_\bC)\rightarrow \MTT(\mathcal{X}_\bC))$. Then the projection embeds $T$ into the center of $\MTT(\mathcal{Y}_\bC)$. 

By \cite[Lemma 2.2]{B01}, a CM abelian variety over $\overline{\bQ}$ has ordinary reduction at $\mathfrak{p}$ if and only if $\mathfrak{p}$ contracts to a prime of height 1 over its reflex field (i.e., the local field at this prime is the trivial extension of $\bQ_p$). Let $E'$ (\textit{resp}. $E$) be the reflex field of $\tilde{y}_{\bC}$  (\textit{resp}. $\tilde{x}_{\bC}$). Suppose that $\tilde{y}_{\overline{K}}$ has non-ordinary reduction. Since $\tilde{x}_{\overline{K}}$ has ordinary reduction, we have $E\subsetneq E'$. In particular, $T$ is nontrivial, and more is true: let $T'$ be the image of $T$ in $\MTT(\mathcal{W}_\bC)^{\mathrm{ab}}$, then  the Hodge cocharacter $h_{\tilde{y}_{\bC}}$ factors nontrivially through $T'_{\bC}$, and the smallest number field that the image can be descent to is $E'$. If $z_{\bC}$ is another CM point of $\mathcal{W}_\bC$, {then $h_{z_{\bC}}$ is conjugate to $h_{\tilde{y}_{\bC}}$ over $\bC$}. Since $T$ is central, this implies that every CM point in $\mathcal{W}_\bC$ has reflex field containing $E'$, hence having non-ordinary reduction. Since the projection embeds $T$ isomorphically into the center of $\MTT(\mathcal{Y}_\bC)$, a similar argument shows that every CM point in $\mathcal{Y}_\bC$ has non-ordinary reduction. This contradicts the fact that $\mathscr{Y}^{\ord}_{\Fpbar}\neq\emptyset$, since by Theorem~\ref{thm:noottheorem} we can quasi-canonically lift an element in $\mathscr{Y}^{\ord}_{\Fpbar}$ to a CM point in $\mathcal{Y}_\bC$ with ordinary reduction. This proves the claim and establishes (\ref{lm:extnaiveit01}).\end{enumerate}
$\hfill\square$

\begin{defn}\label{def:quasi-weakly}
Notation as above. Suppose that there is a special subvariety $\mathcal{X}\subseteq \cA_{g,\bK,k}$ that arises from an embedding of Shimura data $(G,\mathcal{D})\hookrightarrow(\GSp(H_\bQ),\mathfrak{H}^{\pm})$. An irreducible subvariety   $X\subseteq\mathscr{X}_{\Fpbar}^{\ord}$ is called \textbf{quasi-weakly special}, if it is either an $\Fpbar$-point, or there are
\begin{enumerate}
    \item for each $\alpha\in\{0,1\}$, an embedding of Shimura data  $(G_\alpha,\mathcal{D}_\alpha)\hookrightarrow(\GSp(H_{\alpha,\bQ}),\mathfrak{H}^{\pm}_\alpha)$, such that  $(G^{\ad},\mathcal{D}^{\ad})= (G_0^{\ad},\mathcal{D}_0^{\ad})\times (G_1^{\ad},\mathcal{D}_1^{\ad})$,
    \item 
    for each $\alpha\in\{0,1\}$, 
    a special subvariety $\mathcal{X}_\alpha\subseteq \cA_{g_\alpha,\bK_\alpha,k}$ that arises from the embedding of Shimura data, as well as its ordinary naïve integral model $\mathscr{X}_\alpha^{\ord}$, such that $\dim\mathscr{X}_{0,\Fpbar}^{\ord}>0$, 
\item an irreducible component $C\subseteq \mathscr{X}_{0,\Fpbar}^{\ord}$ and a subvariety $X_1\subseteq {\mathscr{X}}_{1,\Fpbar}^{\ord}$,
\end{enumerate}
such that $X\subseteq \mathscr{X}_{\Fpbar}^{\ord}$ specially correspond to $C\times X_1\subseteq  \mathscr{X}_{0,\Fpbar}^{\ord}\times \mathscr{X}_{1,\Fpbar}^{\ord}$. In light of Lemma~\ref{lm:extensiontonormal}, $X$ is an \textbf{almost product} of $C$ with $X_1$.
\end{defn}
\subsubsection{Quasi-lifts}\label{subsub:Quasi-lifts} 
We introduce a very simple construction which, however, will play a significant role later. The set up is as follows: 

Let $\varpi:{V}'\rightarrow {V}$ be a morphism of finite type $\Fpbar$-schemes. Let $(X,x)$ be a pointed connected smooth $\Fpbar$-variety that admits a morphism 
$f:(X,x) \rightarrow {V}$. Let $(X',x')$ be a pointed smooth $\Fpbar$-variety. We say that an $f':(X',x')\rightarrow V'$ 
is a \textbf{quasi-lift of $f$ to $V'$ along $\varpi$}, if there is a surjection $g:(X',x')\rightarrow (X,x)$, such that $\varpi\comp f'=f\comp g$. 

\begin{lemma}\label{lm:quasi-lift}
If $\varpi$ is surjective, then a morphism $f:(X,x)\rightarrow V$  admits a quasi-lift to $V'$ along $\varpi$. 
\end{lemma}
\begin{proof}
Let $X''$ be an irreducible component of  $X\times_{V} V'$ that surjects onto $X$. It is equipped with a morphism $f'':X''\rightarrow V'$. Let $X'$ be any smooth variety that admits a surjection $X'\xrightarrow{h} X''$. This always exists, e.g., de Jong's alteration of $X''$. Now let $f'=f''\comp h$, $g=(X''\rightarrow X)\comp h$, and $x'\in g^{-1}(x)$.
\end{proof}

\section{Monodromy of local systems}\label{sec:monodromy} We review the notion of monodromy for étale lisse sheaves and $F$-isocrystals, which will play a fundamental role in our paper. For simplicity, we will mostly stick to local systems with coefficients in $\bQ_u$, where $u$ is a finite place of $\bQ$, but the treatment extends to local systems with coefficients in finite extensions of $\bQ_u$. For a more comprehensive and broader understanding of these notions, readers are suggested to refer to \cite[\S2,\S3]{DA20}. We will always assume that $X_0$ is a geometric connected smooth variety over a finite field $\mathbb{F}_q$ with $X=(X_{0})_{\Fpbar}$, and $x$ is an $\Fpbar$-point of $X_0$. 

\subsection{Monodromy of étale lisse sheaves}\label{subsub:monolisse} Let $u$ be a finite place of $\bQ$, including $p$. Consider $\LS(X_0,\bQ_u)$, the category of étale lisse sheaves of $\bQ_u$-vector spaces over $X_0$. This category is equivalent to the category of continuous $\pi_1^{\et}(X_0,x)$-representations. It is also a neutral Tannakian category with fiber functor 
\begin{align*}
 \omega_x: \LS(X_0,\bQ_u)&\rightarrow \Vect_{\bQ_u}\\
\mathcal{E}&\rightarrow \mathcal{E}_x.
\end{align*}  
The \textbf{monodromy group} of an object $\mathcal{E}$ in $\LS(X_0,\bQ_u)$ at $x$ is the Tannakian fundamental group of the tensor abelian subcategory $\langle\mathcal{E}\rangle^{\otimes}$ with fiber functor $\omega_x$, denoted $G(\mathcal{E},x)$. Since $\LS(X_0,\bQ_u)$ is equivalent to the category of continuous $\pi_1^{\et}(X_0,x)$-representations, $G(\mathcal{E},x)$ is nothing other than the Zariski closure of the image of $\pi_1^{\et}(X_0,x)$ in $\GL(\mathcal{E}_x)$. 

There is also a notion of Weil lisse sheaves, which is more widely used (cf. \cite{DA20}). A Weil lisse sheaf over $X_0$ is a (geometric) étale lisse sheaf $\mathcal{V}$ over $X$, together with a Frobenius structure $F^*\mathcal{V}\xrightarrow{\sim} \mathcal{V}$, where $F$ is the geometric Frobenius of $\bF_q$ with respect to $\bF$. Let $W(X_0,x)\subseteq \pi_1^{\et}(X_0,x)$ be the Weil group of $X_0$. The category of Weil lisse sheaves is equivalent to the category of continuous $W(X_0,x)$-representations, and is a neutral Tannakian category with fiber functor $\omega_x$. The monodromy group at $x$ of a Weil lisse sheaf $\mathcal{V}$ is the Tannakian fundamental group of $\langle\mathcal{V}\rangle^{\otimes}$ with fiber $\omega_x$. It is the Zaraski closure of $W(X_0,x)$ in $\GL(\mathcal{V}_x)$. 

An étale lisse sheaf $\mathcal{E}$ in $\LS(X_0,\bQ_u)$ is automatically a Weil lisse sheaf via pullback, and its subquotients as Weil lisse sheaves are objects in $\LS(X_0,\bQ_u)$. In other words, the monodromy group of $\mathcal{E}$ as an étale lisse sheaf over $X_0$ equals the monodromy group of $\mathcal{E}$ as a Weil lisse sheaf. In this paper, we almost only work in the category $\LS(X_0,\bQ_u)$.

\subsection{Monodromy of $F$-isocrystals}\label{subsec:lgmono} 
In this section, we use $F$ to denote the absolute Frobenius. Let $\FIsoc(X_0)$ (\textit{resp}. $\FIsoc(X)$) be the tensor abelian category of $F$-isocrystals over $X_0$ (\textit{resp}. $X$). Consider an object $\mathcal{M}$ in $\FIsoc(X_0)$ or $\FIsoc(X)$. We denote $\la\mathcal{M}\ra^{\otimes}$ the tensor abelian subcategory generated by $\mathcal{M}$. Let $e$ be the smallest positive integer such that the slopes of $\mathcal{M}_x$ multiplied by $e$ lie in $\mathbb{Z}$. The fibre functor \begin{align*}
 \omega_x: \langle \mathcal{M}\rangle^\otimes_{\mathbb{Q}_{p^e}} &\rightarrow \Vect_{\mathbb{Q}_{p^e}}\\
(\mathcal{N},F)&\rightarrow \{v\in \mathcal{N}_x|\exists i\in \mathbb{Z},\,\,({F^e_{x}}-p^i)v=0\}.
\end{align*}   
makes $\la\mathcal{M}\ra^{\otimes}_{\mathbb{Q}_{p^e}}$ -- the scalar extension of $\la\mathcal{M}\ra^{\otimes}$ by $\mathbb{Q}_{p^e}$ (cf.\cite[\S 1.4.1]{Abe}) -- a neutral Tannakian category. The functor $\omega_x$ is essentially the same as the \textbf{Dieudonné--Manin fiber functor} in \cite[Construction 3.1.4, Definition 3.1.6]{MD20}, and is also a minor generalization of the fiber functor found in \cite{Ch03}. 
The fundamental group $\text{Aut}^\otimes(\omega_x)\subseteq \GL(\omega_x(\mathcal{M}))$ is called the \textbf{monodromy group} of $\mathcal{M}$ at $x$, denoted $G(\mathcal{M},x)$. 

There is also a local version of the monodromy groups. We won't use it in the paper, but still include it for completeness.  Let $\mathcal{M}^{/x}\in \FIsoc(X^{/x})$ be the base change of $\mathcal{M}$. The subcategory $\la\mathcal{M}^{/x}\ra^{\otimes}_{\mathbb{Q}_{p^e}}$ is again a Tannakian category with the fiber functor $\omega_x$, see \cite[\S 3.3]{Drinfeld2022}. The corresponding monodromy group is called the \textbf{local monodromy group} of $\mathcal{M}$ at $x$, denoted $G(\mathcal{M}^{/x},x)$. 

We have $G(\mathcal{M}^{/x},x)\subseteq G(\mathcal{M},x)\subseteq \GL(\omega_x(\mathcal{M}))$. When the base point is clear from the context, we will simply drop it, and write $G(\mathcal{M})$, $G(\mathcal{M}^{/x})$, \textit{etc}.

We will mainly be interested in $F$-isocrystals with constant slopes. If $\mathcal{M}$ has constant slopes, then \cite[Corollary 2.6.2]{K79} and \cite[Corollary 4.2]{Ked22} imply that $\mathcal{M}$ admits the slope filtration $$0=\mathcal{M}_0\subseteq...\subseteq \mathcal{M}_{l}=\mathcal{M},$$ 
where each graded piece $\mathcal{M}_{i}/\mathcal{M}_{i-1}$ has pure slope $s_i\in \bQ$ and $s_1<...<s_l$. We will write $$\gr\mathcal{M}=\bigoplus_{i=1}^l\mathcal{M}_{i}/\mathcal{M}_{i-1}.$$ Let $U(\mathcal{M},x)$ \textit{resp}. $U(\mathcal{M}^{/x},x)$ be the kernel of the natural projection $G(\mathcal{M},x)\rightarrow G(\gr\mathcal{M},x)$ \textit{resp}. $G(\mathcal{M}^{/x},x)\rightarrow G(\gr\mathcal{M}^{/x},x)$. They are all unipotent, with $U(\mathcal{M}^{/x},x)\subseteq U(\mathcal{M},x)$. The monodromy groups for an $F$-isocrystal with constant slopes is relatively easy to understand:
\begin{lemma}\label{strofmon}
Suppose that $\mathcal{M}$ has constant slopes and let $\nu$ be the Newton cocharacter of $\mathcal{M}_x$. Identify $G(\mathcal{M},x), G(\gr\mathcal{M},x), U(\mathcal{M},x)$ and their local counterparts as subgroups of $\GL(\omega_x(\mathcal{M}))$. The following are true: \begin{enumerate}
   \item $G(\gr\mathcal{M}^{/x},{x})=\im\nu$,
    \item 
$G(\mathcal{M},x)=U(\mathcal{M},x)\rtimes G(\gr\mathcal{M},x)$ via adjoint action,
\item $G(\mathcal{M}^{/x},x)= U(\mathcal{M}^{/x},x)\rtimes G(\gr\mathcal{M}^{/x},x)$ via adjoint action.
\end{enumerate}
\end{lemma}
\proof~\begin{enumerate}
    \item follows since a uniroot $F$- isocrystal over $X^{/x}_{\Fpbar}$ is constant.
    \item  There is a map $\gr:\langle \mathcal{M}\rangle^\otimes\subseteq\langle\gr\mathcal{M}\rangle^\otimes$ sending an $F$-isocrystal to its graded object, inducing a section $G(\gr\mathcal{M},x)\inj G(\mathcal{M},x)$ to the natural map $G(\mathcal{M},x)\rightarrow G(\gr\mathcal{M},x)$, hence we have the semi-direct product. The claim that  $G(\gr\mathcal{M},x)$ acts on $U(\mathcal{M},x)$ via adjoint action is clear from the way that they embed into $\GL(\omega_x(\mathcal{M}))$.
    \item is similar to (2). $\hfill\square$
\end{enumerate} 
\subsubsection{Monodromy of ordinary $p$-divisible groups}\label{Sec:2slope} We now review Chai's result on local and global monodromy of ordinary $p$-divisible groups. Let $\mathscr{G}_0$ be an ordinary $p$-divisible group over $X_0$, which is an extension of $\mathscr{G}^{\loc}_0$ and $\mathscr{G}^{\et}_0$. Let $\mathscr{G}$ (\textit{resp}.  $\mathscr{G}^{/x}$) be the base change of $\mathscr{G}_0$ to $X$ (\textit{resp}.  $X^{/x}$). 

Write $\mathcal{M}=\bD(\mathscr{G}_0)_{\bQ_p}$. Let $\mu: \bG_m \rightarrow \GL(\omega_x(\mathcal{M}))$ be the canonical Hodge cocharacter of $\mathcal{M}_x$. Routinely, the notations $ U_{\GL,\mu}$ and $ U_{\GL,\mu^{-1}}$ refer to the unipotent and the opposite unipotent of $\mu$ in $\GL(\omega_x(\mathcal{M}))$. Consider the canonical pairing:
$$q_{\mathscr{G}^{/x}}:\Lie U_{\GL,\mu}\rightarrow \mathbb{G}^{\wedge}_m(X^{/x}).$$
Define $N(\mathscr{G}^{/x})=\ker(q_{\mathscr{G}^{/x}})^{\perp}$, the sub-lattice of $\Lie U_{\GL,\mu^{-1}}$ which pairs to 0 with $\ker(q_{\mathscr{G}^{/x}})$. Via exponential map, it can also be viewed as a unipotent subgroup of $U_{\GL,\mu^{-1}}$. 
\begin{theorem}[Chai]\label{thm:etale-glob} 
Notations as above. Identify $G(\mathcal{M},x), G(\gr\mathcal{M},x), U(\mathcal{M},x)$ and their local counterparts as subgroups of $\GL(\omega_x(\mathcal{M}))$ and regard $N(\mathscr{G}^{/x})_{\bQ_p}$ as a unipotent subgroup of $U_{\GL,\mu^{-1},\bQ_p}$. We have 
\begin{enumerate}
    \item\label{Chai:mono1} $U(\mathcal{M}^{/x},x)= N(\mathscr{G}^{/x})_{\bQ_p}$ and $G(\mathcal{M}^{/x},x)=N(\mathscr{G}^{/x})_{\bQ_p}\rtimes \im \mu$ via adjoint action. 
    \item\label{Chai:mono2} $U(\mathcal{M},x)= N(\mathscr{G}^{/x})_{\bQ_p}$ and $G(\mathcal{M},x)=N(\mathscr{G}^{/x})_{\bQ_p}\rtimes G(\gr\mathcal{M},x)$ via adjoint action.
\end{enumerate}  
\end{theorem}
\proof Point (\ref{Chai:mono1}) follows from \cite[Theorem 3.3]{Ch03}. Now we prove (\ref{Chai:mono2}). 
Recall from \S\ref{sub:globalST} that $\mathscr{G}_0$ admits a canonical pairing $$q_{\mathscr{G}_0}:X^*(\mathscr{G}^{\loc}_0)\otimes_{\mathbb{Z}_p}  T_p(\mathscr{G}^{\et}_0)\rightarrow \nu_{p^\infty,X_0}.$$
By \cite[Proposition 4.2.1]{Ch03}, $\ker q_{\mathscr{G}_0}$ is a saturated lisse subsheaf (Chai only proved it when the base field is $\Fpbar$, but the proof generalizes to $\Fpbar_q$). Let $N(\mathscr{G}_0):=(\ker q_{\mathscr{G}_0})^\perp \subseteq X_*(\mathscr{G}^{\loc}_0)\otimes_{\mathbb{Z}_p}  T_p(\mathscr{G}^{\et}_0)^\vee$ be the saturated lisse subsheaf that pairs to 0 with  $\ker q_{\mathscr{G}_0}$ (cf. \cite[\S4.2]{Ch03}). From \textit{loc.cit} Theorem 4.4, we know that $U(\cM,x)=N(\mathscr{G}_0)_{\bQ_p,x}$ (again, Chai's setting is over $\Fpbar$, but the proof generalizes to $\Fpbar_q$). So it suffices to show that $\rk N(\mathscr{G}_0)=\rk N(\mathscr{G}^{/x})$, which further reduces to showing that $\rk \ker (q_{\mathscr{G}^{/x}})=\rk \ker (q_{\mathscr{G}_0})$.

Let $\Spec R\rightarrow X$ be the étale local ring at $x$ (i.e., the strict Henselianzation at $x$), and let $X^{/x}=\Spf \widehat{R}$. By the functoriality of canonical pairings, the pullback of $q_{\mathscr{G}_0}$ to $\Spec R$ (\textit{resp}. $X^{/x}$) can be identified with $q_{\mathscr{G}_R}$ (\textit{resp}. $q_{\mathscr{G}^{/x}}$), where ${\mathscr{G}_R}$ is the pullback of $\mathscr{G}$ to $\Spec R$. Since $\ker q_{\mathscr{G}_0}$ is lisse, its rank equals the rank of $\ker q_{\mathscr{G}_R}$. So it suffices to show that $\rk \ker (q_{\mathscr{G}^{/x}})=\rk \ker (q_{\mathscr{G}_R})$.

Note that $X^*(\mathscr{G}^{\loc})\otimes_{\mathbb{Z}_p}  T_p(\mathscr{G}^{\et})$ is constant over both $R$ and $\widehat{R}$, and that $\nu_{p^\infty,X}(R)=\varprojlim_n R^*/(R^*)^{p^n}$, and $\bG_{m}^\wedge(X^{/x})=\varprojlim_n \widehat{R}^*/(\widehat{R}^*)^{p^n}$. To show that $\rk \ker (q_{\mathscr{G}^{/x}})=\rk \ker (q_{\mathscr{G}_R})$, it suffices that for every $n>0$, the natural map
$$R^*/(R^*)^{p^n}\rightarrow \widehat{R}^*/(\widehat{R}^*)^{p^n}$$
 is an injection. This follows from \cite[Lemma 2.1.4$\sim$2.1.6] {Ch03}. 
 
 See \cite[Theorem II]{daddezio2022hecke} for a generalization of the theorem.
$\hfill\square$

\subsection{Monodromy of overconvergent $F$-isocrystals}\label{subsec:monoovercon} We use the same setups and notation as in \S\ref{subsec:lgmono}. Let's denote $\FIsoc^\dagger(X_0)$ as the tensor abelian category of overconvergent $F$-isocrystals over $X_0$. There is a forgetful functor $\text{Fgt}: \FIsoc^\dagger(X_0)\rightarrow \FIsoc(X_0)$. Let $\mathcal{M}^\dagger=(\mathcal{M}^\dagger,F^\dagger)\in \FIsoc^\dagger(X_0)$ and let $\mathcal{M}=(\mathcal{M},F)$ be its image in $\FIsoc(X_0)$ forgetting the overconvergent structure. Recall that $e$ is the smallest positive integer such that the slopes of $\mathcal{M}_x$ multiplied by $e$ lie in $\mathbb{Z}$. The fiber functor $$\omega_x^\dagger=\omega_x\comp \text{Fgt}:\langle \mathcal{M}^\dagger\rangle^\otimes_{\mathbb{Q}_{p^e}}\rightarrow \text{Vect}_{\mathbb{Q}_{p^e}}$$ 
makes $\la\mathcal{M}^\dagger\ra^\otimes_{\mathbb{Q}_{p^e}}$ a neutral Tannakian category. The Tannakian fundamental group thus arises is called the (overconvergent) monodromy group of $\mathcal{M}^\dagger$ at $x$, denoted $G(\mathcal{M}^\dagger,x)$. Note that we have $G(\mathcal{M},x)\subseteq G(\mathcal{M}^\dagger,x)\subseteq \GL(\omega_x(\mathcal{M}))$. Again, when the base point is clear from the context, we will simply drop it, and write $G(\mathcal{M}^\dagger)$.
\begin{theorem}[D'Addezio]\label{MDA}
 The following are true:\begin{enumerate}
    \item Suppose $\mathcal{M}^\dagger$ admits slope filtration, then $G(\mathcal{M},x)\subseteq G(\mathcal{M}^\dagger,x)$ is the parabolic subgroup fixing the slope filtration of $\mathcal{M}_x$.
    \item Suppose $g:A\rightarrow X_0$ is an abelian scheme. Let $D^\dagger(A)=R^1g_{*,\cris}\mathcal{O}_{A,\cris}$, then the overconvergent monodromy group $G(D^\dagger(A),x)$ is reductive.
    \item Setup being the same as (2). If $\mathcal{M}^\dagger$ is an overconvergent $F$-isocrystal in $\langle D^\dagger(A)\rangle^\otimes$ that has constant slopes, then $G(\gr\mathcal{M},x)$ is reductive.
\end{enumerate}
\end{theorem}
\proof~\begin{enumerate}
    \item follows from \cite[Theorem 5.1.2]{MD20}. 
    \item is  \cite[Corollary 3.5.2]{DA20}.
    \item  Since $G(D^\dagger(A),x)$ is reductive by (2), the group $G(\mathcal{M}^\dagger,x)$ is also reductive. By (1), $G(\mathcal{M},x)$ is the parabolic subgroup of $G(\mathcal{M}^\dagger,x)$ fixing the slope filtration on $\mathcal{M}_x$. Since $G(\gr\mathcal{M},x)$ is a quotient of $G(\mathcal{M},x)$ by Lemma~\ref{strofmon}(3), it is also reductive.$\hfill\square$
\end{enumerate}

\subsection{Monodromy of compatible systems}\label{subsub:compatible} We refer to \cite[\S 3]{DA20} and \cite[\S 2]{Ked3} for the notion of coefficient objects and compatible systems. Roughly speaking, a coefficient object over $X_0$ is an $l$-adic lisse sheaf or an overconvergent $F$-isocrystal. A coefficient object $\mathcal{E}_0$ over $X_0$ is equipped with a characteristic polynomial $P(\mathcal{E}_{x_0},t)$ for every closed point $x_0\in |X_0|$. 
For a number field $E$, together with a set of finite places $\Sigma$ containing all places that do not divide $p$, an \textbf{$E$-compatible system of coefficient objects}  $\{\mathcal{E}_{\lambda,0}\}_{\lambda\in \Sigma}$ is a collection of coefficient objects $\mathcal{E}_{\lambda,0}$ over $E_\lambda$ indexed by $\lambda\in \Sigma$, such that the characteristic polynomials of all $\mathcal{E}_{\lambda,0}$ at all closed points of $X_0$ are equal and are $E$-polynomials. In \cite[Theorem 1.2.1, 1.2.4]{DA20}, D'Addezio proved that the monodromy groups of semi-simple coefficient objects in an $E$-compatible system are $\lambda$-independent. 

In this paper, we will mainly consider the coefficient objects that arise from the first $l$-adic and the first crystalline cohomology of an abelian scheme $A_0$ over $X_0$. We will use $\mathrm{fpl}(E)$ to denote the set of finite places of $E$. Let $\mathcal{E}_{l,0}:=H^1_{l,\et}(A_0)[\frac{1}{l}]$ ($l\neq p$) and $\mathcal{E}_{p,0}:=H^1_{\cris}(A_0)[\frac{1}{p}]$. Then $\{\mathcal{E}_{l,0}\}_{l\in \mathrm{fpl}(\bQ)-\{p\}}$ is a $\bQ$-compatible system. 
Indeed, for a closed point $x_0$ in $X_0$, the characteristic polynomials $P(\mathcal{E}_{l,0,x_0}, t)$ are all equal to the characteristic polynomial $P({A_{0,x_0}},t)\in \bQ[t]$ defined by  $P({A_{0,x_0}},n)=\deg([n]_{A_{0,x_0}}-F_{A_{0,x_0}})$, where $F_{A_{0,x_0}}$ is the geometric Frobenius, see \cite[Theorem 12.18]{EGMAV}. For the place $p$, it is a bit more subtle. Indeed, the Frobenius that appears in  $H^1_{\cris}(A_0)[\frac{1}{p}]$ as an $F$-isocrystal is the absolute Frobenius. In order to fit $\mathcal{E}_{p,0}:=H^1_{\cris}(A_0)[\frac{1}{p}]$ in a compatible system with various $\mathcal{E}_{l,0}$'s, we need to consider $\mathcal{E}_{p,0}'$, the $F^{N}$-isocrystal associated to $H^1_{\cris}(A_0)[\frac{1}{p}]$, where $N=\log_p q$. Then $\{\mathcal{E}_{l,0}\}_{l\in \mathrm{fpl}(\bQ)-\{p\}}\cup \{\mathcal{E}_{p,0}'\}$ is a $\bQ$-compatible system. 

Nevertheless, we will still call the collection $\{\mathcal{E}_{u,0}\}_{u\in \mathrm{fpl}(\bQ)}$ a \textbf{weakly compatible system}. In general, a collection of coefficient objects $\{\mathcal{E}_{\lambda,0}\}_{\lambda\in \Sigma}$ will be termed a \textbf{weakly $E$-compatible system}, if for every $\lambda|p$, after replacing $\mathcal{E}_{\lambda,0}$ by an associated  $\mathcal{E}_{\lambda,0}'$ obtained by raising the Frobenius to a power, the new collection $\{\mathcal{E}_{\lambda,0}\}_{\lambda\in \Sigma,\lambda\nmid p}\cup \{\mathcal{E}'_{\lambda,0}\}_{\lambda\in\Sigma,\lambda|p}$ is an $E$-compatible system. 
\begin{lemma}[D'Addezio]\label{lm:compatible} If $\{\mathcal{E}_{\lambda,0}\}_{\lambda\in \Sigma}$ is a weakly compatible system of semi-simple coefficient objects over $X_0$, then the neutral components of the monodromy groups of $\mathcal{E}_{\lambda,0}$ are $\lambda$-independent in the sense of \cite[Theorem 1.2.1]{DA20}. Furthermore, there is a finite étale cover of $X_0'/X_0$, such that the monodromy groups of the elements in the pullback system $\{\mathcal{E}_{\lambda,0}|_{X'_0}\}_{\lambda\in \Sigma}$ are connected. 
\end{lemma}
\proof Let $\{\mathcal{E}_{\lambda,0}\}_{\lambda\in \Sigma,\lambda\nmid p}\cup \{\mathcal{E}'_{\lambda,0}\}_{\lambda\in\Sigma,\lambda|p}$ be the $E$-compatible system obtained from $\{\mathcal{E}_{\lambda,0}\}_{\lambda\in \Sigma}$ by replacing $\mathcal{E}_{\lambda,0}$ where $ \lambda|p$ with $\mathcal{E}'_{\lambda,0}$, as in the definition. Up to base change, the monodromy groups of $\mathcal{E}_{p,0}$ and $\mathcal{E}'_{p,0}$ have the same neutral component. So by \cite[Theorem 1.2.1]{DA20}, the neutral components of the monodromy groups of $\mathcal{E}_{\lambda,0}$ are $\lambda$-independent. For the second assertion, we know from \cite[Theorem 1.2.4]{DA20} that the $\pi_0$ of the monodromy groups for elements in $\{\mathcal{E}_{\lambda,0}\}_{\lambda\in \Sigma,\lambda\nmid p}$ are $\lambda$-independent. Since there are only finitely many $\lambda$ dividing $p$, there is a finite étale cover $X_0'/X_0$ satisfying the assertion made in the lemma. $\hfill\square$

\section{Constructions and conjectures}\label{sec:2}
This section is the general setup for the Tate-linear, the Mumford--Tate, and the André--Oort conjectures. We start from a rather general situation, defining the mod $p$ analogue of the Mumford--Tate group and reveal several properties.
Implications between the conjectures are also made. We then specialize to the case of products of GSpin Shimura varieties. 


\subsection{The construction of $\mathcal{X}_{f}$ and $\MTT(f)$}\label{subsubsec:MTgp} 
In the following, we will always assume that the hyperspecial level structure of
$\mathcal{A}_{g,\bK}$ is sufficiently small, and we will drop it from the notation. Let $\mathscr{A}_{g}$ be the canonical integral model. Suppose that $(X,x)$ is a smooth and connected pointed variety over $\Fpbar$ with a morphism $$f:(X,x)\rightarrow {\mathscr{A}}^{\ord}_{g,\Fpbar}.$$
\subsubsection{Minimal special subvarieties}\label{subsub:smallest} 
If $\mathcal{X}\subseteq \mathcal{A}_{g}$ is a special subvariety, 
Let $\mathrm{MS}_f(\mathcal{X})$ be the set of minimal irreducible special subvarieties  $\mathcal{Y}\subseteq\mathcal{X}_{\overline{\bQ}}$ such that $f$ factors through the naïve integral model of $\mathcal{Y}$. To distinguish it from other special subvarieties, an element in $\mathrm{MS}_f(\mathcal{X})$ is usually written as $\mathcal{X}_f$.

A \textbf{$p$-power Hecke correspondence} is a special correspondence $\mathcal{H}\subseteq \mathcal{A}_{g}\times \mathcal{A}_{g}$  (\S\ref{subsub:specialcorr}) parametrizing pairs of abelian schemes with an isogeny of $p$-power degree. Let $\mathscr{H}$ be the naïve integral model of $\mathcal{H}$ in $\mathscr{A}_{g}\times \mathscr{A}_{g}$, and let $\pi_1,\pi_2$ be the projections. A \textbf{$p$-power Hecke translate} of an irreducible subvariety $Y$ of $\mathscr{A}_{g}$ is an irreducible component of $\pi_1\pi_2^{-1}(Y)$.

\begin{lemma}\label{lm:minimality}
Any two elements in $\mathrm{MS}_f(\mathcal{A}_{g})$ are $p$-power Hecke translates of each other. 
\end{lemma}
\begin{proof}
Let $\mathcal{X}_f, \mathcal{Y}_f\in\mathrm{MS}_f(\mathcal{A}_{g})$, and let $\mathscr{X}_f, \mathscr{Y}_f$ be their naïve integral models. Let $\overline{W}$ be the integral closure of $W$ in $\overline{K}$ (see \S\ref{notation} for what they mean). By Theorem~\ref{thm:noottheorem}, the completions $\mathscr{X}^{/x}_{f,\overline{W}}$ and $\mathscr{Y}^{/x}_{f,\overline{W}}$ are all unions of torsion translates of formal subtori (see \S\ref{notation} for our convention on formal completion). Let $\mathcal{C}$ \textit{resp}. $\mathcal{D}$ be a component of $\mathscr{X}^{/x}_{f,\overline{W}}$  \textit{resp}. $\mathscr{Y}^{/x}_{f,\overline{W}}$ through which $f^{/x}$ factors. Now $\mathcal{C}$ and $\mathcal{D}$ both contain a quasi-canonical lift of $x$. Lemma~\ref{lm:uptoHecke} below implies that, there is a $p$-power Hecke translate $\mathcal{X}_f'$ of $\mathcal{X}_f$, such that ${\mathcal{X}'}^{/x}_{f,\overline{W}}$ contains an irreducible component $\cC'$, with the property that (1) $f^{/x}$ factors through $\cC'$ and (2) $\mathcal{C}'$, $\mathcal{D}$ contain the same quasi-canonical lift $\Tilde{x}_{\overline{W}}$. In this case, the intersection $\mathcal{C}'\cap \mathcal{D}$ is a formal subtorus translated by $\Tilde{x}_{\overline{W}}$, and $f^{/x}$ factors through $\mathcal{C}'\cap \mathcal{D}$. 

We claim that $\mathcal{Y}_f\subseteq \mathcal{X}_f'$. Let $\mathcal{Z}=\mathcal{X}_f' \cap \mathcal{Y}_f$, which is non-empty since $\tilde{x}_\bC\in \mathcal{Z}(\bC)$. If $\mathcal{Y}_f$ was not contained in $\mathcal{X}_f'$, then $\mathcal{Z}$ is a special subvariety with dimension 
 smaller than that of $\mathcal{Y}_f$. In order to deduce a contradiction, it suffices to show that $f$ factors through the naïve integral model $\mathscr{Z}$ of $\mathcal{Z}$. By Theorem~\ref{thm:noottheorem},  $\mathscr{Z}^{/x}_{\overline{W}}$ is a union of torsion translates of formal tori, and $\mathcal{C}'\cap \mathcal{D}$ is one of the components. Therefore $f^{/x}$ factors through  $\mathscr{Z}^{/x}_{\overline{W}}$, hence $f$ factors through $\mathscr{Z}$, contradicting the minimality of $\mathcal{Y}_f$. 

 Now we inverse the $p$-power Hecke translate to get a $\mathcal{Y}_f'\subseteq \mathcal{X}_f$, where $\mathcal{Y}_f'$ is the translate of $\mathcal{Y}_f$. Then $f$ factors through $ \mathscr{Y}'_{f,\overline{W}}$. By minimality of $\mathcal{X}_f$, we deduce that $\mathcal{Y}_f'=\mathcal{X}_f$ (and also $\mathcal{Y}_f=\mathcal{X}_f'$). 
\end{proof}

Therefore, we will often refer to $\mathcal{X}_f$ as ``\textit{the smallest special subvariety whose reduction contains the image of $f$}'', in the sense that it is the smallest such special subvariety up to $p$-power Heck translates.  

\begin{lemma}\label{lm:uptoHecke}
    Let $\mathcal{V}\subseteq\mathcal{A}_{g,\overline{\bQ}}$ be a special subvariety with naïve integral model $\mathscr{V}$, such that $\mathscr{V}^{\ord}_{\Fpbar}$ is nonempty. Let $v$ be an ordinary $\Fpbar$-point of $\mathscr{V}_{\Fpbar}$. Suppose that $\tilde{v}_{\overline{W}}$ and $\tilde{v}_{\overline{W}}'$ are two quasi-canonical lifts of $v$. Further suppose that $\mathscr{V}_{\overline{W}}^{/v}$ admits an irreducible component $\cC$ passing through $\tilde{v}_{\overline{W}}$. Then $\mathcal{V}$ admits a $p$-power Hecke translate ${\mathcal{V}}'$, such that 
 ${\mathscr{V}}'^{/x}_{\overline{W}}$ (the formal germ of the naïve integral model of ${\mathcal{V}}'$) admits a component $\mathcal{C}'$ passing through $\tilde{v}_{\overline{W}}'$, with the further property that $\mathcal{C}_{\Fpbar}=\mathcal{C}'_{\Fpbar}$. 
\end{lemma}
\begin{proof}
There is some $n$ such that the multiplication by $p^n$ map on $v$ lifts to an isogeny $\tilde{v}_{\overline{W}}'\rightarrow  \tilde{v}_{\overline{W}}$. Let $\mathcal{H}\subseteq \mathcal{A}_{g,\overline{\bQ}}\times \mathcal{A}_{g,\overline{\bQ}}$ be the $p$-power Heck correspondence parametrizing pairs with a degree $p^{2gn}$ isogeny. We claim that pull-push along this Hecke correspondence gives us the desired Hecke translate. 

 Let $\mathscr{H}$ be the naïve integral model of $\mathcal{H}
$. Then $\mathscr{H}_{\overline{W}}$ passes through the pair $(\tilde{v}_{\overline{W}}',\tilde{v}_{\overline{W}})$. Theorem~\ref{thm:noottheorem} implies that $\mathscr{H}_{\overline{W}}^{/(v,v)}$ is a union of torsion translates of formal subtori, and there is a unique component $\cM$ passing through $(\tilde{v}_{\overline{W}}',\tilde{v}_{\overline{W}})$. Let $\mathscr{D}_{(v,v)}([p^n])\in \mathscr{A}_{g,{W}}^{/v}\times \mathscr{A}_{g,{W}}^{/v}$ be the deformation space of pairs with an isogeny that lifts the endomorphism $[p^n]$ on $(v,v)$. Explicit computation using canonical coordinates (cf. Remark~\ref{subsub:canonicaldefendo}) tells us that $\mathscr{D}_{(v,v)}([p^n])_{\overline{W}}$ is a union of torsion translates of the diagonal torus  $\Delta\subseteq\mathscr{A}_{g,\overline{W}}^{/v}\times \mathscr{A}_{g,\overline{W}}^{/v}$. Now $\cM$ is contained in $\mathscr{D}_{(v,v)}([p^n])_{\overline{W}}$, and is identified with one of the irreducible components of $\mathscr{D}_{(v,v)}([p^n])_{\overline{W}}$ by dimension reasons. In particular, $\cM$ is a translate of $\Delta$ by $(\tilde{v}_{\overline{W}}',\tilde{v}_{\overline{W}})$, and $\cM_{\Fpbar}=\Delta_{\Fpbar}$. The lemma follows easily.
\end{proof}

\subsubsection{The group $\MTT(f)$}\label{subsub:thegroupMT}
Fix an $\mathcal{X}_f\in \mathrm{MS}_f(\mathcal{A}_{g})$. Define $\MTT(f)$ as the generic Mumford--Tate group of $\mathcal{X}_f$ with respect to the abelian scheme pulled back from $\mathcal{A}_{g}$. By Lemma~\ref{lm:minimality}, $\MTT(f)$ does not depend on the choice made. Let $\tilde{x}$ be the canonical lift of $x$, and let $\tilde{x}_{\bC}$ be the base change of  $\tilde{x}$ along $W\hookrightarrow \bC$ (\S\ref{notation}). Identify $H$ with $\bbH_{\mathrm{B},\tilde{x}_\bC}$ as $\bZ_{(p)}$-lattices. Let $\mu_{\tilde{x}_\bC}:\bG_m\rightarrow \GL(H_{\bQ})$ be the Hodge cocharacter of $\tilde{x}_\bC$.

Pick a quasi-canonical lift of $x$, whose base change to $\bC$ lies in $\mathcal{X}_f(\bC)$. Identify its rational 
Hodge structure with that of $\tilde{x}_{\bC}$. This allows us to identify $\MTT(f)$ as a subgroup of $\GL(H_{\bQ})$. The identification does not depend on the quasi-canonical lift chosen. 
In summary, we have the following chain of embeddings $$\MTT(\tilde{x}_{\bC})\hookrightarrow \MTT(f)\hookrightarrow \GL(H_\bQ).$$
Since $x$ is ordinary, it is classically known that $\MTT(\tilde{x}_\bC)_{\bQ_l}$ coincides with the neutral component of the $l$-adic étale (crystalline when $l=p$) monodromy of ${x}$.

\begin{lemma}\label{lm:endetaepsilon}
Let $\eta$ be the generic point of $X_0$, and let $\epsilon$ be the generic point of $\mathcal{X}_f$. Then under suitable identifications, we have $$\End_{\MTT(f)}(H_\bQ)=\End^0(A_{\overline{\epsilon}})=\End^0(A_{\overline{\eta}}).$$ 
\end{lemma}
\proof In the proof let's fix $\tilde{x}_{\overline{W}}\in \mathscr{X}_f(\overline{W})$ to be a quasi-canonical lift of $x$ and let $\tilde{x}_{\bC}$ be its base change to $\bC$ (the notation $\tilde{x}$ and $\tilde{x}_{\bC}$ is already taken for the canonical lift. But let's make a deliberate abuse of notation, since the rational Hodge structures are identified).

The first equality in the lemma is Torelli for abelian varieties. Let $W'$ be the ring of integers of a sufficiently large finite extension of $K$ so that we can form the base change $\mathscr{X}_{W'}$. Let $\overline{\epsilon}'$ be an algebraic closure of the generic point of $\mathscr{X}_{W'}$. Applying Lemma~\ref{lm:avext} to the following two situations: $(S,F)=(\mathscr{X}_{W'},\overline{\epsilon}')$ and $(S,F)=(X_0,\overline{\eta})$, we obtain
$\End^0(A_{\overline{\epsilon}})\subseteq \End^0(A_{\overline{\eta}})$. For the reversed inclusion, we first have  $\End^0(A_{\overline{\eta}})\subseteq \End^0(A_{x})=\End^0(A_{\tilde{x}_{\bC}})$. This means that $\End^0(A_{\overline{\eta}})$ can be lifted to, and identified as a subspace of $\End^0(A_{\tilde{x}_{\bC}})$. The strategy is to construct a special subvariety $\mathcal{Y}$ that is ``cut out by the Hodge cycles corresponding to $\End^0(A_{\overline{\eta}})$ '', then show that the naïve integral model of $\mathcal{Y}$ contains the image of $f$, and then conclude by the minimality of $\mathcal{X}_f$ that $\mathcal{X}_f\subseteq \mathcal{Y}$. And this implies that  $\End^0(A_{\overline{\eta}})\subseteq \End^0(A_{\overline{\epsilon}})$.

 We write $\{s_\alpha\}\subseteq H^{\otimes(1,1)}_{\bQ}$ for a set of endomorphisms that span $\End^0(A_{\overline{\eta}})$. There is a maximal irreducible special subvariety $\mathcal{Y}_{\overline{\bQ}}\subseteq \mathcal{A}_{g,\overline{\bQ}}$ admitting $\tilde{x}_{\mathbb{C}}$ as a geometric point, such that a suitable multiple of $\{s_\alpha\}$ extends globally to $\mathcal{Y}_{\overline{\bQ}}$. One construction of $\mathcal{Y}_{\overline{\bQ}}$ is as follows: Let $G\subseteq \GSp(H_\bQ)$ be the reductive group fixing $\{s_\alpha\}$. With suitable level structure, $G$ defines a Shimura subvariety of $\mathcal{A}_{g,\overline{\bQ}}$, which in turn gives rise to a special subvariety of $\mathcal{A}_{g,\overline{\bQ}}$. Let $\mathcal{Y}$ 
be a Hecke translate of an irreducible component such that $\mathcal{Y}(\bC)$ contains $\tilde{x}_\bC$. Then $\mathcal{Y}$ satisfies the desired property. In the following, let $\mathscr{Y}$ be the naïve integral model of $\mathcal{Y}$.
\begin{claim}
    $f$ factors through $\mathscr{Y}$.
\end{claim}It suffices to show that $f^{/x}$ factors through $\mathscr{Y}^{/x}_{\overline{W}}$. 
Let $\mathscr{D}_n:=\Def((A_x,\lambda_x),p^n\End(A_{\overline{\eta}})/W)\subseteq \mathscr{A}_{g,W}^{/x}$ be the deformation space parametrizing abelian schemes deforming the abelian variety with polarization $(A_x,\lambda_x)$ such that $p^n$-multiples of endomorphisms in $\End(A_{\overline{\eta}})$ also deform. Now deformation of endomorphisms can be easily understood via the canonical pairing: from the last paragraph of Remark~\ref{subsub:canonicaldefendo}, one finds that $\mathscr{D}_{n}$ is the product of a formal subtori of $\mathscr{A}_{g,{W}}^{/x}$ with a finite flat group scheme over $W$. When $n$ is sufficiently large, there is a unique irreducible component $\mathcal{D}\subseteq \mathscr{D}_{n,\overline{W}}$ that passes through $\tilde{x}_{\overline{W}}$. Then $f^{/x}$ factors through $\mathcal{D}$. Let $\cC$ be the irreducible component of $\mathscr{Y}^{/x}_{\overline{W}}$ that contains $\tilde{x}_{\overline{W}}$. Then $\cC$ factors through $\mathscr{D}_{n,\overline{W}}$ (since $p^n\End(A_{\overline{\eta}})$ deforms to $\mathscr{Y}^{/x}$ for $n\gg 0$), hence must factors through $\mathcal{D}$. Now by Theorem~\ref{thm:noottheorem}, $\dim \mathcal{C}=\dim \mathcal{Y}$, and the later equals the dimension of the opposite unipotent of the Hodge cocharacter $\mu_{\tilde{x}_{\bC}}$ in $G$.  Let $G_p$ be the subgroup of $\GSp(\omega(A_x[p^\infty])_{\bQ_p})$ fixing $\End(A_{\overline{\eta}})\otimes \bQ_p$. Then $\dim \cD$ equals the dimension of the opposite unipotent of the canonical Hodge cocharacter $\mu$ of $x$ in $G_p$. From (\ref{eq:comparison111}) and the paragraphs after that, we can identify $G\otimes \bQ_p$ with $G_p$ as subgroups of $\GL(\omega(A_x[p^\infty])_{\bQ_p})$  (since both groups are cut out by the cycles spanning $\End^0(A_{\overline{\eta}})$), and $\mu_{\tilde{x}_{\bC}}$ also gets identified with $\mu$. It follows that $\dim\cD=\dim \mathcal{C}$, hence $\cC=\cD$. Therefore $f^{/x}$ factors through $\cC$, hence through $\mathscr{Y}^{/x}_{\overline{W}}$, as desired.


We also give a short but not self-contained proof of the claim in the spirit of Noot. By Ogus (\cite[VI. \S4.1]{DP82}), $\{s_\alpha\}$ and the polarization give rise to a collection of Tate classes $\{t_\beta\}\in \Fil^0\bbH_{\mathrm{dR},\tilde{x}_{\overline{W}}}^{\otimes}$. Let $\cN$ be the local moduli of Tate classes $\{t_\beta\}$ passing through $\tilde{x}_{\overline{W}}$ as per \cite[Theorem 2.8]{N96}. Then $f^{/x}$ factors through $\cN$ (in fact, $\cN=\cD$). Since $\{t_\beta\}$ extends over $\cC$ (more precisely, extends over the $\Fil^0$ part of its relative de Rham cohomology), we have $\cC\subseteq\cN$.  It is shown in 
\cite[Theorem 3.7]{N96} using dimension arguments that $\cC=\cN$. This implies that $f^{/x}$ factors through $\cC$, hence through $\mathscr{Y}^{/x}_{\overline{W}}$, as desired.

Finally, let $\mathcal{Z}=\mathcal{X}_f\cap\mathcal{Y}$, which is non-empty since $\tilde{x}_\bC\in \mathcal{Z}(\bC)$. If $\mathcal{X}_f\not\subseteq\mathcal{Y}$, then we can use a similar argument as in Lemma~\ref{lm:minimality} to deduce a contradiction against the minimality of $\mathcal{X}_f$. Therefore we must have $\mathcal{X}_f\subseteq \mathcal{Y}$. This means $\MTT(f)\subseteq G$. Therefore, $\MTT(f)$ fixes $\End^0(A_{\overline{\eta}})$, hence $\End^0(A_{\overline{\eta}})\subseteq \End^0(A_{\overline{\epsilon}})$.$\hfill\square$

\begin{lemma}\label{lm:avext}
Suppose that $A$ is an abelian scheme over a Noetherian integral scheme $S$. Let $F$ be a field extension of $K(S)$ and $w\in \End(A_F)$. Then $w$ extends to a finite cover $S'/S$ with $K(S')\subseteq F$. If $S$ is normal, one can further require $S'/S$ to be étale.
   \end{lemma}
\proof By 
\cite[Proposition 2.3]{MR75}, every component of the group End-scheme $\mathcal{E}nd(A/S)$ is finite and unramified over $S$. 
Note that $w$ corresponds to an $F$-point of $\mathcal{E}nd(A/S)$ dominating $S$. We can let $S'$ be the irreducible component of $\mathcal{E}nd(A/S)$ containing the image of $F$. If $S$ is normal, then a irreducible component of  $\mathcal{E}nd(A/S)$ that dominants $S$ is flat over $S$ (cf. \cite[Theorem 18.10.1]{EGAIV}), hence étale.
$\hfill\square$

\subsection{ $\MTT(f)$ and monodromy groups}\label{subsub:idfi} 
Let $X_0/\bF_q$ be a geometrically connected smooth variety over a sufficiently large finite field with a map $f_0:(X_0,x)\rightarrow \mathscr{A}_{g,\mathbb{F}_q}^{\ord}$. Let $f: (X,x)\rightarrow \mathscr{A}_{g,\mathbb{F}}^{\ord}$ be the base change of $f_0$ to $\Fpbar$, as in the setup of \S\ref{subsubsec:MTgp}.

Let $\bH_{\bullet}$ be various local systems on $\mathscr{A}_g$ (\S\ref{Sec:GSpinDef}). To simultaneously handle $l$-adic and crystalline local systems, we use the following convention: let $u\in \mathrm{fpl}(\mathbb{Q})$, then
\begin{equation}\label{eq:Hlp}
\mathbb{H}_{u}:=\left\{\begin{aligned}
    & \mathbb{H}_{l,\et},\,\,&u=l, \\
    & \mathbb{H}_{\cris},\,\,&u=p.
\end{aligned}\right.
\end{equation}
The convention will be carried through the rest of this paper. We will also write $\bH_{u,f_0}$ (or $\bH_{u,X_0}$, when the notation is more convenient) for the pullback of $\bH_u$ along $f_0$, and write \begin{equation}\label{eq:G_uG_p}
    G_u(f):=G(\mathbb{H}_{u,f_0},x)^{\circ}=G(\mathbb{H}_{u,X_0},x)^{\circ}.
\end{equation}
when $\mathbb{H}_{u,f_0}$ and $x$ is clear from the context. This notation is reasonable, since $G(\mathbb{H}_{u,f_0},x)^{\circ}$ only depends on $f=f_{0,\Fpbar}$ but not a finite field model. We fix an $\mathcal{X}_f\in \mathrm{MS}_{f}(\mathcal{A}_g)$, and identify $H$ with $\bbH_{\mathrm{B},\tilde{x}_\bC}$ as in \S\ref{subsub:thegroupMT}. Let $u\in \mathrm{fpl}(\bQ)$. Recall from \S\ref{sec:monodromy} that  $\omega_x(\mathbb{H}_{u})$ is a $\bQ_u$-space. We have canonical identifications
\begin{equation}\label{eq:canidentifications}
    H_{\bQ_u}\simeq \omega_x(\mathbb{H}_{u}),
\end{equation}
which follows from the Betti--étale comparison when $u\neq p$, and from (\ref{eq:comparison111}) when $u=p$. So we can identify ${\MTT}(f)_{\bQ_u}$ as a subgroup of $\GL(\omega_x(\mathbb{H}_{u}))$ (when $u\neq p$ this is classic, when $u=p$ this is by the paragraphs after (\ref{eq:comparison111})).

Therefore, it makes sense to compare ${\MTT}(f)_{\bQ_u}$ and $G_u(f)$ as subgroups of $\GL(\omega_x(\mathbb{H}_{u}))$, or as subgroups of $\GL(H_{\bQ_u})$ by further composing the inverse of (\ref{eq:canidentifications}). We will freely regard ${\MTT}(f)_{\bQ_u}$ and $G_u(f)$ as subgroups of $\GL(\omega_x(\mathbb{H}_{u}))$, or as subgroups of $\GL(H_{\bQ_u})$, whenever it is more convenient.

\begin{proposition}
\label{lm:containedin}$G_u(f)\subseteq {\MTT}(f)_{\bQ_u}$.
\end{proposition}
\begin{proof}
We regard both groups as subgroups of $\GL(\omega_x(\mathbb{H}_{u}))$. Descend $\mathcal{X}_f$ to a sufficiently large number field $k$. We again let $\mathfrak{p}$ be the prime of $O_{k}$ induced from $\bC\simeq \bC_p$. Let $\Fpbar_q$ be the residue field of $O_{k,(\mathfrak{p})}$, $\mathcal{K}$ be the $\mathfrak{p}$-adic completion of $k$ and let $\mathcal{V}$ be its ring of integers. Let $\mathscr{X}_f$ be the naïve integral model over $O_{k,(\mathfrak{p})}$.  Let $\epsilon$ be the generic point of $\mathcal{X}_f$, and let $\overline{\epsilon}\rightarrow \epsilon$ be a geometric point. Let $\widetilde{\mathscr{X}}_f$ be the normalization of $\mathscr{X}_f$, then $\epsilon$ is also the generic point of $\widetilde{\mathscr{X}}_f$. 

We first treat the case where $u\neq p$. Let $\eta_0$ be the generic point of $X_0$ and let $\overline{\eta}_0\rightarrow \eta_0$ be a geometric point.  The point $ {\eta}_0\xrightarrow{f} \mathscr{X}_f$ lifts to $\eta_0' \rightarrow \widetilde{\mathscr{X}}_f$ for some finite cover $\eta_0'\rightarrow \eta_0$. By Lemma~\ref{lm:etcrisopensub}, we have $G(\mathbb{H}_{u,X_0},\overline{\eta}_0)^\circ=G(\mathbb{H}_{u,\eta'_0},\overline{\eta}_0)^\circ$. On the other hand, we also have the following relations between several étale fundamental groups: 
$$ {\pi_1^{\et}(\eta_0',\overline{\eta}_0)}\rightarrow \pi_1^{\et}(\widetilde{\mathscr{X}}_f,\overline{\eta}_0) \xleftarrow{\sim}\pi_1^{\et}(\widetilde{\mathscr{X}}_f,\overline{\epsilon})\xleftarrow{g} \pi_1^{\et}(\epsilon,\overline{\epsilon}).$$
The map $g$ is surjective since $\widetilde{\mathscr{X}}_f $ is normal (\cite[0BQM]{stacks-project}). Recall that the monodromy groups are the Zariski closures of the images of étale fundamental groups. Choosing appropriate identifications of the fibers, we see that $G_u(f)\subseteq G(\mathbb{H}_{u,X_0},x)^\circ\subseteq G(\mathbb{H}_{u,\epsilon},\overline{\epsilon})^\circ$. The later group is known to be contained in $\MTT(f)_{\bQ_u}$ (\cite[ I]{DP82}). 

Now we treat the crystalline case. Let $\{s_{\alpha}\}\subseteq H_\bQ^\otimes$ be a finite collection of tensors cutting out $\MTT(f)$. Via the process in \S\ref{subsub:pHodge}, an $s_{\alpha}$ gives rise to a Galois invariant tensor $s_{\alpha,p,\et,\tilde{x}}\in \bH_{p,\et,\tilde{x}}^{\otimes}$, an $F$-invariant tensor $s_{\alpha,\cris,x}\in \Fil^0\bH_{\cris,x}^\otimes$, and a tensor $\overline{s}_{\alpha,\cris,x}\subseteq \omega_x(\bH_{\cris})^\otimes$. Then (\ref{eq:comparison111}) takes $s_\alpha$ to $ s_{\alpha,p,\et,\tilde{x}}$ then to $\overline{s}_{\alpha,\cris,x}$. It suffices to show that $G_p(f)$ fixes $\{\overline{s}_{\alpha,\cris,x}\}$. 

By Theorem~\ref{thm:noottheorem} $\widetilde{\mathscr{X}}_{f}$ is smooth along the ordinary locus $Z_0:=\widetilde{\mathscr{X}}_{f,\Fpbar_q}^{\ord}$. Fix an $\Fpbar$-point of $Z_0$ that lifts $f_0(x)$ and by abuse of notation, still call it $x$. Define $G_p(Z):=G(\bH_{\cris,Z_0},x)^{\circ}$. Then $G_p(f)\hookrightarrow G_p(Z)$, so it suffices to show that $G_p(Z)$ fixes $\{\overline{s}_{\alpha,\cris,x}\}$.

Take the rigid generic fiber of the formal completion of $\mathscr{A}_{g,O_{k,(\mathfrak{p})}}$ along its special fiber, and call it $\mathfrak{A}$. {By \cite{BO83,Rigid_coh},} the overconvergent $F$-isocrystal $\bH_{\cris}\otimes \mathcal{K}$ is a bundle with connection over $\mathfrak{A}$ with Frobenius structure and overconvergence conditions, and there is a canonical isomorphism of bundles with connection between $\bH_{\mathrm{dR},\mathfrak{A}}$ and $\bH_{\cris}\otimes \mathcal{K}$. Take completion of $\widetilde{\mathscr{X}}_{f}$ along the special fiber, and take the rigid generic fiber of the resulting formal scheme. Then we take the rigid tube $]Z_0[$ of the ordinary locus (i.e., the admissible open subspace that specializes to the ordinary locus). Then $\bH_{\cris,Z_0}\otimes \mathcal{K}$ is the pullback of $\bH_{\cris}\otimes \mathcal{K}$ over $]Z_0[$ (see \cite[\S 4]{Rigid_coh} for the formation of overconvergent $F$-crystals over the rigid tube), and is isomorphic to
$\bH_{\mathrm{dR},{]Z_0[}}$ as a bundle with connection.

Via Betti--de Rham comparison, $\{s_\alpha\}$  {give rise to overconvergent horizontal sections $\{s_{\alpha,\mathrm{dR}}\}\subseteq \Fil^0(\bH_{\mathrm{dR},{]Z_0[}}\otimes\bC_p)^\otimes$.} The restriction of $s_{\alpha,\mathrm{dR}}$ to $(\bH_{\cris,x}\otimes\bC_p)^\otimes$ is exactly $s_{\alpha,\cris,x}\otimes \bC_p$. Since $s_{\alpha,\cris,x}$ is $F$-invariant, we see that $s_{\alpha,\mathrm{dR}}$ is also $F$-invariant under the Frobenius on $\bH_{\cris,Z_0}\otimes\bC_p$. As a result, $\{s_{\alpha,\mathrm{dR}}\}$ are $F$-invariant overconvergent horizontal sections of $(\bH_{\cris,Z_0}\otimes \bC_p)^\otimes$. Therefore, $G_p(Z)\otimes \bC_p$ fixes $ \{s_{\alpha,\cris,x}\}\otimes \bC_p$. This immediately implies that  $G_p(Z)\otimes W$ fixes $ \{s_{\alpha,\cris,x}\}$, hence $G_p(Z)$ fixes $ \{\overline{s}_{\alpha,\cris,x}\}$. The proof is complete. 
\end{proof}
\begin{lemma}\label{lm:etcrisopensub}The following are true: \begin{enumerate}
    \item\label{it:1etcrisopensub} Suppose that $(X'_0,x)$ is a geometrically connected smooth pointed variety over a finite field that 
    admits a dominant morphism $X'_0\rightarrow X_0$. Then $G(\mathbb{H}_{u,X_0'},x)^\circ=G(\mathbb{H}_{u,X_0},x)^\circ$. 
    \item\label{it:2etcrisopensub} Suppose that $u\neq p$. Let $\eta_0$ be the generic point of $X_0$ and let $\xi$ the spectrum of a finitely generated field with a morphism $\xi\rightarrow \eta_0$, then $G(\mathbb{H}_{u,\xi},\overline{\xi})^\circ=G(\mathbb{H}_{u,\eta_0},\overline{\xi})^\circ=G(\mathbb{H}_{u,X_0},\overline{\xi})^\circ$.
\end{enumerate}
\end{lemma}
\begin{proof}
Using Theorem~\ref{lm:compatible}, we can reduce the $u=p$ case of (\ref{it:1etcrisopensub}) to the $u\neq p$ case, which can be further reduced to (\ref{it:2etcrisopensub}). Recall that the monodromy groups are the Zariski closures of the images of étale fundamental groups. The first equality of (\ref{it:2etcrisopensub}) follows since $\pi_1^{\et}(\mu,\overline{\xi})\rightarrow \pi_1^{\et}(\eta_0,\overline{\xi})$ is of finite index (if $\xi$ is purely transcendental, then the map is surjective; if $\xi$ is algebraic, then it is of finite index. The general case is a composition of the two). The second equality follows from the fact that $\pi_1^{\et}(\eta_0,\overline{\xi})\twoheadrightarrow \pi_1^{\et}(X_0,\overline{\xi})$, since $X_0$ is normal (\cite[0BQM]{stacks-project}).\end{proof}

\begin{proposition}\label{prop:compareEnd} Let $\eta$ be the generic point of $X_0$, then
$$\End(A_{\overline{\eta}})_{\bQ_u}=\End_{G_u(f)}(H_{\bQ_u})=\End_{\MTT(f)_{\bQ_u}}(H_{\bQ_u}).$$
\end{proposition}
\begin{proof}
    This follows from Lemma~\ref{lm:endetaepsilon}, Lemma~\ref{lm:avext}, and Tate's conjecture over function fields.
\end{proof}

Now we establish a characteristic $p$ analogue of \cite[Theorem 1.3.1]{V08}. 
In the following, for an algebraic group $G$, we denote by $Z^0G$ the connected center of $G$. 
\begin{proposition}\label{prop:centercontain}
$Z^0G_u(f)= Z^0{\MTT}(f)_{\bQ_u}$.
\end{proposition}
\proof Note that $G_u(f)$ and $\MTT(f)_{\bQ_u}$ both embed into $\End(H_{\bQ_u})$. Combining Proposition~\ref{lm:containedin} and Proposition~\ref{prop:compareEnd}, we find that
\begin{align*}
ZG_u(f)= G_u(f)\cap \End_{G_u(f)}(H_{\bQ_u})\subseteq \MTT(f)_{\bQ_u}\cap \End_{\MTT(f)_{\bQ_u}}(H_{\bQ_u})=Z{\MTT}(f)_{\bQ_u}.
\end{align*}
This shows one direction.

Now by \cite[Theorem 4.2.11]{DA20}, there is a point $y$ of $X$ such that the Frobenius torus of $y$ identifies with a maximal torus of $G_u(f)$. Denote by $M\subseteq \MTT(f)$ the group $\MTT(\bH_{\mathrm{B},\tilde{y}_\bC})$, here $\tilde{y}$ is a quasi-canonical lift of $y$ with $\tilde{y}_{\bC}\in \mathcal{X}_f $. Then $M_{\bQ_u}$ coincides with the Frobenius torus of $y$. We have 
$Z^0G_u(f)\subseteq M_{\bQ_u}\subseteq G_u(f)\subseteq \MTT(f)_{\bQ_u}$. As a result, $M_{\bQ_u}$ is an almost direct product:$$M_{\bQ_u}=Z^0G_u(f)\cdot T_u,$$ where $T_u$ is a torus in $G_u(f)^{\der}$. On the other hand, the abelianization $\MTT(f)^{\mathrm{ab}}$ corresponds to a variation of CM Hodge structures over $\mathcal{X}_f $, which is constant. Therefore 
$M$ surjects onto $\MTT(f)^{\mathrm{ab}}$. Since $Z^0 G_u(f)\subseteq Z^0{\MTT}(f)_{\bQ_u}$ and $T_u\subseteq G_u(f)^{\der}\subseteq {\MTT}(f)_{\bQ_u}^{\der}$, we see that $Z^0 G_u(f)$ also surjects onto $\MTT(f)^{\mathrm{ab}}_{\bQ_u}$. This forces $Z^0 G_u(f)=Z^0{\MTT}(f)_{\bQ_u}$. $\hfill\square$

\subsection{Statement of the conjectures and the first reduction}\label{subsub:tmIMPLICA}  Let $X_0/\bF_q$ be a geometrically connected smooth variety over a sufficiently large finite field with a map $f_0:(X_0,x)\rightarrow \mathscr{A}_{g,\mathbb{F}_q}^{\ord}$. Let $f: (X,x)\rightarrow \mathscr{A}_{g,\mathbb{F}}^{\ord}$ be the base change of $f_0$ to $\Fpbar$. Let $\mathscr{T}_{f,x}$ be the smallest subtorus of the Serre--Tate formal torus $\mathscr{A}_{g, \Fpbar}^{/x}$ containing the image of $f^{/x}$. Let  \begin{equation}\label{eq:ttlattices}
    \mathcal{T}_{f,x}\subseteq X_*(\Def(A_x[p^\infty]/\Fpbar))=\Lie U_{\GL,\mu^{-1}}
\end{equation} be the cocharacter lattice of $\mathscr{T}_{f,x}$ ($\mu$ being the canonical Hodge cocharacter of $x$). Let $$T_{f,x}=\mathcal{T}_{f,x}\otimes \bQ_p\subseteq \Lie U_{\GL,\mu^{-1},\bQ_p}.$$


Using the language introduced in the previous sections, we rephrase the conjectures that we stated in the introduction.
\begin{conj}[$=$ Conjecture~\ref{conj:MTforAg}]\label{conj:MT} For every $u\in \mathrm{fpl}(\bQ)$,  $G_u(f)={\MTT}(f)_{\bQ_u}$.
\end{conj}
\begin{conj}[$=$ Conjecture~\ref{Conj:AOforSV}]\label{conj:AOAO}
Suppose $X\subseteq \mathscr{A}_{g,{\Fpbar}}^{\ord}$ is a closed subvariety which contains a Zariski dense collection $\Xi$ of positive dimensional special subvarieties, then $X$ is quasi-weakly special.     
\end{conj}

\begin{conj}[$=$ Conjecture~\ref{conj:AxSchanuel}]\label{conj:Ttlinear}
Let $\mathcal{X}_f\in \mathrm{MS}_f(\mathcal{A}_g)$, then $\dim \mathcal{X}_{f}=\dim T_{f,x}$.
\end{conj}
To see why this is equivalent to Conjecture~\ref{conj:AxSchanuel}, just note that Conjecture~\ref{conj:AxSchanuel} is equivalent to the statement that $\mathscr{X}_{f,\Fpbar}^{/x}$ admits $\mathscr{T}_{f,x}$ as an irreducible component, which is equivalent to $\dim \mathcal{X}_{f}=\dim T_{f,x}$ by Theorem~\ref{thm:noottheorem}.

\begin{remark}\label{rmk:ASimpliesTtl}
 Let $f:X\hookrightarrow {\mathscr{A}}^{\ord}_{g,\Fpbar}$ be an immersion. Then $X$ is special if and only if for one (hence any) $\mathcal{X}_f\in \mathrm{MS}_f(\mathcal{A}_g)$, we have $\dim X=\dim \mathcal{X}_f$. Therefore   Conjecture~\ref{conj:TTl} is a special case of Conjecture~\ref{conj:Ttlinear} when $f$ is an immersion. 
\end{remark}

\begin{theorem}\label{Thm:Tatelocal}
Notation as above.  Let $\bH_{p,X_0}^-$ be the underlying $F$-isocrystal of $\bH_{p,X_0}$. Then 
\begin{enumerate}
    \item\label{itcor1} 
$G(\bH_{p,X_0},x)$ and $G(\gr\bH_{p,X_0},x)$ are reductive, and $G(\bH_{p,X_0}^-,x)$ is the opposite parabolic of $G(\bH_{p,X_0},x)$ with respect to $\mu$.  
\item\label{itcor2}$G(\bH_{p,X_0}^-,x)=U(\bH_{p,X_0}^-,x)\rtimes G(\gr\bH_{p,X_0},x)$, and $U(\bH_{p,X_0}^-,x)$ is the opposite unipotent of $G(\bH_{p,X_0},x)$ with respect to $\mu$.
\item\label{itcor3} $U(\bH_{p,X_0}^-,x)=T_{f,x}$, where $T_{f,x}$ is viewed as a subgroup of $U_{\GL,\mu^{-1},\bQ_p}$ via exponential map.
\end{enumerate}
\end{theorem}
\proof Points (\ref{itcor1}) and (\ref{itcor2}) follow from Lemma~\ref{strofmon} and Theorem~\ref{MDA}. Point (\ref{itcor3}) follows from Theorem~\ref{thm:etale-glob} if we can show that $T_{f,x}=N(\mathscr{G}^{/x})$, where $\mathscr{G}$ is the pullback of $A[p^{\infty}]$ to $X$. Following \S\ref{sub:ST}, let $\Lambda=\Lie U_{\GL,{\mu}}$, and let $X^{/x}=\Spf R$. Consider the Serre--Tate pairing $q\in \Hom(\Lambda, \bG_m^{\wedge}(R))$ corresponding to $\mathscr{G}^{/x}$. Since $R$ is reduced, $\ker(q)$ is saturated in $\Lambda$. Let $\ker(q)^{\perp}$ be the sub-lattice of $\Lambda^{\vee}$ that pairs to 0 with $\ker(q)$. Then $\ker(q)^{\perp}\otimes \bG_m^\wedge$ is the smallest subtorus of $\Lambda^{\vee}\otimes \bG_m^\wedge=\mathscr{A}_{g,\Fpbar}^{/x}$ through which $f^{/x}$ factors. This shows that $X_*(\mathscr{T}_{f,x})=\ker(q)^{\perp}$. $\hfill\square$ 
\begin{proof}[Proof of Proposition~\ref{prop:cocharofSTspecial}]
It suffices to show that $\mathcal{C}_\Fpbar$ is the formal subtorus of $\Def(A_x[p^{\infty}]/\Fpbar)$ with rational cocharacter $\Lie U_{\MTT_{\bQ_p},\mu^{-1}}$. Let $f:(X,x)\rightarrow{\mathscr{X}}^{\ord}_{\Fpbar}$ be a morphism from a smooth $\Fpbar$-variety $X$, such that $\mathcal{C}_\Fpbar=\mathscr{T}_{f,x}$ (e.g., one can take $X$ to be a suitable component of $\widetilde{\mathscr{X}}^{\ord}_{\Fpbar}$ by Theorem~\ref{thm:noottheorem}). Then $\MTT(\mathcal{X}_\bC)=\MTT(f)$. We then combine 
Proposition~\ref{lm:containedin} with Theorem~\ref{Thm:Tatelocal}(\ref{itcor3}) to show that $T_{f,x}\subseteq \Lie U_{\MTT_{\bQ_p},\mu^{-1}}$. Now Theorem~\ref{thm:noottheorem}
implies that $T_{f,x}= \Lie U_{\MTT_{\bQ_p},\mu^{-1}}$.
\end{proof}

\begin{proposition}[$=$ Theorem~\ref{thm:implications}]\label{prop:MTimpliesTl}
Conjecture~\ref{conj:MT} $\Rightarrow$     Conjecture~\ref{conj:Ttlinear}.
\end{proposition}
\begin{proof}
 Let $d=\dim \mathscr{T}_{f,x}$. By Theorem~\ref{Thm:Tatelocal}, $d$ equals the dimension of the opposite unipotent of $G_p(f)$ with respect to $\mu$. If Conjecture~\ref{conj:MT} is true, then $G_p(f)=\MTT(f)_{\bQ_p}$. Under identification (\ref{eq:comparison111}), we see that 
 $d$ equals the dimension of the opposite unipotent of $\MTT(f)$ with respect to $\mu_{\tilde{x}_\bC}$. This shows that $d=\dim \mathcal{X}_f$. So Conjecture~\ref{conj:Ttlinear} holds.
\end{proof}
\subsection{Intrinsic perspective and modifications}\label{sub:inpmM} Let the setup be the same as \S\ref{subsubsec:MTgp}. The goal of this section is to show that the construction of $\mathcal{X}_f$ and $\MTT(f)$ are essentially ``in dependent of'' the ambient space $\mathscr{A}_{g}$, and use this observation to simplify the conjectures.

\subsubsection{Changing ambient spaces}\label{subsub:cas} Let $\{a,b\}$ be an index set. In the following, suppose that: 
\begin{enumerate}
    \item\label{indpXM1} There are two ambient spaces $\mathcal{A}_{g_a}$ and $\mathcal{A}_{g_b}$, with canonical integral models $\mathscr{A}_{g_a}$ and $\mathscr{A}_{g_b}$.
    \item\label{indpXM2} There are special subvarieties $\mathcal{X}_a\subseteq \mathcal{A}_{g_a,\overline{\bQ}}$ and $\mathcal{X}_b\subseteq \mathcal{A}_{g_b,\overline{\bQ}}$, with naïve integral models $\mathscr{X}_a$ and $\mathscr{X}_b$, 
    whose reductions admit non-empty ordinary loci. 
 \item For $\alpha\in \{a,b\}$, there are connected smooth varieties $X_\alpha$ together with morphisms $f_\alpha:X_\alpha\rightarrow \mathscr{A}_{g_\alpha,\Fpbar}^{\ord}$ that factors through $\mathscr{X}_\alpha^{\ord}$.
\end{enumerate}
\begin{defn}\label{def:modifications}
    In the above situation, we will say that $f_b$ is a \textbf{modification of $f_a$}, if moreover: \begin{enumerate}\setcounter{enumi}{3}
      \item\label{indpXM4} There is a special correspondence $\mathcal{X}_{c}$ between $ \mathcal{X}_a$ and $\mathcal{X}_b$ with naïve integral model $\mathscr{X}_{c}\subseteq \mathscr{X}_a\times \mathscr{X}_b$ and projections $\pi_a: \mathscr{X}_{c}\rightarrow\mathscr{X}_a$ and $\pi_b: \mathscr{X}_{c}\rightarrow\mathscr{X}_b$ \footnote{Note that $\pi_\alpha$ induces a quasi-finite surjection on ordinary loci by Lemma~\ref{lm:extensiontonormal}.}.
       \item\label{indpXM3}  
    There is a connected smooth variety $X_c$ together with a morphism $f_c:X_c\rightarrow \mathscr{A}_{g_a,\Fpbar}^{\ord}\times \mathscr{A}_{g_b,\Fpbar}^{\ord}$ that factors through $\mathscr{X}^{\ord}_{c,\Fpbar}$, such that there are dominant morphisms $g_\alpha:X_c\rightarrow X_\alpha$ fitting into the following commutative diagram \begin{center}
        \begin{tikzcd}
{X_a} \arrow[d, "f_a"]    & {X_c} \arrow[l, "g_a"'] \arrow[r, "g_b"] \arrow[d, "f_c"]                       & {X_b} \arrow[d, "f_b"]    \\
{\mathscr{X}^{\ord}_{a,\Fpbar}} & {\mathscr{X}^{\ord}_{c,\Fpbar}} \arrow[r, "\pi_b"] \arrow[l, "\pi_a"'] & {\mathscr{X}^{\ord}_{b,\Fpbar}}
\end{tikzcd}
    \end{center}
    \end{enumerate}
    If we need to mention $f_c$, we will say that $f_b$ is a \textbf{modification of $f_a$ via $f_c$}.
\end{defn}

A modification can be understood as changing ambient space from $\mathcal{A}_{g_a}$ to $\mathcal{A}_{g_b}$. It is clear that modification is an equivalence relation. Let $f_b$ be a modification of $f_a$ via $f_c$ as above. To study how the conjectures behave under modifications, we need to further introduce the following:
  \begin{enumerate}\setcounter{enumi}{5}
    \item\label{indpXM5} For $\alpha\in \{a,b,c\}$, let  $\mathrm{MS}_{f_\alpha}(\mathcal{X}_{\alpha})$ be the set of 
    minimal special subvarieties of $\mathcal{X}_{\alpha}$ whose reduction contain the image of $f_\alpha$. 
     \item There are universal abelian schemes $A_a$ and $A_b$ over $\mathscr{A}_{g_a}$ and $\mathscr{A}_{g_b}$, and various sheaves $\bH_{a,\bullet}$ and $\bH_{b,\bullet}$ that arise from various relative cohomology of the universal abelian schemes. We equip $\mathscr{X}_{c}$ with the universal family pulled back from $\mathscr{A}_{g_a}\times \mathscr{A}_{g_b}$, as well as cohomology sheaves $\bH_{c,\bullet}=\bH_{a,\bullet}\oplus \bH_{b,\bullet}$.  Let $\MTT(f_\alpha)$ be the generic Mumford--Tate groups of an element in $\mathrm{MS}_{f_\alpha}(\mathcal{X}_{\alpha})$ with respect to the sheaf $\bH_{\alpha,\mathrm{B}}$. The group is well defined by 
    Lemma~\ref{lm:minimality}. 
\item\label{indpXM6.5} We put a marked point $x_c\in X_c(\Fpbar)$, which induces, via the diagram in (\ref{indpXM3}), marked points $x_a,x_b$ on $X_a,X_b$. This makes the diagram in (\ref{indpXM3}) commute in the category of pointed $\Fpbar$-schemes.
    \item\label{indpXM7} For $\alpha\in \{a,b,c\}$, let $\mathscr{T}_{f_\alpha,x_\alpha}$ be the smallest formal subtorus of  $\mathscr{A}_{g_\alpha,\Fpbar}^{/x_\alpha}$ containing the image of $f^{/x_\alpha}_\alpha$. Let $\mathscr{T}_{f_c,x_c}$ be the smallest subtorus of $\mathscr{A}_{g_a,\Fpbar}^{/x_a}\times \mathscr{A}_{g_b,\Fpbar}^{/x_b}$ containing the image of $f^{/x_c}_c$.
\end{enumerate}
The following lemma justifies our claim at beginning that $\mathcal{X}_f$ and $\MTT(f)$ are ``in dependent of'' the ambient space.
\begin{lemma}\label{lm:finitefinite}Notation as above. Let $f_b$ be a {modification of $f_a$ via $f_c$}.
\begin{enumerate}
    \item\label{modi:1} For $\alpha\in \{a,b\}$, $\pi_\alpha$ induces a finite surjective homomorphism $\mathscr{T}_{f_c,x_c}\rightarrow \mathscr{T}_{f_\alpha,x_\alpha}$.
    \item\label{modi:2} There is an $\mathcal{X}_{f_\alpha}\in \mathrm{MS}_{f_\alpha}(\mathcal{X}_{\alpha})$ for each $\alpha\in \{a,b,c\}$ such that 
 ${\pi}_a,\pi_b$ induce finite surjections $ \mathcal{X}_{f_a}\twoheadleftarrow\mathcal{X}_{f_c}\twoheadrightarrow \mathcal{X}_{f_b}$. 
 \item  $\pi_a,\pi_b$ induce surjctions $ \MTT(f_a)\twoheadleftarrow\MTT(f_c)\twoheadrightarrow\MTT(f_b)$ which become isomorphisms after taking adjoints.  
\end{enumerate}
\end{lemma}
\begin{proof} 
\begin{enumerate}
    \item Let $\alpha\in \{1,2\}$. The fact that $\mathscr{T}_{f_c,x_c}\rightarrow \mathscr{A}_{g_\alpha,\Fpbar}^{/x_\alpha}$ is a finite homomorphism follows from the proof of Lemma~\ref{lm:extensiontonormal}(\ref{lm:extnaiveit00}) (note that, after modulo $p$, the finite linear map $\pi_\alpha^{/x_c}:\mathscr{X}_{c,\overline{W}}^{/x_c}\rightarrow \mathscr{X}_{\alpha,\overline{W}}^{/x_\alpha}$ induces a finite homomorphism over each irreducible component). It remains to show that the image is exactly $\mathscr{T}_{f_\alpha,x_\alpha}$. But this is easy from the construction.
    \item  Pick an $\mathcal{X}_{f_c}\in \mathrm{MS}_{f_c}(\mathcal{X}_{c})$. Note that for $\alpha\in \{a,b\}$,  $\pi_\alpha(\mathcal{X}_{f_c})$ is a special subvariety of $\mathcal{X}_\alpha$ whose reduction contains $f_\alpha(X_\alpha)$. Hence there exists $\mathcal{X}_{f_\alpha}\in \mathrm{MS}_{f_\alpha}(\mathcal{X}_{\alpha})$ such that $\pi_\alpha(\mathcal{X}_{f_c})\supseteq\mathcal{X}_{f_\alpha}$. It suffices to show that $\pi_\alpha(\mathcal{X}_{f_c})=\mathcal{X}_{f_\alpha}$. Once this is done, the finiteness is clear.

In the following we can fix an $\alpha\in \{a,b\}$.  Since $\pi_\alpha(\mathcal{X}_{f_c})\supseteq\mathcal{X}_{f_\alpha}$, we can take preimage of $\mathcal{X}_{f_\alpha}$ in $\mathcal{X}_{f_c}$. Let $C_1,...,C_N\subseteq \mathcal{X}_{f_c}$ be the irreducible components of the preimage. They are all special subvarieties. For $1\leq i\leq N$, let $\overline{C}_i$ be the closure of $C_i$ in $\mathscr{X}_{f_c}^{\ord}$. Our goal is to show that  $\exists i$ such that $f_c$ factors through $\overline{C}_i$. 

To show this, let $z$ be an $\Fpbar$-point in the image of $f_c$, and let $y=\pi_\alpha(z)$. Applying Theorem~\ref{thm:noottheorem} to  $\mathscr{X}_{f_c,\overline{W}}^{/z}$, and using the linearity of $\pi_\alpha^{/z}:\mathscr{X}_{c,\overline{W}}^{/z}\rightarrow \mathscr{X}_{\alpha,\overline{W}}^{/y}$ as observed in the proof of Lemma~\ref{lm:extensiontonormal}(\ref{lm:extnaiveit00}, we can find an irreducible component $\mathcal{C}\subseteq\mathscr{X}_{f_c,\overline{W}}^{/z}$ whose image under $\pi_\alpha^{/z}$ contains an irreducible component $\mathcal{D}\subseteq \mathscr{X}_{f_\alpha,\overline{W}}^{/y}$, as well as a quasi-canonical lift $\tilde{z}_{\overline{W}}\in \mathcal{C}(\overline{W})$ of $z$, whose image under $\pi_\alpha^{/z}$ is a quasi-canonical lift $\tilde{y}_{\overline{W}}\in \mathcal{D}(\overline{W})$ of $y$. It then follows that $\tilde{z}_{\overline{K}}\in \mathcal{X}_{f_c}(\overline{K})$ is a preimage of $\tilde{y}_{\overline{K}}\in  \mathcal{X}_{f_\alpha}(\overline{K})$ that specializes to $z$. In particular, $\exists j$ such that $\tilde{z}_{\overline{K}}\in C_j(\overline{K})$ and $z\in \overline{C}_j(\Fpbar)$.  Let $z$ vary over all $\Fpbar$-points of the image of $f_c$, we find that $f_c(X_c)$ is contained in the union of $\overline{C}_1,...,\overline{C}_N$, hence contained in some $\overline{C}_i$.  

Since ${C}_i$ is a special subvariety of $\mathcal{X}_{f_c}$, and $\mathcal{X}_{f_c}$ is a minimal subvariety whose naïve integral model contains the image of $f_c$, we conclude that $C_i=\mathcal{X}_{f_c}$. This shows that $\pi_\alpha(\mathcal{X}_{f_c})=\mathcal{X}_{f_\alpha}$.
    \item Surjectivity follows from the definition of $\MTT(f_c)$ using $\bH_{c,\mathrm{B}}=\bH_{a,\mathrm{B}}\oplus \bH_{b,\mathrm{B}}$. The isomorphism after taking adjoint follows from (\ref{modi:2}).
\end{enumerate} \end{proof}

\begin{proposition}\label{prop:redalongquasilifts}
Notation as above. If $f_a$ is a modification of $f_b$, then \begin{enumerate}
    \item Conjecture~\ref{conj:MT} is true for $f_a$ if and only if it is true for $f_b$,
    \item Conjecture~\ref{conj:Ttlinear} is true for $f_a$ if and only if it is true for $f_b$. 
\end{enumerate}
\end{proposition}
\begin{proof} 
\begin{enumerate}
    \item Let $f_c$ be as in Definition~\ref{def:modifications}. It suffices to show that Conjecture~\ref{conj:MT} is true for $f_c$ if and only if it is true for $f_a$. By Lemma~\ref{lm:finitefinite}, $\pi_a$ induces a surjection \begin{equation}\label{eq:mttred1}
        \MTT(f_c)\twoheadrightarrow \MTT(f_a)
    \end{equation} which becomes an isomorphism after taking adjoint. For $\alpha\in\{a,c\}$, let $X_{\alpha,0},f_{\alpha,0}$ be a finite field model of $X_\alpha,f_\alpha$. Let $u\in \mathrm{fpl}(\bQ)$. The fact $\bH_{c,u}=\bH_{a,u}\oplus \bH_{b,u}$ together with Lemma~\ref{lm:etcrisopensub} imply that $\pi_\alpha$ also induces a surjection \begin{equation}\label{eq:mttred2}
        G_{u}(f_c)\twoheadrightarrow G_{u}(f_a).
    \end{equation} By Proposition~\ref{lm:containedin}, (\ref{eq:mttred2}) also injects into the base change of (\ref{eq:mttred1}) to $\bQ_u$. If Conjecture~\ref{conj:MT} holds for $f_c$, then $\MTT(f_c)_{\bQ_u}=G_{u}(f_c)$, which implies that  $\MTT(f_a)_{\bQ_u}=G_{u}(f_a)$, so Conjecture~\ref{conj:MT}  holds for $f_a$. Conversely, if Conjecture~\ref{conj:MT} holds for $f_a$, then by passing to adjoints and counting dimension, we have $\MTT(f_c)_{\bQ_u}^{\ad}=G_{u}(f_c)^{\ad}$. By Proposition~\ref{prop:centercontain}, Conjecture~\ref{conj:MT}  holds for $f_c$. 
    \item Let $f_c$ be as in Definition~\ref{def:modifications}. It suffices to show that Conjecture~\ref{conj:Ttlinear} is true for $f_c$ if and only if it is true for $f_a$. By Lemma~\ref{lm:finitefinite}(\ref{modi:1}), we have $\rk\mathscr{T}_{f_c,x_c}=\rk\mathscr{T}_{f_a,x_a}$. 
By Lemma~\ref{lm:finitefinite}(\ref{modi:2}), we have $\dim\mathcal{X}_{f_a}=\dim\mathcal{X}_{f_c}$. Therefore $\dim\mathcal{X}_{f_c}=\rk\mathscr{T}_{f_c,x_c}$ if and only if $\dim\mathcal{X}_{f_a}=\rk\mathscr{T}_{f_a,x_a}$.
\end{enumerate}
\end{proof}

\subsubsection{Morphisms that separate factors}\label{subsub:spefactor}
Let $\bbI$ be a finite index set. For each $i\in \bbI$, let $\mathcal{A}_{g_i}$ be a Siegel modular variety and let $\mathscr{A}_{g_i}$ be its canonical integral model. We use $\mathcal{A}_{g_\bbI}$ to denote the product $\prod_{i\in \bbI}\mathcal{A}_{g_i}$, regarded as a Shimura subvariety of an ambient space $\mathcal{A}_{g}$. Consider a map $$f:X\rightarrow \mathscr{A}_{g, \Fpbar}^\ord$$
that factors through $\mathscr{A}_{g_\bbI, \Fpbar}^\ord$. Let $f_i$ be the composition $X\rightarrow \mathscr{A}_{g_\bbI, \Fpbar}^\ord\rightarrow \mathscr{A}_{g_i, \Fpbar}^{\ord}$. 
\begin{defn}
    We say that $f$ \textbf{separates factors over $\mathscr{A}_{g_\bbI, \Fpbar}^{\ord}$}, or simply \textbf{separates factors}, if for $i\in\bbI$ and $\mathcal{X}_{f_i}\in \mathrm{MS}_{f_i}(\mathcal{A}_{g_i})$: (1) $\mathcal{X}_{f_i}$ is almost simple (i.e. $\MTT(f_i)^{\ad}$ is simple), and (2) $\prod_{i\in \bbI}\mathcal{X}_{f_i}\in\mathrm{MS}_f(\mathcal{A}_{g_\bbI})$ (i.e., $\prod_{i\in \bbI}\mathcal{X}_{f_i}$ is a minimal special subvariety of $\mathcal{A}_{g_\bbI,\overline{\bQ}}$ whose naïve integral model contains the image of $f$).
\end{defn}
If $f$ {separates factors over $\mathscr{A}_{g_\bbI, \Fpbar}^{\ord}$}, then $\MTT(f)^{\ad}$ splits into a product of simple factors $\MTT(f_i)^{\ad}$, $i\in \bbI$. This is where the name ``separating factors'' come from. The following construction says that every $f:X\rightarrow \mathscr{A}_{g, \Fpbar}^\ord$ can be modified to one that separates factors:
\begin{construction}\label{constr:sepa}
Let $\mathcal{A}_{g_a}$ be a Siegel modular variety and let $f_a:X_a\rightarrow \mathscr{A}_{g_a,\Fpbar}^{\ord}$. We will construct a modification $f_b:X_b\rightarrow \mathscr{A}_{g_b,\Fpbar}^{\ord}$ with the following two properties:
\begin{itemize}
    \item $f_b$ is an immersion. 
    \item $f_b$ separates factors over some $\mathscr{A}_{g_\bbI, \Fpbar}^{\ord}\subseteq \mathscr{A}_{g_b,\Fpbar}^{\ord}$ for some finite index set $\bbI$. 
\end{itemize}

More precisely, we will construct $\mathcal{A}_{g_b}$, $\mathcal{X}_a,\mathcal{X}_b,\mathcal{X}_c,f_b,f_c$ (together with the given $\mathcal{A}_{g_a}$ and $f_a$) that fit in \S\ref{subsub:cas}(\ref{indpXM1})$\sim$(\ref{indpXM3}), as well as a product of Siegel modular varieties $\mathcal{A}_{g_\bbI}$ such that $\mathcal{X}_b\subseteq\mathcal{A}_{g_\bbI}\subseteq \mathcal{A}_{g_b}$, and that $f_b$ is an immersion that separates factors over $\mathscr{A}_{g_\bbI, \Fpbar}^{\ord}$. The construction goes as follows:   
\begin{enumerate}[label=(\alph*)]
\item\label{constr:sepa1} Take a point $x_a\in X_a(\Fpbar)$. Take an $\mathcal{X}_{f_a}\in \mathrm{MS}_{f_a}(\mathcal{A}_{g_a})$ and let $\mathcal{X}_a:=\mathcal{X}_{f_a}$. Let $\tilde{x}_{a,\bC}$ be a quasi-canonical lift of $x_a$ that lies in $\mathcal{X}_{f_a}(\bC)$. 
\item\label{constr:sepa2} There is some index set $\bbI$ such that $$\MTT(f_a)^{\ad}=\prod_{i\in \bbI} G_{a,i},$$
where each $G_{a,i}$ is simple and adjoint over $\bQ$. 
\item\label{constr:sepa3}  {Let $M_{a,i}$ be the $\bQ$-subgroup of $\MTT(f_a)$ generated by $\MTT(
{\tilde{x}_{a,\bC}})$ and the preimage of $G_{a,i}$ in $\MTT(f_a)$.}  Then $M_{a,i}^{\ad}=G_{a,i}$, and it defines an almost simple special subvariety $\mathcal{Y}_{a,i}\subseteq\mathcal{X}_a$ that passes through $\tilde{x}_{a,\bC}$. The naïve integral model of $\mathcal{Y}_{a,i}$ has non-empty ordinary locus.
    \item\label{constr:sepa4} Let $\mathcal{A}_{g_\bbI}:=\prod_{i\in \bbI}\mathcal{A}_{g_a}$ (a product of $|\bbI|$ copies of $\mathcal{A}_{g_a}$), viewed as a subvariety of some $\mathcal{A}_{g_b}$.
    \item\label{constr:sepa5} Let $\mathcal{X}_b:=\prod_{i\in \bbI}\mathcal{Y}_{a,i}$. There is a special correspondence $\mathcal{X}_c$ between $\mathcal{X}_a$ and $ \mathcal{X}_b$. There are projections $\pi_a,\pi_b$ from $\mathscr{X}_c$ to $\mathscr{X}_a, \mathscr{X}_b$, inducing quasi-finite surjections on ordinary loci by  Lemma~\ref{lm:extensiontonormal}. 
    \item\label{constr:sepa6} Let $f'_c:X'_c\rightarrow \mathscr{X}_{c,\Fpbar}^{\ord}$ be a quasi-lift of $f_a$ to $\mathscr{X}_{c,\Fpbar}^{\ord}$ along $\pi_a$ (\S\ref{subsub:Quasi-lifts}). 
    \item\label{constr:sepa7} Take closure of $\pi_b\comp f'_c(X'_c)$ in $\mathscr{X}_{b,\Fpbar}^{\ord}$, then take $X_b$ to be a smooth open dense subset of the closure. Take $f_b$ to be the immersion of $X_b$ into $\mathscr{A}_{g_b,\Fpbar}^{\ord}$.
    \item\label{constr:sepa8} Trace back and take $X_c$ to be the preimage of $X_b$ in $X_c'$, and take $f_c$ as the restriction of $f_c'$ to $X_c$. This finishes the construction. 
\end{enumerate} 
Note that $f_b$ is by construction an immersion. Let $f_{b,i}$ be the composition of $f_b$ with the projection of $\mathscr{A}_{g_\bbI}$ to the $i$-th copy of $\mathscr{A}_{g_a}$.  Then $\mathcal{Y}_{a,i}\in \mathrm{MS}_{f_{b,i}}(\mathcal{A}_{g_a})$ and by Lemma~\ref{lm:finitefinite}(\ref{modi:2}) and dimension counting, we have $\prod_{i\in \bbI}\mathcal{Y}_{a,i}=\mathcal{X}_b\in \mathrm{MS}_{f_{b}}(\mathcal{A}_{g_\bbI})$. So $f_b$
separates factors over $\mathscr{A}_{g_\bbI, \Fpbar}^{\ord}$.
\end{construction}

\subsubsection{André--Oort for subvarieties that separate factors}\label{subsub:AOsepfact}
We state the mod $p$ André--Oort for $X$ that separates factors, which has an explicit description of the special factor of $X$ as compared to the qualitative version~\ref{conj:AOAO}. In the following, let $\bbI$ be a finite index set. Suppose that we have an immersion $f:X\hookrightarrow \mathscr{A}_{g, \Fpbar}^\ord$ that separates factors over $\mathscr{A}_{g_\bbI, \Fpbar}^\ord\subseteq \mathscr{A}_{g, \Fpbar}^\ord$. Let $f_i$ and $\mathcal{X}_{f_i}$ be as in the beginning of \S\ref{subsub:spefactor}, and let $\mathscr{X}_{f_i}$ be the naïve integral model of $\mathcal{X}_{f_i}$ in $\mathscr{A}_{g_i}$. 

\begin{conj}[Mod $p$ André--Oort for subvarieties that separate factors]\label{conj:varAO}
Notation as above. Suppose that $X$ contains a Zariski dense collection $\Xi$ of positive dimensional special subvarieties. 
Then there exists an index $i_0\in \bbI$, such that  $\overline{X}$, the Zariski closure of $X$ in $\mathscr{A}_{g_\bbI,\Fpbar}^\ord$, is the product of a component of $\mathscr{X}_{f_{i_0}, \Fpbar}^\ord$ with a subvariety $Y\subseteq\prod_{i\in \bbI\setminus\{i_0\}}\mathscr{X}_{f_{i}, \Fpbar}^\ord$. 
\end{conj}
\begin{remark}\label{rmk:quantatative}
Conjecture~\ref{conj:varAO} admits a quantitative version: Let $\mathbf{I}_{\Xi}\subseteq \bbI$ be the set of indices $i$ such that $\Xi$ contains a Zariski dense sub-collection of special subvarieties whose projections to $\mathscr{A}_{g_i, \Fpbar}^\ord$ are positive dimensional, then $\overline{X}$ is the product of a component of $\prod_{i\in \bbI_{\Xi}}\mathscr{X}_{f_{i}, \Fpbar}^\ord$ with a subvariety of $\prod_{i\in \bbI\setminus\bbI_{\Xi}}\mathscr{X}_{f_{i}, \Fpbar}^\ord$. In particular, if $\bbI=\mathbf{I}_{\Xi}$, then $X$ is special. 

The quantitative version follows directly from Conjecture~\ref{conj:varAO} by projecting the collection $\Xi$ down to $Y\subseteq \prod_{i\in \bbI\setminus\{i_0\}}\mathscr{X}_{f_{i}, \Fpbar}^\ord$ and then iterate.
\end{remark}

\begin{proposition}\label{prop:AOquasilifts}
   Suppose that $f_a:X_a\hookrightarrow \mathscr{A}^{\ord}_{g_a,\Fpbar}$ is an immersion, such that $X_a$ contains a Zariski dense collection $\Xi_a$ of positive dimensional special subvarieties. Let $f_b:X_b\hookrightarrow \mathscr{A}^{\ord}_{g_b,\Fpbar}$ be a modification of $f_a$ which is also an immersion. Then \begin{enumerate}
       \item $X_b$ contains a Zariski dense collection $\Xi_b$ of positive dimensional special subvarieties that specially correspond to elements in $\Xi_a$.
       \item  Conjecture~\ref{conj:AOAO} is true for $X_a$ if and only it is true for $X_b$. 
       \item If $f_b$ separates factors, then Conjecture~\ref{conj:AOAO} holds for $X_a$ if Conjecture~\ref{conj:varAO} holds for $X_b$.
   \end{enumerate}
\end{proposition}
\begin{proof}
    \begin{enumerate}
        \item This is by pull-pushing the elements in $\Xi_a$ along the special correspondence $\mathscr{X}_c$. Let $f_c:X_c\rightarrow \mathscr{X}_{c,\Fpbar}^{\ord}$ and $g_a,g_b$ be as in Definition~\ref{def:modifications}. Consider a sub-collection $\Xi'_a\subseteq \Xi_a$ consisting of elements which lie generically in $g_a(X_c)$. Since $g_a$ is dominant, $\Xi'_a$ is again a Zariski dense collection.       Let $Z_a\in \Xi'_a$. Up to a shrinking, we can assume that it lies in $g_a(X_c)$. Lift it to an irreducible subvariety $Z_c\subseteq X_c$, then take $Z_b$ to be the closure of the image of $Z_c$ in $X_b$. 
        
        We claim that $Z_b$ is a positive dimensional special subvariety of $X_b$. Up to a shrinking, we can assume that $Z_a,Z_b,Z_c$ are all smooth, and 
        $g_a$ maps $Z_c$ to a dense subset of $Z_a$, while $g_b$ maps $Z_c$ to a dense subset of $Z_b$. For $\alpha\in \{a,b,c\}$, let $f'_\alpha$ be the restriction of $f_\alpha$ to $Z_\alpha$. Then $f'_b$ is a modification of $f'_a$ via $f'_c$. Apply Lemma~\ref{lm:finitefinite}(\ref{modi:2}) to $f'_a,f'_b,f'_c$, we see that there is an $\mathcal{X}_{f'_a}\in \mathrm{MS}_{f'_a}(\mathcal{A}_{g_a})$ and $\mathcal{X}_{f'_b}\in \mathrm{MS}_{f'_b}(\mathcal{A}_{g_b})$ such that $ \dim\mathcal{X}_{f'_a}=\dim\mathcal{X}_{f'_b}$. Since $\dim Z_a=\dim Z_b$ by Lemma~\ref{lm:extensiontonormal}, we conclude that $$(\dim\mathcal{X}_{f'_a}=\dim Z_a)\Leftrightarrow(\dim\mathcal{X}_{f'_b}= \dim Z_b).$$ From the observation made in Remark~\ref{rmk:ASimpliesTtl}, this shows that $Z_b$ is special. When $Z_a$ runs over $\Xi'_a$, the corresponding $Z_b$ forms a Zariski dense collection of positive dimensional special subvarieties of $X_b$, as desired. 
\item\label{it:2343}follows easily from the definition of quasi-weakly special subvarieties and the definition of modifications (note that they both use special correspondences). 
\item is similar (\ref{it:2343}): a splitting of $\overline{X}_b$ gives to an almost product splitting of $\overline{X}_a$ via the special correspondence $\mathscr{X}_c$.
    \end{enumerate}
\end{proof}

\subsection{Setups for products of GSpin Shimura varieties}\label{subsub:prodsetup}  
We set up the notation for working with products of GSpin Shimura varieties. The general constructions and results in the previous sections remain valid in this special case. But there are a few more things to add. 

\subsubsection{Two setups} The following are two slightly different setups. Setup~\ref{setup1} is more general, and is essentially the setup where the main theorems are stated in the introduction. Setup~\ref{setup2} is technically simpler. We will reduce the conjectures in the first setup to the second setup using results in \S\ref{sub:inpmM} as well as a variant of Construction~\ref{constr:sepa}.

\begin{setup}\label{setup1}
    Suppose that $\mathcal{A}_{g}$ is an ambient Siegel modular variety, with canonical integral model $\mathscr{A}_g$. Let $\bbI$ be a finite index set. Suppose that $\mathcal{S}_{\bbI}:=\prod_{i\in \bbI}\mathcal{S}_{i}\subseteq \mathcal{A}_{g}$ is a finite product of GSpin Shimura varieties, such that each $\mathcal{S}_i$ is almost simple, and the naïve integral model $\mathscr{S}_{\bbI}\subseteq \mathscr{A}_g$ has nonempty ordinary locus. Let $(X,x)$ be a smooth pointed connected variety over $\Fpbar$ with a morphism \begin{equation}\label{eq:GSpinsetupthing}
    f:(X,x)\rightarrow \mathscr{A}_{g,\Fpbar}^{\ord}
\end{equation}
that factors through ${\mathscr{S}^{\ord}_{\bbI,\Fpbar}}$. Let $f_0:X_0\rightarrow \mathscr{A}_{g,\Fpbar_q}^{\ord}$ be a suitable finite field model of $f$. Equip $X_0$ with the abelian scheme that is pulled back from the universal family over $\mathscr{A}_{g}$, as well as various cohomology sheaves $\bH_{\bullet}$. Fix an $\mathcal{X}_f\in \mathrm{MS}_{f}(\mathcal{S}_\bbI)$. We define $\MTT(f)$ (\textit{resp}. $G_u(f)$) using $\bH_{\bullet}$.
\end{setup}
\begin{setup}\label{setup2}   
Suppose that $\bbI$ is a finite index set, and for each $i\in \bbI$, there is an almost simple GSpin Shimura variety $\mathcal{S}_i$ attached to a self-dual $\bZ_{(p)}$-lattice $(L_i,Q_i)$, with integral canonical model $\mathscr{S}_i$, and a Kuga--Satake embedding of integral canonical models $\mathscr{S}_i\hookrightarrow\mathscr{A}_{g_i}$ associated to the lattice $H_i=\Cl(L_i)$ (\S\ref{subsub:Selfdual}). Consider the product $\mathscr{A}_{g_\bbI}:=\prod_{i\in \bbI} \mathscr{A}_{g_i}$ as a subvariety of an ambient Siegel modular scheme $\mathscr{A}_{g}$. Let $(X,x)$ be a smooth connected pointed variety over $\Fpbar$ with a morphism \begin{equation}\label{eq:GSpinsetupthing2}
    f:(X,x)\rightarrow \mathscr{A}_{g,\Fpbar}^{\ord}
\end{equation}
that factors through ${\mathscr{S}^{\ord}_{\bbI,\Fpbar}}$, and also \textbf{separates factors} over $\mathscr{A}_{g_\bbI,\Fpbar}^{\ord}$. 

Let $f_0:X_0\rightarrow \mathscr{A}_{g,\Fpbar_q}^{\ord}$ be a suitable finite field model of $f$. Instead of equipping $X_0$ with the pullback abelian scheme from $\mathscr{A}_{g}$, we equip each ${\mathscr{S}_{i}}$ with a Kuga--Satake abelian scheme $A_i$ that is pulled back from $\mathscr{A}_{g_i}$, as well as various cohomology sheaves $\bH_{i,\bullet}$, and then equip $X_0$ with the pullback of the product abelian scheme $A_{\bbI}:=\prod_{i\in \bbI}A_i$ and pullback of the product cohomology sheaves $\bH_{\bbI,\bullet}:=\prod_{i\in \bbI}\bH_{i,\bullet}$. Fix an $\mathcal{X}_f\in \mathrm{MS}_{f}(\mathcal{S}_\bbI)$. We define $\MTT(f)$ (\textit{resp}. $G_u(f)$) by $\bH_{\bbI,\bullet}$, it has the same adjoint with the Mumford--Tate group (\textit{resp}. $u$-adic monodromy group) defined using the cohomology sheaves pulled back from $\mathscr{A}_{g}$. 
\end{setup}
In Construction~\ref{constr:varsepa} below we will show that a map in Setup~\ref{setup1} can be modified to a map in Setup~\ref{setup2}. Since the Conjecture~\ref{conj:MT}, \ref{conj:AOAO}, \ref{conj:Ttlinear} are insensitive to modifications by Proposition~\ref{prop:redalongquasilifts},\ref{prop:AOquasilifts}, to prove Conjecture~\ref{conj:MT}, \ref{conj:AOAO}, \ref{conj:Ttlinear} in Setup~\ref{setup1}, it suffices to prove Conjecture~\ref{conj:MT}, \ref{conj:varAO}, \ref{conj:Ttlinear} in Setup~\ref{setup2}.

\textbf{$\bullet$ For this reason, if not otherwise specified, we will always use Setup~\ref{setup2} in the rest of this paper.} 

\begin{construction}\label{constr:varsepa} We run a variant of the construction~\ref{constr:sepa} to modify the map $f$ in Setup~\ref{setup1} to a map in Setup~\ref{setup2}. To make the notation match that of \ref{constr:sepa}, we will attach an subscript ``$a$'' to  (\ref{eq:GSpinsetupthing}), so to write it as $f_a:X_a\rightarrow \mathscr{A}_{g_a,\Fpbar}^{\ord}$.

\begin{enumerate}[label=(\alph*)]
\item\label{constr:varsepa1} Take a point $x_a\in X_a(\Fpbar)$. Take an $\mathcal{X}_{f_a}\in \mathrm{MS}_{f_a}(\mathcal{S}_{\bbI})$ and let $\mathcal{X}_a:=\mathcal{X}_{f_a}$. Let $\tilde{x}_{a,\bC}$ be a quasi-canonical lift of $x_a$ that lies in $\mathcal{X}_{f_a}(\bC)$. 

 \item\label{constr:varsepa2} There is some index set $\bbJ$ such that $$\MTT(f_a)^{\ad}=\prod_{j\in \bbJ} G_{a,j},$$
where each $G_{a,j}$ is simple and adjoint over $\bQ$. 
\item\label{constr:varsepa3}  {Let $M_{a,j}$ be the $\bQ$-subgroup of $\MTT(f_a)$ generated by $\MTT(
{\tilde{x}_{a,\bC}})$ and the preimage of $G_{a,j}$ in $\MTT(f_a)$.}  Then $M_{a,j}^{\ad}=G_{a,j}$, and it defines an almost simple special subvariety $\mathcal{Y}_{a,j}\subseteq\mathcal{X}_a$ that passes through $\tilde{x}_{a,\bC}$. The naïve integral model of $\mathcal{Y}_{a,j}$ has non-empty ordinary locus. There exists a function $s:\bbJ\rightarrow \bbI$, such that the projection $\mathcal{Y}_{a,j}\rightarrow \mathcal{S}_{s(j)}$ is finite. Taking the image of this projection, we get a special subvariety $\mathcal{Y}'_{a,j}\subseteq \mathcal{S}_{s(j)}$. 
 \item\label{constr:varsepa4} For each $j\in \bbJ$, by Remark~\ref{subsub:Nonselfdual}, we can embed the Shimura datum of $\mathcal{S}_{s(j)}$ into 
 a GSpin Shimura datum attached to a self-dual $\bZ_{(p)}$-lattice, then do a Kuga--Satake embedding into a Siegle Shimura datum (\S\ref{subsub:Selfdual}). Let $\mathcal{S}_{j}'$ and $\mathcal{A}_{g_j}$ be a self-dual GSpin Shimura variety and a Kuga--Satake Siegel modular variety that arise from the aforementioned Shimura data, equipped with sufficiently small hyperspecial level structures, so that there is an embedding $\mathscr{S}_{j}'\hookrightarrow\mathscr{A}_{g_j}$ between canonical integral models. Possibly passing $\mathcal{S}_{s(j)}$ to a smaller level structure, and lifting $\mathcal{Y}'_{a,j}$, we also have a chain of embeddings $\mathcal{Y}'_{a,j}\hookrightarrow \mathcal{S}_{s(j)} \hookrightarrow \mathcal{S}_{j}'$.  Let $\mathcal{A}_{g_\bbJ}:=\prod_{i\in \bbJ}\mathcal{A}_{g_j}$, viewed as a subvariety of some large $\mathcal{A}_{g_b}$.  In summary, this step yields embeddings: 
 $$\prod_{j\in \bbJ} \mathcal{Y}'_{a,j} \hookrightarrow \prod_{j\in \bbJ} \mathcal{S}_{s(j)} \hookrightarrow \prod_{j\in \bbJ} \mathcal{S}_{j}'\hookrightarrow \prod_{j\in \bbJ}\mathcal{A}_{g_j}\hookrightarrow \mathcal{A}_{g_b}.$$
    \item\label{constr:varsepa5} Let $\mathcal{X}_b:=\prod_{j\in \bbJ}\mathcal{Y}'_{a,j}$. Then by construction $\mathscr{X}_{b}$ has non-empty ordinary locus. There is a special correspondence $\mathcal{X}_c$ between $\mathcal{X}_a$ and $ \mathcal{X}_b$. There are projections $\pi_a,\pi_b$ from $\mathscr{X}_c$ to $\mathscr{X}_a, \mathscr{X}_b$, inducing quasi-finite surjections on ordinary loci.
    \item\label{constr:varsepa6} Let $f'_c:X'_c\rightarrow \mathscr{X}_{c,\Fpbar}^{\ord}$ be a quasi-lift of $f_a$ to $\mathscr{X}_{c,\Fpbar}^{\ord}$ along $\pi_a$. 
    \item\label{constr:varsepa7} Take closure of $\pi_b\comp f'_c(X'_c)$ in $\mathscr{X}_{b,\Fpbar}^{\ord}$, then take $X_b$ to be a smooth open dense subset of the closure. Take $f_b$ to be the immersion of $X_b$ into $\mathscr{A}_{g_b,\Fpbar}^{\ord}$.
    \item\label{constr:varsepa8} Trace back and take $X_c$ to be the preimage of $X_b$ in $X_c'$, and take $f_c$ as the restriction of $f_c'$ to $X_c$. This finishes the construction. 
\end{enumerate} 
Similar to \ref{constr:sepa}, $f_b$ is by construction an immersion, and 
separates factors over $\mathscr{A}_{g_\bbJ, \Fpbar}^{\ord}$. By \ref{constr:varsepa4}, its projection to $\mathscr{A}_{g_j}$ further factors through the canonical integral model of a GSpin Shimura variety for a self-dual $\bZ_{(p)}$-lattice, whose Kuga--Satake Siegel modular variety is exactly $\mathcal{A}_{g_j}$. Note that $f_b$ is exactly a map in Setup~\ref{setup2}.
\end{construction}

\subsubsection{More notation}\label{subsub:X0fS}
In the Setup~\ref{setup2}, let $(L_{\mathbf{I}},Q_\mathbf{I})=\bigoplus_{i\in\mathbf{I}}(L_i,Q_i)$ and $H_\bbI=\bigoplus_{i\in\mathbf{I}}H_i$. If $R$ is a $\bZ_{(p)}$-algebra, let  $$\SO'(L_{\mathbf{I},R}):=\prod_{i\in\mathbf{I}}\SO(L_{i,R}),\,\,\GSpin'(L_{\bbI,R})=\prod_{i\in\mathbf{I}}\GSpin(L_{i,R}).$$
Let $c:\GSpin'(L_{\bbI,\bZ_{(p)}})\rightarrow \SO'(L_{\mathbf{I},\bZ_{(p)}})$ be the canonical ``conjugation'' map. Note that it induces an isomorphism on adjoints. We identify $H_{\bbI}= \bbH_{\bbI,\mathrm{B},\tilde{x}_\bC}$, where $\tilde{x}_\bC$ is the canonical lift of $x$ base changed to $\bC$. For each $u\in \mathrm{fpl}(\bQ)$, we make an identification $$ H_{\bbI,\bQ_p}\simeq \omega_x(\bH_{\bbI,p})$$ following (\ref{eq:canidentifications}). Let $\mu:\bG_{m}\rightarrow \GL(H_{\bbI,\bQ_p})$ be the canonical Hodge cocharacter of $x$, then it factors through $\GSpin'(L_{\bbI,\bQ_p})$.  The cocharacter  $$\bG_m\xrightarrow{\mu}\GSpin'(L_{\bbI,\bQ_p})\xrightarrow{c}  \SO'(L_{\bbI,\bQ_p})$$ will again be called the canonical Hodge cocharacter of $x$ (but for the $F$-isocrystal $\bL_{p,x}$, see below), and will again be denoted by $\mu$.

\subsubsection{Monodromy groups in the GSpin setting}\label{subsub:monoGsetting} GSpin Shimura varieties are also equipped with sheaves $\bL_{\bullet}$ (\S\ref{subsub:Selfdual}). In Setup~\ref{setup2}, let $\bL_{i,\bullet}$ be the corresponding sheaves associated to $\mathscr{S}_i$, and let $\bL_{\bbI,\bullet}:=\prod_{i\in \bbI}\bL_{i,\bullet}$ be their product over $\mathscr{S}_{\bbI}$. We make the same convention for $\bL_{\bbI,\bullet}$ as in (\ref{eq:Hlp}). Let $${\Hg}(f)\subseteq \GL(L_{\bbI,\bQ})$$ be the generic Hodge group of $\bL_{\bbI,\mathrm{B},\mathcal{X}_{f}}$. It is the image of ${\MTT}(f)\subseteq \GSpin'(L_{\bbI,\bQ})$ in $\SO'(L_{\bbI,\bQ})$ under the map $c$. The kernel of the map $\MTT(f)\twoheadrightarrow {\Hg}(f)$ is the rank one torus $\bbT$ of homotheties.

For $u\in \mathrm{fpl}(\bQ)$, there is an identification $${L}_{\bbI,\bQ_u}\simeq \omega_x(\mathbb{L}_{\bbI,u})$$ that is compatible with the identification $H_{\bbI,\bQ_u}\simeq \omega_x(\mathbb{H}_{\bbI,u})$. Let $$G^L_u(f):=G(\bL_{\bbI, u,X_0},x)^\circ\subseteq \GL(L_{\bbI,\bQ_u}).$$ It is the image of $G_u(f)$ in $\SO'(L_{\bbI,\bQ_u})$. The kernel of the map $G_u(f)\twoheadrightarrow G^L_u(f)$ is $\bbT_{\bQ_u}$. 

All results in \S\ref{subsub:idfi} remain valid for $\bH_{\bbI,u}$. In addition, Proposition~\ref{lm:containedin} and \ref{prop:centercontain} also remain valid, if we replace $\MTT(f),G_u(f)$ by $\Hg(f),G^L_u(f)$. It is also easy to see that Conjecture~\ref{conj:MT} is equivalent to the assertion that $\Hg(f)_{\bQ_u}=G^L_u(f)$.


\subsubsection{Tate linearity in the GSpin setting}\label{subsub:TTGSpinsetting} 
Let $U_{\GSpin',{\mu}^{-1}}$ be the unipotent of $\GSpin'(L_{\bbI,\bQ_p})$ with respect to $\mu$. Let $U_{\SO',{\mu}^{-1}}$ be the unipotent of $\SO'(L_{\bbI,\bQ_p})$ with respect to $\mu$. Then $c$ induces an isomorphism $$c:U_{\GSpin',{\mu}^{-1}}\simeq U_{\SO',{\mu}^{-1}}$$
By Proposition~\ref{prop:cocharofSTspecial} (as well as the paragraph after that), the rational cocharacter lattice of $\mathscr{S}_{\bbI,W}^{/x}$ is $\Lie U_{\GSpin',{\mu}^{-1}}$.

Let $\mathscr{T}_{f,x}$ and $T_{f,x}$ be as in \S\ref{subsub:tmIMPLICA}. We have $T_{f,x}\subseteq \Lie U_{\GSpin',{\mu}^{-1}}$. We can also view $T_{f,x}$ as a subspace of $\Lie U_{\SO',{\mu}^{-1}}$ via $c$. Taking the exponential map, we can view ${T}_{f,x}$ as a subgroup of both $U_{\GSpin',{\mu}^{-1}}$ and $U_{\SO',{\mu}^{-1}}$. Then Theorem~\ref{Thm:Tatelocal} holds if we replace $G(\bH_{p,X_0},x)$ by $G(\bH_{\bbI,p,X_0},x)$ or $G(\bL_{\bbI,p,X_0},x)$. 

\begin{lemma}\label{lm:subHodgeT} 
In Setup~\ref{setup2}, let $\eta$ be the generic point of $X$ and identify $\End(A_{\mathbf{I},\overline{\eta}})_{\bQ_p}$ as a subset of $\End(\omega_x(\bH_{\bbI,p}))$. Let $U\subseteq U_{\GSpin',\mu^{-1}}$ be the subgroup that  
commutes with $\End(A_{\mathbf{I},\overline{\eta}})_{\bQ_p}$. Then $\dim \mathcal{X}_f\leq \dim U$. In particular, if $U\subseteq T_{f,x}$, then Conjecture~\ref{conj:Ttlinear} holds.
\end{lemma}
\proof Let $\tilde{x}_{\bC}\in\mathcal{X}_f(\bC)$ be a quasi-canonical lift of $x$. Let $\mathcal{Y}$ be an irreducible special subvariety of $\mathcal{S}_{\bbI,\overline{\bQ}}$ that is ``cut out by the Hodge cycles $\End^0(A_{\bbI,\overline{\eta}})$ '' and contains $\tilde{x}_{\bC}$ as a $\bC$-point, similar to the proof of Lemma~\ref{lm:endetaepsilon}. More precisely, let $G\subseteq\GSpin'(L_{\bbI})_{\bQ}$ be the subgroup fixing the Hodge tensors in $H_{\bbI,\bQ}^{\otimes (1,1)}$ that span $\End^0(A_{\bbI,\overline{\eta}})\subseteq \End(H_{\bbI,\bQ})$. Then $G$ defines the desired special subvariety $\mathcal{Y}$. Let $\mathscr{Y}$ be the naïve integral model of $\mathcal{Y}$. Using the same argument as Lemma~\ref{lm:endetaepsilon}, we find that (1) $f$ factors through $\mathscr{Y}$ and (2) the opposite unipotent of $G_{\bC}$ with respect to $\mu_{\tilde{x}_\bC}$ has dimension $\dim U$. So $\mathcal{Y}$ has dimension $\dim U$. By minimality of $\mathcal{X}_f$, we have $\dim \mathcal{X}_f\leq \dim U$.

If $U\subseteq T_{f,x}$, then $\dim  \mathcal{X}_f\leq \dim U\leq \dim T_{f,x}$. On the other hand $\mathcal{X}_f\geq \dim T_{f,x}$ by Theorem~\ref{thm:noottheorem}. So we have $\dim  \mathcal{X}_f= \dim T_{f,x}$.
$\hfill\square$

\section{Conjecture~\ref{conj:MT} and \ref{conj:Ttlinear} for products of GSpin Shimura varieties}\label{Sec:ttlinearsingle}
We prove Conjecture~\ref{conj:Ttlinear} \ref{conj:MT} in the context of Setup~\ref{setup2}. As we have seen, this solves the conjectures in the more general Setup~\ref{setup1}. We first treat the single case where $\bbI=1$. Then use the single case to establish the general case.

\subsection{Isocrystals and their monodromy in the single GSpin case}\label{sec:hearts} Let $\#\bbI=1$ in Setup~\ref{setup2}. We will call it the ``single GSpin case''. To ease notation, we drop the subscript $\bbI$. To make it clearer, let's summarize the setup for the single GSpin case in the following:
\begin{setup}\label{setup3}
We have a single GSpin Shimura variety $\mathcal{S}$ attached to a self-dual $\bZ_{(p)}$-lattice $(L,Q)$, with integral canonical model $\mathscr{S}$, as well as a map
\begin{equation}\label{eq:GSpinsetupthing3}
    f:(X,x)\rightarrow \mathscr{S}_{\Fpbar}^{\ord}\hookrightarrow \mathscr{A}_{g,\Fpbar}^{\ord}
\end{equation}
that separate factors (it is equivalent to $\MTT(f)$ being almost simple), where $\mathscr{A}_g$ on the right hand side is from the Kuga--Satake construction.

Let $f_0:X_0\rightarrow \mathscr{S}_{\Fpbar_q}^{\ord}$ be a finite field model. Equipped $X_0$ with the pullback Kuga--Satake abelian scheme. Let $b+2=\rk L$ and $d=\dim {T}_{f,x}$. 
Possibly replacing $X_0$ by a finite étale cover, we may assume that the monodromy groups of $\bH_{u,X_0}$ are connected for all $u\in \mathrm{fpl}(\bQ)$. This implies the monodromy groups of $\bL_{u,X_0}$ are also connected. 

\end{setup}

\subsubsection{Graded objects}\label{subsub:gddobj} For technical reasons, we will use $\bL_{\bullet}$ instead of $\bH_{\bullet}$. See \S\ref{subsub:X0fS}$\sim$\ref{subsub:TTGSpinsetting} for basic setups. Recall that we have the canonical map $c:\GSpin(L_{\bQ_p})\rightarrow\SO(L_{\bQ_p})$ and the canonical Hodge cocharacter $\mu:\bG_m\rightarrow\SO(L_{\bQ_p})$. Then $\mu^{-1}$ gives rise to a ``slope decomposition'': 
$${L}_{{\mathbb{Q}}_p}=L_{-1,{\mathbb{Q}}_p}\oplus L_{0,{\mathbb{Q}}_p} \oplus L_{1,{\mathbb{Q}}_p},\,\,\,L_{j,{\mathbb{Q}}_p}=\{v\in L_{\bQ_p}|\mu^{-1}(t)v=t^{j}v\},$$ which splits $(L_{\bQ_p},Q)$ into nondegenerated subspaces: $$(L_{\bQ_p},Q)=(L_{\bQ_p},Q_0)\oplus (L_{-1,\bQ_p}\oplus L_{1,\bQ_p},Q').$$ 
We have $\dim L_{0,\bQ_p}=b $ and $ \dim L_{-1,\bQ_p}=\dim L_{1,\bQ_p}=1$. The one dimensional subspaces $L_{-1,\bQ_p}$ and $L_{1,\bQ_p}$ are totally isotropic and mutually dual. 
Let $$\overline{L}_{ \bQ_p}:= L_{ 0,\bQ_p} \oplus L_{ 1,\bQ_p}.$$ The cocharacter $\mu$ quotients to a cocharacter $\overline{\mu}:\bG_m\rightarrow \GL(\overline{L}_{{\mathbb{Q}}_p})$. We have a chain of isomorphisms of $\bQ_p$-vector spaces: 
\begin{align}
\label{EQ:gsid1} \Lie U_{\GSpin,\mu^{-1}}\stackrel{c}{\simeq}\Lie U_{\SO,\mu^{-1}}{\simeq} \Lie U_{\GL,\overline{\mu}^{-1}}\simeq L_{1,\bQ_p}^\vee\otimes L_{0,\bQ_p}
\simeq L_{0,\bQ_p},
\end{align}
all but the last isomorphism are canonical, while the last isomorphism is canonical up to a scalar. Therefore ${T}_{f,x}\subseteq \Lie U_{\GSpin,\mu^{-1}}$ (cf. \S\ref{subsub:TTGSpinsetting}) can be canonically identified as a subspace of $L_{0,\bQ_p}$. To ease notation, we will write $T_{f,x}$ as $T_f$.

\subsubsection{Construction of some isocrystals}\label{subsub:cissd} Consider the slope filtration over $\mathscr{S}_{\Fpbar}^{\ord}$:$$0\subseteq\Fil_{-1}\mathbb{L}_{ p}\subseteq \Fil_{0}\mathbb{L}_{ p}\subseteq \Fil_{1}\mathbb{L}_{ p}=\mathbb{L}_{ p}.$$ Let $\gr_j\mathbb{L}_{ p}$ be the graded objects and let $\gr\mathbb{L}_{ p}=\bigoplus \gr_j\mathbb{L}_{ p}$. We have $$G(\gr\mathbb{L}_{p,X_0})\subseteq G(\mathbb{L}_{p,X_0})\subseteq \SO(L_{\bQ_p}).$$
We have the following observations:\begin{enumerate}
    \item\label{eq:333332222} By \S\ref{subsub:TTGSpinsetting}, we can apply Theorem~\ref{Thm:Tatelocal} to $G(\mathbb{L}_{p,X_0})$. We find that $G(\gr\mathbb{L}_{p,X_0})$ is the Levi of $G(\mathbb{L}_{p,X_0})$ fixing $\mu$, and the adjoint action $G(\gr{\bL}_{p,X_0})\acts  \Lie U_{\SO,\mu^{-1}}$ admits $T_{f}$ as a subrepresentation.
    \item  $\gr_j\mathbb{L}_{p,X_0}$
corresponds to a representation $G(\gr\mathbb{L}_{p,X_0})\rightarrow G(\gr_j\mathbb{L}_{p,X_0})\hookrightarrow \GL(L_{j,\bQ_p})$ by Tannakian formalism. In particular, $G(\gr_0\mathbb{L}_{p,X_0})\subseteq \SO(L_{0,\bQ_p})$. 
\item The pairing $\bbQ$ forces $\gr_1\mathbb{L}_{ p,X_0}\simeq\gr_{-1}\mathbb{L}_{p,X_0}^{\vee}$, and $$G(\gr_{-1}{\bL}_{p,X_0}\oplus \gr_{1}{\bL}_{p,X_0})=\im \mu\subseteq \SO(L_{-1,\bQ_p}\oplus L_{1,\bQ_p}).$$
\item For slope reasons, we have $$G(\gr{\bL}_{p,X_0})=G(\gr_0{\bL}_{p,X_0})\times G(\gr_{-1}{\bL}_{p,X_0}\oplus \gr_{1}{\bL}_{p,X_0})=G(\gr_0{\bL}_{p,X_0})\times \im\mu.$$ 
\item  The restriction of the adjoint representation in (\ref{eq:333332222}) to $G(\gr_0{\bL}_{p,X_0})$ is isomorphic, via (\ref{EQ:gsid1}), to the standard representation $G(\gr_0{\bL}_{p,X_0})\acts L_{0,\bQ_p}$. 
Therefore $T_{f}$ is a subrep of $G(\gr_0{\bL}_{p,X_0})\acts L_{0,\bQ_p}$.
\end{enumerate}

Now let $T_f^\perp$ be the orthogonal complement of $T_f$ in $L_{0,\bQ_p}$. The totally isotropic subspace $U_f=T_f\cap T_f^\perp$ is another $G(\gr_0\bL_{p,X})$-subrep. Since $G(\gr_0\bL_{p,X})$ is reductive, there are  $G(\gr_0\bL_{p,X})$-subreps $V_f,V'_f,U'_f$ such that 
\begin{align*}
     &T_{f}=V_{f}\oplus U_{f},\\
     &T_{f}^\perp= U_{f}\oplus V_{f}',\\
     &(V_{f}\oplus V_{f}')\perp (U_{f}\oplus U_{f}').
\end{align*}Note that $V_{f},V_{f}'$ and $U_{f}\oplus U_{f}'$ are all nondegenerate, and $U_{f}$ is a maximal totally isotropic subspace of $U_{f}\oplus U_{f}'$. Since $G(\gr_0\bL_{p,X_0})\rightarrow \SO(U_{f}\oplus U_{f}')$ factors through a unitary subgroup $\SU(\varphi_{f,0})$ with $\Tr\varphi_{f,0}=Q|_{U_{f}\oplus U_{f}'}$, we can take $U_{f}'$ to be dual to $U_{f}$. Define \begin{align*}
    M_f&:=(L_{-1,\bQ_p}\oplus L_{1,\bQ_p})\oplus (U_f\oplus U'_f)\oplus V_f,\\
    W_{f}&:= U_{f}\oplus L_{1,\bQ_p},\\
    W'_f&:=U_f'\oplus L_{-1,\bQ_p}.
\end{align*}
The space $W_f\oplus W'_f$ is equipped with an Hermitian form $\varphi_f$ extending $\varphi_{f,0}$, and such that $\Tr\varphi_{f}=Q|_{W_f\oplus W'_f}$. It is easy to check that 
$M_f,V_f'$ are $G(\bL_{p,X_0}^-)$-reps (as usual, $\bL_{p,X_0}^-$ is the underlying $F$-isocrystal of $ \bL_{p,X_0}$), whereas $W_{f},W_{f}'$ are $G(\bL_{p,X_0})$-reps as soon as $V_f=0$.

Note that $L_{\bQ_p}=M_f\oplus V_f'$, so $M_f$ and $V_f'$ can be obtained as the kernel and image of a certain idempotent $\alpha\in \End_{G(\bL_{p,X_0}^-)}(L_{\bQ_p})$. Kedlaya's result \cite{Ked01} implies that $\alpha\in\End_{G(\bL_{p,X_0})}(L_{\bQ_p})$. Therefore $M_f$ and $V_f'$ gives rise to $\bM_{p,f},\bV'_{p,f}\in \FIsoc^{\dagger}(X_0)$, and \begin{equation}\label{eq:L=V+M}
    \bL_{p,X_0}=\bV'_{p,f}\oplus\bM_{p,f}.
\end{equation} Similarly, when $V_f=0$,  $W_f$ and $W_f'$ gives rise to $\bW_{p,f},\bW'_{p,f}\in \FIsoc^{\dagger}(X_0)$, which are dual to each other, and $\bM_{p,f}=\bW_{p,f}\oplus\bW_{p,f}'$.

\begin{proposition}
\label{cor:ShapeT_f}
Either $T_{f}=V_f$ or $T_{f}=U_f$. In other words, $T_{f}$ is either nondegenrate or totally isotropic. In the first case $M_f=(L_{-1,\bQ_p}\oplus L_{1,\bQ_p})\oplus T_f$. In the second case $M_f=W_f\oplus W_f'$, and is equipped with an Hermitian form $\varphi_f$. We have $$G(\bM_{p,f})=\left\{\begin{aligned}
    &\SO(M_{f}),\,\,&T_{f}=V_f,\\
    &\SU(\varphi_f),\,\,&T_{f}=U_f.
\end{aligned}\right. $$  
\end{proposition}
\proof Let $M_{f,0}=U_f\oplus U'_f\oplus V_f$ and let $\gr_0\bM_{p,f}$ be the slope 0 graded part of $\bM_{p,f}$. Let $\nu=\mu^{-1}$. To ease notation, we denote $G(\bM_{p,X_0})$ by $G$, and denote $G(\gr_0\bM_{p,X_0})$ by $G_0$. Whenever $H\subseteq \SO(M_{f})$ is a reductive subgroup containing $\im \nu$, we use $U_{H,\nu}$ to denote the corresponding unipotent. In order to make it easier to track the spaces, we write $\SU(U_f\oplus U_f')$ \textit{resp}. $\SU(W_f\oplus W_f')$ for $\SU(\varphi_{f,0})$ \textit{resp}. $\SU(\varphi_{f})$.

Theorem~\ref{MDA} implies that $G$ is a reductive subgroup of $\SO(M_{f})$ which admits $G(\bM_{p,f}^-)$ as the parabolic subgroup corresponding to $\nu$. The Levi of $G$ corresponding to $\nu$ is $G_0\times \im \nu$. The projection of $G_0$ to $\GL(V_{f})$ \textit{resp}. $\GL(U_{f}\oplus U_{f}')$ lies in $\SO(V_{f})$ \textit{resp}. $\SU(U_{f}\oplus U_{f}')$, so we have
\begin{equation}\label{eq:uv1}
G_0\subseteq \SO(V_{f})\times \SU(U_{f}\oplus U_{f}').
\end{equation}
On the other hand, $U_{G,\nu}=T_{f}$ by Theorem~\ref{Thm:Tatelocal}. Let $N_{f}=V_{f}\oplus L_{1,\bQ_p}\oplus L_{-1,\bQ_p}$. We have  
\begin{equation}\label{eq:uv2}
U_{G,\nu}=U_{\SO(N_f),\nu}\times U_{\SU(W_{f}\oplus W_f'),\nu}.    
\end{equation}

The rest of the proof is purely group theoretical. In the following we base change to $\overline{\bQ}_p$, while keeping the notation unchanged. 

Fix maximal tori $T_{\SO(V_f)}$, $T_{\SU(U_f\oplus U_f')}$ and $T_{G_0}$, such that $T_{G_0}\subseteq T_{\SO(V_f)}\times T_{\SU(U_f\oplus U_f')}= T_{\SO(M_{f,0})}$. Let  $T_{\SO(M_f)}$, $T_{\SO(N_f)}$, $T_{\SU(W_f\oplus W_f')}$ and $T_{G}$ be the products of $\im \nu$ with $T_{\SO(M_{f,0})}$, $T_{\SO(V_f)}$, $T_{\SU(U_f\oplus U_f')}$ and $T_{G_0}$, respectively. 

Note that $T_{\SO(V_f)}$ commutes with $\SU(W_f\oplus W_f')$, whereas $T_{\SU(U_f\oplus U_f')}$ commutes with $\SO(N_f)$. We see that the root sub-algebra of $\mathfrak{so}(N_f)$ \textit{resp}. $\mathfrak{u}(W_f\oplus W_f')$ associated to a positive root is also a root sub-algebra of $\mathfrak{g}:=\Lie G$ associated to a positive root (here positive roots mean the ones associated to the unipotent corresponding to $\nu$). On the other hand, the root sub-algebra of $\mathfrak{so}(N_f)$ \textit{resp}. $\mathfrak{u}(W_f\oplus W_f')$ associated to any root is a root sub-algebra of $\mathfrak{so}(M_f)$.

One deduces that the root sub-algebra of $\mathfrak{g}$ associated to a negative root is also a root sub-algebra of $\mathfrak{so}(N_f)$ or $\mathfrak{u}(W_f\oplus W_f')$ associated to a negative root, hence  
\begin{equation}\label{eq:uv3}
U_{G,\nu^{-1}}=U_{\SO(N_f),\nu^{-1}}\times U_{\SU(W_{f}\oplus W_f'),\nu^{-1}}.    
\end{equation}
Now $U_{\SO(N_f),\nu}$ and $U_{\SO(N_f),\nu^{-1}}$ generates $\SO(N_f)$, while ${U_{\SU(W_{f}\oplus W_f'),\nu}}$, ${U_{\SU(W_{f}\oplus W_f'),\nu^{-1}}}$ and $\im \nu$ generates ${\SU(W_{f}\oplus W_f')}$. We see that $G$ contains both $\SO(N_f)$ and $\SU(W_{f}\oplus W_f')$. In particular, $T_G=T_{\SO(M_f)}$. Therefore the root system of $G$ is a sub root system of $\SO(M_f)$. Easy calculation then shows that either $U_f=0$ and $G=\SO(M_f)$, or $V_f=0$ and $G=\SU(W_f\oplus W_f')$.
$\hfill\square$
\begin{remark}\label{eq:mtbaby}
From the proof of Proposition~\ref{cor:ShapeT_f} we see that 
$$G(\gr_0\bM_{p,f})=\left\{\begin{aligned}
    &\SO(T_{f}),\,\,&T_{f}\text{ is nondegenerate},\\
    &\SU(\varphi_{f,0}),\,\,&T_{f}\text{ is totally isotropic}.
\end{aligned}\right. $$    
\end{remark}

\begin{corollary}
 \label{cor:split1}
$G(\bL_{p,X_0})=G(\bV_{p,f}')\times G(\bM_{p,f}).$
\end{corollary}
\proof By (\ref{eq:L=V+M}), we have $G(\bL_{p,X_0})\subseteq G(\bV_{p,f}')\times G(\bM_{p,f})$, and the projection of $G(\bL_{p,X_0})$ to each factor is surjective. Since any object in the Tannakian subcategory $\la\bV_{p,f}'\ra^\otimes$ has zero slope, it suffices to show the following:
\begin{claim}
 An object in $\la\bM_{p,f}\ra^\otimes$ has zero slope only when it is a direct sum of trivial objects.
\end{claim}
Indeed, by Goursat's lemma, there is a common quotient $H$ of $G(\bV_{p,f}')$ and $G(\bM_{p,f})$, such that
$G(\bL_{p,X_0})=G(\bV_{p,f}')\times_H G(\bM_{p,f})$. A faithful representation of $H$ gives rise to an object that lies in both $\la\bV_{p,f}'\ra^\otimes$ and $\la\bM_{p,f}\ra^\otimes$. If the \textit{Claim} is true, then any faithful representation of $H$ is trivial, so $H=\{1\}$. Therefore $G(\bL_{p,X_0})=G(\bV_{p,f}')\times G(\bM_{p,f})$.

We now prove the \textit{Claim}. Let $G=G(\bM_{p,f})$
and $\nu=\mu^{-1}$. An object $\mathbb{X}\subseteq\la\bM_{p,f}\ra^\otimes$ with zero slope gives rise to a representation $G\rightarrow \GL(\omega_x(\mathbb{X}))$, whose kernel $N$ is a normal subgroup containing $\im\nu$. By Proposition~\ref{cor:ShapeT_f}, $G$ is either $\SO(M_{f})$ or $\SU(\varphi_{f})$. It is easy from the explicit group theory that $N=G$. The proof of the \textit{Claim} is complete. $\hfill\square$

\subsection{Proof of the conjectures in the single GSpin  case}\label{sub:SingleGSpinTl} Let the setup and notation be as in \S\ref{sec:hearts}. Let $\eta$ be the generic point of $X_0$. 
\subsubsection{Proof of Conjecture~\ref{conj:Ttlinear}} Let $\dim V_{f}'=r$. We have constructed an overconvergent $F$-isocrystal $\bV'_{p,f} \subseteq \bL_{p,X_0}$. From Remark~\ref{rmk:embwedge}, we know that there is a chain of embeddings \begin{equation}
\det{\bV}'_{p,f}\subseteq \wedge^{r}\bL_{p,X_0}\subseteq\mathcal{E}nd(\bH_{p,X_0}).
\end{equation}
Since $G(\mathbb{V}'_{p,f})\subseteq \SO(V_{f}')$, the uniroot $F$-isocrystal  $\det\mathbb{V}'_{p,f}$ is constant. Pick a $\bQ_p$-generator ${\bm{\delta}}_{\mathbf{v}}$ of $\det V_{f}'$. By 
the equivalence between the category of Dieudonné crystals and the category of $p$-divisible groups over $X_0$ (\cite{DJ95}) and
crystalline Tate's conjecture (\cite{DJ98}), we have  $${\bm{\delta}}_{\mathbf{v}}\in \End(\bH_{p,\eta})=\End(A_{\eta})\otimes \bQ_p.$$ 

Now suppose $V_{f}=0$, i.e., $T_{f}$ is totally isotropic and $L_{\bQ_p}=W_{f}\oplus W'_{f}\oplus V'_{f}$. Let $\pm1\neq \gamma\in \bZ_p^*$. We define an isometry ${\bm{\delta}}_{\mathbf{w}}\in \End(L_{\bQ_p})$ by \begin{equation}
    {\bm{\delta}}_{\mathbf{w}}(v)=\left\{\begin{aligned}
       & \gamma v, \,\,\, v\in W_{f}, \\
       &\gamma^{-1}v,\,\,\,v\in W_{f}',\\
       & v,\,\,\,v\in V'_{f}.
    \end{aligned}\right.
\end{equation}
Since $G(\bL_{p,X_0})\subseteq \SU(\varphi_f)\times \SO(V_{f}')$, ${\bm{\delta}}_{\mathbf{w}}$ is also an isometric isomorphism of $G(\bL_{p,X_0})$-representations. By funtoriality of the Kuga--Satake construction, the equivalence between the category of Dieudonné crystals and the category of $p$-divisible groups, and the 
crystalline isogeny theorem over finite generated fields, we again have  ${\bm{\delta}}_{\mathbf{w}}\in\End(A_{\eta})\otimes \bQ_p$.

\begin{theorem}\label{thm:TlinSingle}
  Conjecture~\ref{conj:Ttlinear} holds in Setup~\ref{setup3}.
 \end{theorem} \proof 
 When $T_{f}$ is nondegenerate, let $\bm{\Delta}=\{\bm{\delta}_{\mathbf{v}}\}\subseteq \End(A_{\overline{\eta}})_{\bQ_p}$. The subgroup of $U_{\GSpin,\mu^{-1}}$ that commutes with $\bm{\Delta}$ is carried by $c$ to the subgroup of $U_{\SO,\mu^{-1}}$ that preserves $V'_{f}$, which is exactly $T_{f}$. So the subgroup of $U_{\GSpin,\mu^{-1}}$ that commutes with $\End(A_{\overline{\eta}})_{\bQ_p}$ is contained in $T_f$. Therefore Conjecture~\ref{conj:Ttlinear} holds by Lemma~\ref{lm:subHodgeT}. When $T_{f}$ is totally isotropic, let $\bm{\Delta}=\{\bm{\delta}_{\mathbf{v}},\bm{\delta}_{\mathbf{w}}\}$. The subgroup of $U_{\GSpin,\mu^{-1}}$ that commutes with $\bm{\Delta}$ is carried by $c$ to the subgroup of $U_{\SO,\mu^{-1}}$ that preserves both $V'_{f}$ and $W_{f}$. This group is again $T_{f}$. By Lemma~\ref{lm:subHodgeT} again, Conjecture~\ref{conj:Ttlinear} holds.
$\hfill\square$
\subsubsection{Structure of $\Hg(f)$}\label{subsec:MTgpclass}Let $\varrho:{\Hg}(f)\hookrightarrow \SO(L_{\bQ})$ be the standard representation (\S\ref{subsub:monoGsetting}). We write $(L_{\bQ}^\varrho,Q^{\varrho})$ for the subspace of $(L_\bQ,Q)$ over which ${\Hg}(f)$ acts trivially, which is necessarily negative definite. The following proposition describes the structure of $\Hg(f)$ in terms of $T_f$ and $M_f$, and shows that $G(\bL_{p,X_0})$ and $\Hg(f)_{\bQ_p}$ have a factor in common: it is the factor containing the image of $\mu$.

\begin{proposition}[Structure of $\Hg(f)$]\label{Thm:classMTgp}
Let the setup be as in Setup~\ref{setup3}, and let $ T_f, M_f$ be as in \S\ref{sec:hearts}. Let $d=\dim T_f$. The possible structures of ${\Hg}(f)$ lie in the following two cases:
\begin{enumerate}[label=(\alph*)]
    \item\label{it:orthog} $T_{f}$ is nondegenerate, and there is a quadruple $(F,\tau,\mathcal{M}, \mathcal{Q})$, where $F$ 
is totally real, $\tau:F\rightarrow \bC$ is an embedding, and $(\mathcal{M}, \mathcal{Q})$ is a quadratic $F$-subspace of $(L,Q)\otimes F$, such that \begin{enumerate}[label=(\subscript{a}{{\arabic*}})]
    \item
\label{it:nonde1}$\Hg(f)=\Res_{F/\bQ}\SO(\mathcal{Q})$, 
 \item\label{it:nonde2} 
$Q=\Tr_{F/\mathbb{Q}}(\mathcal{Q})\oplus Q^\varrho$,
    \item\label{it:nonde3} $\mathcal{Q}$ has signature $(2,d)$ at place $\tau$ and is negative definite at all other real places,
    \item\label{it:nonde4} Let $\mathfrak{p}$ be the place of $F$ given by $F\xrightarrow{\tau}\bC\simeq\bC_p$. Then $F_\mathfrak{p}=\bQ_p$, $
    (M_{f}, Q|_{M_{f}})=(\mathcal{M}, \mathcal{Q})_{F_{\mathfrak{p}}}$, and $G(\bM_{p,f})=\SO(\mathcal{Q})_{F_\mathfrak{p}}$.  
\end{enumerate}
\item\label{it:unitar} $T_{f }$ is totally isotropic, and there is a hextuple $(F,\tau,\mathcal{M}, \mathcal{Q}; E, {{\Phi}})$, where $(F,\tau,\mathcal{M}, \mathcal{Q})$ is a quadruple as in \ref{it:orthog}, $E/F$ is a quadratic imaginary extension, $(\mathcal{M},{\Phi})$ is a Hermitian space over $E$ with $\mathcal{Q}=\Tr_{E/F}({{\Phi}})$, such that
\begin{enumerate}[label=(\subscript{b}{{\arabic*}})]
    \item\label{it:ttiso1}
    $\Hg(f)=\Res_{F/\bQ}\SU({{\Phi}})$, 
   
    \item\label{it:ttiso2}  $Q=\Tr_{E/\mathbb{Q}}({\Phi})\oplus Q^\varrho$,
     \item\label{it:ttiso3} $\Phi$ has signature $(1,d)$ at place $\tau$ and is negative definite at all other real places,
    \item\label{it:ttiso4} Let $\mathfrak{p}$ be the place of $F$ given by $F\xrightarrow{\tau}\bC\simeq {\bC}_p$. Then $F_{\mathfrak{p}}=\bQ_p$, $\mathfrak{p}$ splits completely in $E$, $(M_{f }, \varphi_{f})=(\mathcal{M}, \Phi)_{F_{\mathfrak{p}}}$, and $G(\bM_{p,f})=\SU({{\Phi}})_{F_\mathfrak{p}}$.
\end{enumerate}
\end{enumerate}
\end{proposition}
\proof  It is well-known that the Hodge group of a special subvariety of a GSpin Shimura variety is the scalar restriction of an orthogonal group or a unitary group over a totally real field with signature conditions as described in the proposition (see for example \cite{Z83}).

Suppose that $\Hg(f)$ is the scalar restriction of an orthogonal group. Then there is a quadruple $(F,\tau,\mathcal{M}, \mathcal{Q})$ where $F$ 
is totally real field, $\tau:F\rightarrow \bC$ is an embedding, and $(\mathcal{M}, \mathcal{Q})$ is a quadratic space over $F$ satisfying \ref{it:nonde1} and \ref{it:nonde2}, and that $\mathcal{Q}$ has signature $(2,l)$ for some integer $l$ at place $\tau$ and is negative definite at all other real places. Theorem~\ref{thm:TlinSingle} implies that $l=\dim\mathcal{X}_{f}=\dim {T}_{f}=d$, so \ref{it:nonde3} holds. Recall that we have identified both $\Hg(f)_{\bQ_p}$ and $G(\bL_{p,X_0})$ as subgroups of $\SO(L_{\bQ_p})$. Let $\mu:\bG_m\rightarrow\SO(L_{\bQ_p})$ be the canonical Hodge cocharacter of $x$, which is also identified with the Hodge cocharacter of the canonical lift $\tilde{x}_{\bC}$ (which lies \textit{a priori} in $\SO(L_{\bC_p})$ but descends to $\bQ_p$). By Corollary~\ref{cor:split1}, $G(\bM_{p,f})$ is the unique factor of $G(\bL_{p,X_0})$ containing $\im \mu$. On the other hand, $\SO(\mathcal{Q})_{\tau,\bC_p}$ is the unique factor of $\Hg(f)_{\bC_p}$ containing $\im\mu$. Since $\mu$ is defined over $\bQ_p$, this implies that $F_\mathfrak{p}=\bQ_p$. The unipotent of $\SO(\mathcal{Q})_{F_\mathfrak{p}}$ with respect to $\mu^{-1}$ has dimension $d$, so must coincide with $T_{f}$. Therefore $T_f$ is nondegenerate. Combining Proposition~\ref{cor:ShapeT_f}, we obtain \ref{it:nonde4}. 

The unitary case is similar, and is left as an exercise.$\hfill\square$

\subsubsection{Proof of Conjecture~\ref{conj:MT}}\label{subsub:cptsystem1} Let $F,\tau$ be the totally real field and its complex embedding from Proposition~\ref{Thm:classMTgp}. There is an $F$-group $\mathcal{G}$, which is either $\SO(\cQ)$ or $\SU(\Phi)$, 
such that
$\Hg(f)=\Res_{F/\bQ}\mathcal{G}$. The group $\mathcal{G}$ is equipped with a standard representation $\phi:\mathcal{G}\rightarrow \GL(\mathcal{M})$. Let $C\subseteq \bC$ be a sufficiently large Galois extension of $\bQ$ containing all complex embeddings of $F$. Let $\Sigma$ be the set of embeddings of $F$ into $C$. By Tannakian formalism, for each $u\in \text{fpl}(\bQ)$, the 
representation $$G(\bL_{u,X_0})\hookrightarrow \Hg(f)_{\bQ_u}\xrightarrow{(\Res_{F/\bQ}\phi)_u} \GL((\Res_{F/\bQ}\mathcal{M})_{\bQ_u})$$
corresponds to a $\bQ_u$-coefficient object $\mathscr{M}_u$. Since $L_{\bQ}=\Res_{F/\bQ}\mathcal{M}\oplus L_\bQ^{\varrho}$, we see that $\bL_{u,X_0}$ is a direct sum of $\mathscr{M}_u$ with a trivial local system $\bL_{u,X_0}^{\varrho}$. Similarly, for each $u\in \text{fpl}(\bQ)$, $v\in\text{fpl}(C)$ dividing $u$, and $\sigma\in \Sigma$, the representation 
$$G(\bL_{u,X_0})_{C_v}\hookrightarrow \Hg(f)_{C_v}\rightarrow \mathcal{G}_{\sigma,C_v}\xrightarrow{\phi_{\sigma,C_v}}\GL(\mathcal{M}_{\sigma,C_v})$$
gives rise to a $C_v$-coefficient object $\mathscr{M}_{\sigma,v}\subseteq \mathscr{M}_u\otimes C_v$. We have $\mathscr{M}_u\otimes C_v=\bigoplus_{\sigma\in \Sigma}\mathscr{M}_{\sigma,v}$. For any pair $(\sigma,v)\in \Sigma\times \mathrm{fpl}(C)$, we have the following relations between monodromy groups: 
\begin{equation}\label{eq:goodequalityilikeit}
    G(\mathscr{M}_{\sigma,v})\subseteq \mathcal{G}_{\sigma,C_v},\,\,\,\,G(\mathscr{M}_{u}\otimes C_v)\subseteq G(\mathscr{M}_{\sigma,v})\times G(\bigoplus_{\sigma'\in \Sigma-\{\sigma\}}\mathscr{M}_{\sigma',v}).
\end{equation}

\begin{lemma}\label{lm:cptsystem}
The collection $\{\mathscr{M}_{u}\}_{u\in \mathrm{fpl}(\bQ)}$ is a weakly $\bQ$-compatible system of coefficient objects in the sense of \S\ref{subsub:compatible}. For each $\sigma\in \Sigma_F$, the collection $\{\mathscr{M}_{\sigma,v}\}_{v\in \mathrm{fpl}(C)}$ is a weakly $C$-compatible system of coefficient objects. 
\end{lemma}
\proof In one word, this follows from the existence of the canonical lift. We only prove the second assertion, and the first follows from that. Let $x_0\in |X_0|$ be a closed point. Since $x_0$ is ordinary, we can canonical lift the geometric Frobenius of $A_{x_0}$ to characteristic 0. This gives rise to a conjugate class $\epsilon$ in $\GL(L_{\bQ})$. We denote by $P(\mathcal{M}_{\sigma,C},t)$ the characteristic polynomial for $\epsilon$ by passing to the quotient $\mathcal{M}_{\sigma,C}$. We denote by $P(\mathscr{M}_{\sigma,v,x_0},t)$ the characteristic polynomial for $\mathscr{M}_{\sigma,v}$ at $x_0$.  Suppose that $v\nmid p$, then we have
\begin{equation}\label{eq:erwerwe}
P(\cM_{\sigma,C},t)=P(\mathscr{M}_{\sigma,v,x_0},t).
\end{equation}
For $v|p$ if we can replace the Frobenius of $\mathscr{M}_{\sigma,v,x_0}$ by a suitable power, (\ref{eq:erwerwe}) still holds. Since the left hand side of (\ref{eq:erwerwe}) depends only on $\sigma$, this proves the weakly-compatibility. 
$\hfill\square$
\begin{lemma}\label{lm:atoneplace}
Let $\mathfrak{p}$ be the place of $F$ induced by the embedding $\tau$ (cf. Proposition~\ref{Thm:classMTgp}). Take a place $\mathfrak{P}|p$ of $C$ lying over $\tau \mathfrak{p}$. Then (\ref{eq:goodequalityilikeit}) are equalities for the pair $(\tau,\mathfrak{P})$.
\end{lemma}
\proof Note that $\mathscr{M}_{\tau,\mathfrak{P}}=\bM_{p,f}\otimes C_{\mathfrak{P}}$, and 
$$(\bL_{p,X_0}^{\varrho}\otimes C_{\mathfrak{P}})\oplus\bigoplus_{\sigma'\in \Sigma-\{\tau\}}\mathscr{M}_{\sigma',\mathfrak{P}}=\bV'_{p,f}\otimes C_{\mathfrak{P}}.$$ 
So the lemma follows from  Corollary~\ref{cor:split1} and Proposition~\ref{Thm:classMTgp}.$\hfill\square$


\begin{construction}[Galois twists]
  Note that $\Gal({C}/\mathbb{Q})$ acts on $\Sigma$. Let $g\in \Gal({C}/\mathbb{Q})$, the \textbf{$g$-twist} of a $C_v$-scheme $X$ is the $C_{gv}$-scheme $X\times_g \Spec C_{gv}$, denoted by $gX$. There is also an obvious notion of $g$-twist of morphisms between $C_v$-schemes. There is also a notion of \textbf{$g$-twist of $C_v$-coefficient objects} defined as $g\mathcal{E}:= \mathcal{E}\otimes_{g}C_{gv}$. We have 
\begin{equation}\label{lm:Galtwist}
g\mathcal{G}_{\sigma,{C}_v}=  \mathcal{G}_{g\sigma,{C}_{gv}},\,\,\,g\mathscr{M}_{\sigma,v}=\mathscr{M}_{g\sigma,gv},\,\,\,,gG(\mathscr{M}_{\sigma,v})=G(g\mathscr{M}_{\sigma,v}).
\end{equation}  
  
\end{construction} 

\begin{theorem}\label{thm:MTinSingle}
    Conjecture~\ref{conj:MT}  holds in Setup~\ref{setup3}.
 \end{theorem}
\begin{proof}
As observed in \S\ref{subsub:monoGsetting}, it suffices to show that $\Hg(f)_{\bQ_u}=G(\bL_{u,X_0})$ for $u\in \mathrm{fpl}(\bQ)$. In the following, we will use (\ref{eq:goodequalityilikeit})$_{(\sigma,v)}$ to denote 
(\ref{eq:goodequalityilikeit}) for the pair $(\sigma,v)$.

Let $\mathfrak{p}$ and $\mathfrak{P}$ be as in Lemma~\ref{lm:atoneplace}, then 
(\ref{eq:goodequalityilikeit})$_{(\tau,\mathfrak{P})}$ are equalities. Now we claim that for all $v\in \mathrm{fpl}(C)$, (\ref{eq:goodequalityilikeit})$_{(\tau,v)}$ are equalities: using weakly-compatibility for $\{\mathscr{M}_{\tau,v}\}_{v\in \mathrm{fpl}(C)}$ (Lemma~\ref{lm:cptsystem}) and  Lemma~\ref{lm:compatible}, there is a finite extension $\mathfrak{C}/C$ with the property that for each $v\in \mathrm{fpl}(C)$ and  $\mathfrak{v}\in \mathrm{fpl}(\mathfrak{C})$ with $\mathfrak{v}|v$, the base changes of (\ref{eq:goodequalityilikeit})$_{(\tau,v)}$ to $\mathfrak{C}_{\mathfrak{v}}$ are equalities. By dimension counting, (\ref{eq:goodequalityilikeit})$_{(\tau,v)}$ are already equalities. 

Now let $\sigma\in \Sigma$, Chebotarev's density theorem guarantees the existence of a place $v_0$ and an element $g\in \Gal(C/\bQ)$ fixing $v_0$ such that $g\tau=\sigma$. Now take $g$-twist of (\ref{eq:goodequalityilikeit})$_{(\tau,v_0)}$, we see that (\ref{eq:goodequalityilikeit})$_{(\sigma,v_0)}$ are equalities. Using weakly-compatibility for $\{\mathscr{M}_{\sigma,v}\}_{v\in \mathrm{fpl}(C)}$ and  Lemma~\ref{lm:compatible} again, we find that (\ref{eq:goodequalityilikeit})$_{(\sigma,v)}$ are equalities for any $v\in \mathrm{fpl}(C)$. Let $\sigma$ run over $\Sigma$, we see that (\ref{eq:goodequalityilikeit})$_{(\sigma,v)}$ are equalities for all $(\sigma,v)\in \Sigma\times \mathrm{fpl}(C)$.  It then follows that $$G(\bL_{u,X_0}\otimes C_v)=G(\mathscr{M}_{u}\otimes C_v)=\prod_{\sigma\in \Sigma}\mathcal{G}_{\sigma,{C}_v}=(\Res_{F/\bQ}\mathcal{G})_{C_v}.$$
As a result, $$G(\bL_{u,X_0} )=(\Res_{F/\bQ}\mathcal{G})_{\bQ_u} =\Hg(f)_{\bQ_u}.$$  
This finishes the proof. 
\end{proof}
\subsection{Proof of the conjectures in the product GSpin case} Now let's go back to Setup~\ref{setup2}. Possibly replacing $X_0$ by a finite étale cover, we may assume that the monodromy groups of $\bH_{\bbI,u,X_0}$ are connected for all $u\in \mathrm{fpl}(\bQ)$. Let $f_i:X\rightarrow \mathscr{A}_{g_i,\Fpbar}^{\ord}$ be the composition of $f$ with $\mathscr{A}_{g_{\bbI},\Fpbar}^{\ord}\rightarrow \mathscr{A}_{g_i,\Fpbar}^{\ord}$.  Each $f_i$ fits in the single GSpin Setup~\ref{setup3}. Since $f$ separates factors, we have $$\MTT(f)^{\ad}=\prod_{i\in \bbI}\MTT(f_i)^{\ad},$$ and each $\MTT(f_i)^{\ad}$ is simple. 

\subsubsection{Compatible adjoint coefficient objects} For $i\in \bbI$, there exists an  absolutely simple reductive group over a number field $N_i$, such that $$\MTT(f_i)^{\ad}=\Res_{N_i/\bQ}\mathcal{H}_i.$$ Since $\MTT(f)^{\ad}=\Hg(f)^{\ad}$, it is easy to see from Proposition~\ref{Thm:classMTgp} what $N_i$ and $\mathcal{H}_i$ should be: in most cases $N_i$ equals the totally real field $F_i$ as in \textit{loc.cit} and $\mathcal{H}_i$ is the adjoint of either an $\SO(\mathcal{Q}_i)^{\ad}$ or a $\SU(\Phi_i)^{\ad}$, except when $\mathcal{Q}_i$ is of rank 4: in that case  $\SO(\mathcal{Q}_i)$ won't be absolutely simple, and instead, $N_i$ is a quadratic extension of $F_i$, and $\mathcal{H}_{i}$ is a descent of $\SO(2,1)$.

Let $\Sigma'_{i}$ be the set of complex embeddings of $N_i$. Fix $C\subseteq\bC$ to be a sufficiently large number field which contains the image of all complex embedding of all $N_i$. For each $i\in \bbI$, let $\phi_i:\mathcal{H}_i\hookrightarrow \GL(\mathcal{V}_i)$ be a faithful irreducible representation over $N_i$. By Tannakian formalism, for each $i\in \bbI$ and $u\in \mathrm{fpl}(\bQ)$, the representation $$G_u(f)\hookrightarrow \MTT(f)_{\bQ_u}\xrightarrow{(\Res_{N_i/\bQ}\phi_{i})_{\bQ_u}} \GL((\Res_{N_i/\bQ} \mathcal{V}_i)_{\bQ_u})$$
gives rise to a $\bQ_u$-coefficient object $\mathscr{E}_{i,u}\in \langle \mathbb{H}_{i,u,X_0}\rangle^{\otimes}$. For each $i\in \bbI$, $\sigma\in \Sigma'_i$, $u\in \mathrm{fpl}(\bQ)$ and $v\in\text{fpl}(C)$ dividing $u$, the representation
$$G_u(f)_{C_v}\hookrightarrow \MTT(f)_{C_v}\rightarrow (\mathcal{H}_i)_{\sigma,C_v}\xrightarrow{(\phi_i)_{\sigma,C_v}}\GL((\mathcal{V}_{i})_{\sigma,C_v})$$
gives rise to a $C_v$-coefficient object $\mathscr{E}_{i,\sigma,v}$. Arguing as Lemma~\ref{lm:cptsystem}, for each $i\in \bbI$, $\{\mathscr{E}_{i,u}\}_{u\in\mathrm{fpl}(\bQ)}$ is a weakly $\bQ$-compatible system, and for each $\sigma\in \Sigma'_i$ the collection $\{\mathscr{E}_{i,\sigma,v}\}_{v\in\mathrm{fpl}(C)}$ is a weakly $C$-compatible system.    
\begin{lemma}\label{lm:GC_v}
Notation as above, we have 
 $G(\mathscr{E}_{i,u}, x)=(\Res_{N_i/\bQ}\mathcal{H}_i)_{\bQ_u}$, and $G(\mathscr{E}_{i,\sigma,v}, x)=(\mathcal{H}_i)_{\sigma,C_v}$. 
 \end{lemma}
\begin{proof}
This follows from Theorem~\ref{thm:MTinSingle} for each $f_i$.
\end{proof}
As a consequence, the projection $G_u(f)\rightarrow (\Res_{N_i/\bQ}\mathcal{H}_i)_{\bQ_u}$ is surjective. This gives rise to a surjection $G_u(f)^{\ad}\twoheadrightarrow (\Res_{N_i/\bQ}\mathcal{H}_i)_{\bQ_u}$. The induced map $G_u(f)^{\ad}\rightarrow \prod_{i\in\bbI}(\Res_{N_i/\bQ}\mathcal{H}_i)_{\bQ_u}= \MTT(f)_{\bQ_u}^{\ad}$ is an injection, and identifies $G_u(f)^{\ad}$ as the monodromy group of $\bigoplus_{i\in \bbI}\mathscr{E}_{i,u}$. 

\subsubsection{Reducing Conjecture~\ref{conj:MT} to Conjecture~\ref{conj:Ttlinear}} In the following, let $\bbJ$ be a subset of $\bbI$. We use $\sigma_i$ to denote an element in $\Sigma'_i$, and use $\sigma_\bbJ=(\sigma_{j})_{j\in \bbJ}$ to denote an element in $\prod_{j\in \bbJ} \Sigma'_{j}$. To compactify the notation, we define, for $u\in \mathrm{fpl}(Q)$ and $v\in \mathrm{fpl}(C)$ dividing $u$,
$$\mathscr{E}_{\bbJ,u}:= \bigoplus_{j\in \bbJ}\mathscr{E}_{j,u},\,\,\,\,G_{\bbJ,u}:= G(\mathscr{E}_{\bbJ,u},x).
$$$$\mathscr{E}_{\bbJ,\sigma_{\bbJ},v}:=\bigoplus_{j\in \bbJ}\mathscr{E}_{j,\sigma_{j},v},\,\,\,\,G_{\bbJ,\sigma_{\bbJ},v}:= G(\mathscr{E}_{\bbJ,\sigma_{\bbJ},v},x).$$
We will also write $G_{i,u}$ (\textit{resp}. $G_{i,\sigma_i,v}$) for $G(\mathscr{E}_{i,u},x)$ (\textit{resp}. $G(\mathscr{E}_{i,\sigma_{i},v},x)$). The group $G_{i,\sigma_i,v}$ is absolutely simple by Lemma~\ref{lm:GC_v}. It is also clear that for fixed $\bbJ$ and $\sigma_\bbJ$, $\{\mathscr{E}_{\bbJ,u}\}_{u\in \mathrm{fpl}(\bQ)}$ is a weakly $\bQ$-compatible system and  $\{\mathscr{E}_{\bbJ,\sigma_{\bbJ},v}\}_{v\in \mathrm{fpl}(C)}$ is a  weakly $C$-compatible system.

\begin{defn}
    Let $\{G_a\}_{a\in \bbA}$ be a finite collection of (algebraic) groups over a field. Suppose there is a subgroup $G\subseteq \prod_{a\in \bbA} G_a$ such that each projection $G\rightarrow G_a$ is an isomorphism. We will say that $G$ is an \textbf{iso-graph} over the collection $\{G_a\}_{a\in \bbA}$. 
\end{defn}

\begin{lemma}\label{lm:Gfgpstructure}
Take $\bbA=\bigsqcup_{i\in \bbI} (\{i\}\times \Sigma'_i)$. There is a partition $\bbA=\bigsqcup_{w\in \mathcal{W}}\bbA_w$, with the property 
$$(i,\sigma), (i,\sigma')\in \bbA_w \Leftrightarrow \sigma=\sigma',$$
such that for any $u\in\mathrm{fpl}(\bQ)$ and $v\in \mathrm{fpl}(C)$ dividing $u$,  $G_u(f)^{\ad}_{C_v}$ is a product of iso-graphs over sub-collections $\{(\mathcal{H}_i)_{\sigma,C_v}\}_{(i,\sigma)\in \bbA_w}$, where $w$ runs over $\mathcal{W}$.
\end{lemma}
\begin{proof}
First, fix one $u\in\mathrm{fpl}(\bQ)$ and $v\in \mathrm{fpl}(C)$ dividing $u$. Lemma~\ref{lm:GC_v} and the paragraph after that imply that $G_u(f)^{\ad}_{C_v}\subseteq \prod_{a\in \bbA} (\mathcal{H}_i)_{\sigma,C_v}$, and $G_u(f)^{\ad}_{C_v}$ surjects onto $\prod_{\sigma\in \Sigma'_i}(\mathcal{H}_i)_{\sigma,C_v}$ for each $i\in \bbI$. Since each $(\mathcal{H}_i)_{\sigma,C_v}$ is absolutely simple, the existence of partition $\bbA=\bigsqcup_{w\in \mathcal{W}}\bbA_w$ and the product of iso-graph claim follow from Goursat's lemma.

Now we need to show that the partition and the product of iso-graph claim is independent of $u$ and $v$ chosen. This is an easy consequence of Lemma~\ref{lm:compatible} and the weakly compatibility of the collections $\{\mathscr{E}_{\bbJ,u}\}_{u\in \mathrm{fpl}(\bQ)}$  and $\{\mathscr{E}_{\bbJ,\sigma_{\bbJ},v}\}_{v\in \mathrm{fpl}(C)}$ for any $\bbJ\subseteq \bbI$. 
\end{proof}
\begin{lemma}
  \label{lm:goodlemma2}
Notation as in Lemma~\ref{lm:Gfgpstructure}. If $\bm{A}_w$ contains $(i,\sigma_i)$, $(j,\sigma_j)$, then $\sigma_i N_i=\sigma_j N_j$.
\end{lemma}
\begin{proof}
Suppose the conclusion does not hold. Without loss of generality, we can assume that  $\sigma_i N_i\not\subseteq\sigma_j N_j$. Chebotarev's density theorem implies that there is a $g\in \Gal(C/\bQ)$, such that $g\sigma_i=\sigma_i$ but $g\sigma_j\neq \sigma_j$. Let $\bbJ=\{i,j\}$ and $\sigma_{\bbJ}=\{\sigma_i,\sigma_j\}$. Since $\{(i,\sigma_i),(j,\sigma_j)\}\subseteq \bm{A}_w$, $G_{\bbJ,\sigma_{\bbJ},v}$ is an iso-graph over $\{G_{i,\sigma_{i},v},G_{j,\sigma_{j},v}\}$ for any $v\in \mathrm{fpl}(C)$. A $g$-twist shows that $G_{\bbJ,g\sigma_{\bbJ},gv}$ is an iso-graph over $\{G_{i,\sigma_{i},gv},G_{j,g\sigma_{j},gv}\}$. As a result, we also have $(j,g\sigma_j)\in\bm{A}_w$. Contradiction. 
\end{proof}
\begin{lemma}\label{lm:goodlemma}
There is a partition $\bbI=\bigsqcup_{h\in \mathcal{H}}\bbI_h$, such that for any $u\in \mathrm{fpl}(\bQ)$, we have $G_u(f)^{\ad}=\prod_{h\in \mathcal{H}} G_{\bbI_h,u}$, where each $G_{\bbI_h,u}$ is an iso-graph over $\{G_{i,u}\}_{i\in \mathbf{I}_h}$. 
\end{lemma}
\begin{proof}
Fix one $u$. It suffices to show that for any $\bbJ=\{i,j\}\subseteq \bbI$, $G_{\bbJ,u}$ is either a product $G_{i,u}\times G_{j,u}$ for all $u$, or an iso-graph over $\{G_{i,u},G_{j,u}\}$ for all $u$. Once this is done, the lemma follows from induction and weak-compatibility. 

For any $\bbJ$, let $\bbA=(\{i\}\times \Sigma'_i) \sqcup (\{j\}\times \Sigma'_j)$. We can deduce from Lemma~\ref{lm:Gfgpstructure} that there is a partition $\bbA=\bigsqcup\bbA_w$ satisfying the property as described in \textit{loc.cit}, such that for $v\in \mathrm{fpl}(C)$ dividing $u$, $(G_{\bbJ,u})_{C_v}$ is a product of iso-graphs over the partition. Suppose that  $G_{\bbJ,u}$ is neither the product of $G_{i,u}$ and $G_{j,u}$, nor an iso-graph over them. Then there is a subset $\bbA_{w}$ which consists of two elements $(i,\sigma_i)$ and $(j,\sigma_j)$, as well as another subset $\bbA_{w'}$ that consists of a single element. Without loss of generality, we can assume this single element is $(i,\sigma_i')$. Let $\sigma_{\bbJ}=\{\sigma_i,\sigma_j\}$. By Chebotarev's density theorem, there exists a place $v'\in \mathrm{fpl}(C)$ (that may not divide $u$) and an element $g\in \Gal(C/\bQ)$ fixing $v'$, such that $g\sigma_i=\sigma_i'$. By $g$-twist, we see that $G_{\bbJ,g\sigma_{\bbJ},v'}$ is an iso-graph over 
$\{G_{i,\sigma_i',v'}, G_{j,g\sigma_j,v'}\}$. This violates the fact that $\bbA_{w'}$ consists of a single element. 
\end{proof}
\begin{corollary}\label{lm:conjimplicationMTTl}
   In Setup~\ref{setup2}, Conjecture~\ref{conj:MT} holds if Conjecture~\ref{conj:Ttlinear} holds. 
\end{corollary}
\begin{proof} 
By Lemma~\ref{lm:GC_v}, $\MTT(f)^{\ad}_{\bQ_p}= \prod_{i\in \bbI}(\Res_{N_i/\bQ}\mathcal{H}_i)_{\bQ_p}=\prod_{i\in \bbI} G_{i,p}$. By Lemma~\ref{lm:goodlemma}, $G_p(f)^{\ad}$ is a product of iso-graphs over subcollections of $\{G_{i,p}\}_{i\in \bbI}$. Now $G_p(f)^{\ad}\subseteq \MTT(f)^{\ad}_{\bQ_p}$. If Conjecture~\ref{conj:Ttlinear} holds, then 
$\MTT(f)^{\ad}_{\bQ_p}$ and $G_p(f)^{\ad}$ have the same unipotent with respect to $\mu^{-1}$. It then follows that $\MTT(f)^{\ad}_{\bQ_p}=G_p(f)^{\ad}$. We conclude by Proposition~\ref{prop:centercontain}. 
\end{proof}
\subsubsection{Reducing Conjecture~\ref{conj:Ttlinear} to $\#\bbI=2$}\label{subsub:redI2} Let $T_{f_i}$ be the rational cocharacter lattice of the smallest formal tori containing the image of $f_i^{/x}$. We have \begin{equation}\label{eq:TplusTf_i}
    T_{f}\subseteq \bigoplus_{i\in \bbI} T_{f_i}.
\end{equation} Let $\mathcal{X}_{f_i}\in \mathrm{MS}_{f_i}(\mathcal{A}_{g_i})$.
Theorem~\ref{thm:TlinSingle} implies that $\dim T_{f_i}=\dim \mathcal{X}_{f_i}$ for each $i$. Since $f$ separates factors, we have $\prod_{i\in \bbI} \mathcal{X}_{f_i}\in \mathrm{MS}_{f}(\mathcal{A}_{g_\bbI})$, so to show Conjecture~\ref{conj:Ttlinear} it suffices to show that  (\ref{eq:TplusTf_i}) is an equality. For this, one can show the contrapositive statement: \begin{claim}\label{cm:oneredu}
If (\ref{eq:TplusTf_i}) is not an equality, then $\prod_{i\in \bbI} \mathcal{X}_{f_i}\notin \mathrm{MS}_{f}(\mathcal{A}_{g_\bbI})$ ($f$ does not separate factors). 
\end{claim} 
By Theorem~\ref{Thm:Tatelocal}, $T_{f}$ can be identified as the (Lie algebra of the) unipotent of $G_p(f)$ with respect to $\mu^{-1}$. By Lemma~\ref{lm:goodlemma}, we see that $T_{f,x}=\prod_{h\in \mathcal{H}}T_{{\bbI_h}}$, where each $T_{{\bbI_h}}$ is an iso-graph over $\{T_{f_i}\}_{i\in \bbI_h}$. Therefore, if (\ref{eq:TplusTf_i}) is not an equality, then there exist two indices $i,j$, such that the projection of $T_{f}$ to $T_{f_i}\oplus T_{f_j}$ is an iso-graph over $\{T_{f_i}, T_{f_j}\}$. Let $\bbJ=\{i,j\}$. Consider $$f_{\bbJ}:X\xrightarrow{f} \mathscr{A}_{g_{\bbI},\Fpbar}^{\ord}\rightarrow \mathscr{A}_{g_{\bbJ},\Fpbar}^{\ord}=\mathscr{A}_{g_{i},\Fpbar}^{\ord}\times \mathscr{A}_{g_{j},\Fpbar}^{\ord}.$$
Then $T_{f_{\bbJ}}$ is an iso-graph over $\{T_{f_i}, T_{f_j}\}$. Therefore to show \textit{Claim}~\ref{cm:oneredu}, it suffices to show that  $\mathcal{X}_{f_i}\times \mathcal{X}_{f_j}\notin \mathrm{MS}_{f_\bbJ}(\mathcal{A}_{g_{\bbJ}})$. This amounts to showing Conjecture~\ref{conj:Ttlinear} for the case where $\#\bbI=2$.
\subsubsection{Proof of Conjecture~\ref{conj:Ttlinear}} As we have seen in \S\ref{subsub:redI2}, to prove Conjecture~\ref{conj:Ttlinear} in Setup~\ref{setup2}, it suffices to do the case when $\#\bbI=2$. In the following, let $\bbI=\{1,2\}$. From \S\ref{subsub:redI2}, we see that it suffices to show that 
\begin{claim}\label{claimmoremore}
    If $T_{f}$ is an iso-graph over $\{T_{f_1}, T_{f_2}\}$, then $\mathcal{X}_{f_1}\times \mathcal{X}_{f_2}\notin \mathrm{MS}_{f}(\mathcal{A}_{g_{\bbI}})$.
\end{claim}
In the following, we will always assume that $T_{f}$ is an iso-graph over $\{T_{f_1}, T_{f_2}\}$. Then for any $u\in \mathrm{fpl}(\bQ)$, $G_u(f)^{\ad}$ is an iso-graph over $\{\Hg(f_i)^{\ad}_{\bQ_u}\}_{i\in \bbI}$. For each $i\in \bbI$, apply Theorem~\ref{Thm:classMTgp} to $\Hg(f_i)$ to get $F_i,\tau_i,\mathfrak{p}_i$ and a group $\mathcal{G}_i/F_i$ which is either $\SO(\mathcal{Q}_i)$ or $\SUU(\Phi_i)$. Let $\Sigma_i$ be the set of embeddings of $F_i$ into $\bC$. It follows from Lemma~\ref{lm:goodlemma2} that $\tau_1F_1=\tau_2F_2$. Therefore, we can identify $(F_1,\tau_1,\mathfrak{p}_1,\Sigma_1)$ with $(F_2,\tau_2,\mathfrak{p}_2,\Sigma_2)$ and simply call them $(F,\tau,\mathfrak{p},\Sigma)$. Since $\mathcal{G}_{1}^{\ad}$ is a form of $\mathcal{G}_{2}^{\ad}$, we are in three possible situations: \begin{enumerate}
    \item  (Orthogonal case) both $\mathcal{G}_1$ and $\mathcal{G}_2$ are orthogonal,
    \item (Unitary case) both $\mathcal{G}_1$ and $\mathcal{G}_2$ are unitary, 
    \item (Mixed case)
one of $\mathcal{G}_1$ and $\mathcal{G}_2$ is a descent of $\SU(1,1)$, while the other is a descent of $\SO(2,1)$. Moreover, one of $\mathcal{X}_{f_1}$ and $\mathcal{X}_{f_2}$ is a unitary Shimura curve, while the other is an orthogonal Shimura curve. Now a unitary Shimura curve specially corresponds to an orthogonal Shimura curve. We can do a further modification as in Construction~\ref{constr:varsepa} to replace the unitary Shimura curve by an orthogonal one. Therefore, we end up in the orthogonal case. So we can assume that this case does not happen. 
\end{enumerate}

Fixing $C\subseteq \bC$ to be a sufficiently large Galois extension of $\bQ$ containing all complex embeddings of $F$. For each $i\in \bbI$ and $\sigma\in \Sigma$, let $\bM_{p,f_i}$ be the overconvergent $F$-isocrystal we constructed in \S\ref{subsub:cissd}. Let
$\mathscr{M}_{i,u}$ and $\mathscr{M}_{i,\sigma,v}$ be the coefficient objects that we constructed in \S\ref{subsub:cptsystem1}. Let $\bM_{p,f}=\bM_{p,f_1}\oplus \bM_{p,f_2}$, let $\mathscr{M}_{u}=\mathscr{M}_{1,u}\oplus \mathscr{M}_{2,u}$, and let 
$\mathscr{M}_{\sigma,v}=\mathscr{M}_{1,\sigma,v}\oplus \mathscr{M}_{2,\sigma,v}$. We know from Theorem~\ref{thm:MTinSingle} that $G(\mathscr{M}_{i,\sigma,v})=(\mathcal{G}_{i})_{\sigma,C_v}$. 

\begin{lemma}\label{lm:gooogogogogogogo}
In both orthogonal and unitary cases, $G(\bL_{\bbI,u,X_0})$ is an iso-graph over $\{\Hg(f_i)_{\bQ_u}\}_{i\in\bbI}$, inducing an isomorphism $b_u:\mathscr{M}_{1,u}\xrightarrow{\sim} \mathscr{M}_{2,u}$ whose base change to $C_v$ splits into isomorphisms $b_{\sigma,v}:\mathscr{M}_{1,\sigma,v}\xrightarrow{\sim}\mathscr{M}_{2,\sigma,v}$. \end{lemma}
\begin{proof}
We know that $G(\bL_{\bbI,u,X_0})^{\ad}=G_u(f)^{\ad}$ is an iso-graph over 
$\{\Hg(f_i)^{\ad}_{\bQ_u}\}_{i\in \bbI}$. Now we treat case by case.
\begin{enumerate}
    \item (Orthogonal case) Since the center of $\Hg(f_i)$ is finite and since $G(\bL_{\bbI,u,X_0})$ is connected, we deduce from Goursat's lemma that 
$G(\bL_{\bbI,u,X_0})$ is an iso-graph over $\{\Hg(f_i)_{\bQ_u}\}_{i\in\bbI}$. Now we claim that $\bM_{p,f_1}\simeq \bM_{p,f_2}$. It suffices to show that $\omega_x(\bM_{p,f_1})$ and  $\omega_x(\bM_{p,f_2})$ are isomorphic as $G(\bM_{p,f})$-reps. 

Let $\gr\bM_{p,f}$ be the graded object of $\bM_{p,f}$, and let $\gr_j\bM_{p,f}$ be the graded piece of slope $j$. Let's take a maximal torus $B\subseteq G(\bM_{p,f})$ that is of the form $\im \mu\times B_0$, where $B_0$ is a maximal torus of $G(\gr_0\bM_{p,f})$. Consider the adjoint action $B\acts T_{f}$. Since $T_{f}$ is an iso-graph over $\{T_{f_i}\}_{i\in\bbI}$, we see that $T_{f_1}$ is isomorphic to $T_{f_2}$ as $B$-reps. This soon implies that $\omega_x(\bM_{p,f_1})$ and  $\omega_x(\bM_{p,f_2})$ are isomorphic as $B$-reps, hence isomorphic as $G(\bM_{p,f})$-reps. 

Now for $\mathfrak{P}\in \mathrm{fpl}(C)$ lying above $\tau \mathfrak{p}$, we have $\mathscr{M}_{i,\tau,\mathfrak{P}}\simeq \bM_{p,f_i}\otimes C_{\mathfrak{P}}$. So  $\mathscr{M}_{1,\tau,\mathfrak{P}}$ and  $\mathscr{M}_{2,\tau,\mathfrak{P}}$ are isomorphic. Equivalently, they are $C$-compatible. Since $\{\mathscr{M}_{i,\sigma,v}\}_{v\in \mathrm{fpl}(C)}$ is a weakly $C$-compatible collection, we see that $\mathscr{M}_{1,\sigma,v}$ and $ \mathscr{M}_{2,\sigma,v}$ are $C$-compatible for all $v$. It follows that $\mathscr{M}_{1,u}$ and $\mathscr{M}_{2,u}$ are $\bQ$-compatible. By \cite[Theorem 4.2.11]{DA20}, there is a closed point of $X_0$ whose Frobenius torus is a maximal torus of $G_u(f)$. Therefore $\omega_x(\mathscr{M}_{1,u})$ and $\omega_x(\mathscr{M}_{2,u})$ are isomorphic as reps of this maximal torus, hence as $G_u(f)$-reps.  This implies the existence of the desired isomorphism. 
\item (Unitary case) To show that $G(\bL_{\bbI,u,X_0})$ is an iso-graph over $\{\Hg(f_i)_{\bQ_u}\}_{i\in\bbI}$, it suffices to show that $G(\bM_{p,f})$ is an iso-graph over $\{G(\bM_{p,f_i})\}_{i\in \bbI}$. 
Since $G(\bM_{p,f})^{\der}$ is already an iso-graph over $\{G(\bM_{p,f_i})^{\der}\}_{i\in \bbI}$, it suffices to show that $G(\bM_{p,f})$ does not contain the product of the central tori of $G(\bM_{p,f_1})$ and $G(\bM_{p,f_2})$. This is equivalent to showing that $G(\gr_{-1}\bM_{p,f}\oplus \gr_{1}\bM_{p,f})=\im \mu$, which follows from the fact that $T_{f}$ is invariant under the conjugation action of $G(\gr_{-1}\bM_{p,f}\oplus \gr_{1}\bM_{p,f})$. The isomorphisms between the coefficient objects are similar to the orthogonal case. 
\end{enumerate}
\end{proof}
\begin{theorem}\label{thm:MTTllinear}
   Conjecture~\ref{conj:Ttlinear} holds in Setup~\ref{setup2} (hence in Setup~\ref{setup1}).
\end{theorem}
\begin{proof}
 By \S\ref{subsub:redI2}, it suffices to show Conjecture~\ref{conj:Ttlinear} in Setup~\ref{setup2} with $\#\bbI=2$. For this, it suffices to show that \textit{Claim}~\ref{claimmoremore}. Let $\eta$ be the generic point $X_0$. By Lemma~\ref{lm:subHodgeT}, it suffices to show that the subgroup $U\subseteq U_{\GSpin',\mu^{-1}}$ that  
commutes with $\End(A_{\mathbf{I},\overline{\eta}})_{\bQ_p}$ has the property \begin{equation}\label{eq:dimcompare}
    \dim U<\dim (\mathcal{X}_{f_1}\times \mathcal{X}_{f_2}) =\dim {T}_{f_1}+\dim {T}_{f_2}.
\end{equation}
The key point is to construct a nontrivial homomorphism $A_{1,\overline{\eta}}\rightarrow A_{2,\overline{\eta}}$ using Kuga--Satake construction. In the following, we treat case by case. 
\begin{enumerate}
    \item(Orthogonal case) Let $\mathfrak{P}\in \mathrm{fpl}(C)$ denote a place lying above $\tau \mathfrak{p}$. Let $b_p:\mathscr{M}_{1,p} \xrightarrow{\sim}\mathscr{M}_{2,p}$ and $b_{\sigma,\mathfrak{P}}:\mathscr{M}_{1,\sigma,\mathfrak{P}} \xrightarrow{\sim}\mathscr{M}_{2,\sigma,\mathfrak{P}}$ be as in Lemma~\ref{lm:gooogogogogogogo}. Representation theoretically, we have isomorphisms of $G(\mathscr{M}_{p})$-reps $b_p:\omega_x(\mathscr{M}_{1,p}) \xrightarrow{\sim}\omega_x(\mathscr{M}_{2,p})$ and 
  $G(\mathscr{M}_{\sigma,\mathfrak{P}})$-reps $b_{\sigma,\mathfrak{P}}:\omega_x(\mathscr{M}_{1,\sigma,\mathfrak{P}}) \xrightarrow{\sim}\omega_x(\mathscr{M}_{2,\sigma,\mathfrak{P}})$. 
  
  Recall that $\omega_x(\bL_{i,p})=L_{i,\bQ_p}$ is equipped with a pairing $Q_i$ that comes from the GSpin Shimura datum. By Proposition~\ref{Thm:classMTgp} we have  $Q_i=\Tr_{F/\bQ}(\mathcal{Q}_i)\oplus Q^{\varrho}_i$. For each $\sigma\in \Sigma$, the pairing $\mathcal{Q}_i$ induces pairings $\mathcal{Q}_{i,\sigma,\mathfrak{P}}:=\mathcal{Q}_i\otimes_{\sigma} C_{\mathfrak{P}}$ on $\omega_x(\mathscr{M}_{i,\sigma,\mathfrak{P}})$, such that $G(\mathscr{M}_{i,\sigma,\mathfrak{P}})=\SO(\mathcal{Q}_{i,\sigma,\mathfrak{P}})$. Since $G(\mathscr{M}_{\sigma,\mathfrak{P}})$ is an iso-graph over $\{G(\mathscr{M}_{i,\sigma,\mathfrak{P}})\}_{i\in \bbI}$, we find by Lemma~\ref{lm:homothetyQ} that each $b_{\sigma,\mathfrak{P}}$ preserves the pairings up to homothety:  
  \begin{equation}\label{eq:er34sdf}
      \exists \lambda_{\sigma,\mathfrak{P}} \in C_{\mathfrak{P}}^*, \,\,\,\mathcal{Q}_{2,\sigma,\mathfrak{P}}(b_{\sigma,\mathfrak{P}}\,-,b_{\sigma,\mathfrak{P}}\,-)=\lambda_{\sigma,\mathfrak{P}} Q_{1,\sigma,\mathfrak{P}}(-,-).
  \end{equation} 
  The Galois group $\Gal(C_{\mathfrak{P}}/\bQ_p)$ acts on $b_{\sigma,\mathfrak{P}}$'s and $\lambda_{\sigma,\mathfrak{P}}$'s by $gb_{\sigma,\mathfrak{P}}=b_{g\sigma,\mathfrak{P}}$ and $g\lambda_{\sigma,\mathfrak{P}}=\lambda_{g\sigma,\mathfrak{P}}$. 
  
 Now we pass to spin representations to construct a morphism $\bH_{1,p,X_0} \rightarrow \bH_{2,p,X_0}$. Let $\Cl(-)$ be the Clifford algebra construction. First, for each $\sigma\in \Sigma$, we can construct an isomorphism of \textit{vector spaces} 
\begin{equation}\label{eq:eeqeqeclifoord}
\begin{aligned}
\theta_{\sigma,\mathfrak{P}}:\Cl(\omega_x(\mathscr{M}_{1,\sigma,\mathfrak{P}}))&\rightarrow \Cl(\omega_x(\mathscr{M}_{2,\sigma,\mathfrak{P}})),\\
s_1s_2...s_{m}&\rightarrow \lambda^{-\lfloor m/2\rfloor}_{\sigma,\mathfrak{P}} b_{\sigma,\mathfrak{P}}(s_1)b_{\sigma,\mathfrak{P}}(s_2)...b_{\sigma,\mathfrak{P}}(s_{m}).
   \end{aligned}
\end{equation}
Define an isomorphism of \textit{vector spaces} $\theta_{\mathfrak{P}}:\Cl(\omega_x(\mathscr{M}_{1,p}\otimes C_{\mathfrak{P}})\rightarrow \Cl(\omega_x(\mathscr{M}_{2,p}\otimes C_{\mathfrak{P}})$ by
\begin{equation}
\begin{tikzcd}
{\bigotimes_{\sigma\in \Sigma}\Cl(\omega_x(\mathscr{M}_{1,\sigma,\mathfrak{P}}))} \arrow[d, Rightarrow,no head] \arrow[rr, "{\theta_{\mathfrak{P}}=\bigotimes_{\sigma\in \Sigma}\theta_{\sigma,\mathfrak{P}}}"] &  & {\bigotimes_{\sigma\in \Sigma}\Cl(\omega_x(\mathscr{M}_{2,\sigma,\mathfrak{P}}))} \arrow[d, Rightarrow,no head] \\
{\Cl(\omega_x(\mathscr{M}_{1,p}\otimes C_{\mathfrak{P}}))}                                                                                                                                                       &  & {\Cl(\omega_x(\mathscr{M}_{2,p}\otimes C_{\mathfrak{P}}))}                                                      
\end{tikzcd}
\end{equation}
It is easy to check that $\theta_{\mathfrak{P}}$ is invariant under $\Gal(C_{\mathfrak{P}}/\bQ_p)$-action, so it descends to a map $$\theta_p: \Cl(\omega_x(\mathscr{M}_{1,p}))\rightarrow \Cl(\omega_x(\mathscr{M}_{2,p})).$$
Since $\bL_{i,p,X_0}$ is a direct sum of $\mathscr{M}_{1,p}$ with a trivial coefficient object $\bL_{i,p,X_0}^{\varrho}$ (say, of rank $r_i$), we see that $\omega_x(\bH_{i,p,X_0})=\Cl(\omega_x(\bL_{i,p,X_0}^{\varrho}))\otimes\Cl(\omega_x(\mathscr{M}_{i,p}))$ is a direct sum of $2^{r_i}$ copies of $\Cl(\omega_x(\mathscr{M}_{i,p}))$ (even as $G(\bH_{i,p,X_0})$-reps). Pick a copy of $\Cl(\omega_x(\mathscr{M}_{1,p}))$ in $\omega_x(\bH_{1,p,X_0})$ and a copy of $\Cl(\omega_x(\mathscr{M}_{2,p}))$ in $\omega_x(\bH_{2,p,X_0})$, then extend $\theta_p$ to a map $$\Theta_p:\omega_x(\bH_{1,p,X_0})\rightarrow \omega_x(\bH_{2,p,X_0})$$ by filling in 0 maps on the other copies.

We claim that $\Theta_p$ is a morphism of $G_p(f)$-reps. First, $\Spin(Q_{i,\sigma,\mathfrak{P}})$ acts on $\Cl(\omega_x(\mathscr{M}_{i,\sigma,\mathfrak{P}}))$, and we have an isomorphism $c_{\sigma}:\Spin(Q_{1,\sigma,\mathfrak{P}})\rightarrow \Spin(Q_{2,\sigma,\mathfrak{P}})$ which reduces to the isomorphism $G(\mathscr{M}_{1,\sigma,\mathfrak{P}}) )\rightarrow G(\mathscr{M}_{2,\sigma,\mathfrak{P}})$ whose graph is $G(\mathscr{M}_{\sigma,\mathfrak{P}})$. It is easy to check that (\ref{eq:eeqeqeclifoord}) is compatible with $c_\sigma$ and the spin representations. It then follows that $\theta_{p}$ is a map of $G_p(f)$-reps, hence $\Theta_p$ is a map of $G_p(f)$-reps.

This $\Theta_p$ then corresponds to an element $\bm{\delta}_{\bbI}\in\Hom(A_{1,\overline{\eta}},A_{2,\overline{\eta}})\otimes\bQ_p$ by crystalline Tate conjecture (\cite{DJ98}). By construction, the existence of such an element will force the subgroup $U$ in (\ref{eq:dimcompare}) to be contained in the iso-graph $T_f$. Therefore the inequality in (\ref{eq:dimcompare}) holds. This finishes the proof in the orthogonal case.  

\item (Unitary case) Let $\mathfrak{P}\in \mathrm{fpl}(C)$ denote a place lying above $\tau \mathfrak{p}$. Similar to the orthogonal case, we start by isomorphisms of $G(\mathscr{M}_{p})$-reps $b_p:\omega_x(\mathscr{M}_{1,p}) \xrightarrow{\sim}\omega_x(\mathscr{M}_{2,p})$ and 
  $G(\mathscr{M}_{\sigma,\mathfrak{P}})$-reps $b_{\sigma,\mathfrak{P}}:\omega_x(\mathscr{M}_{1,\sigma,\mathfrak{P}}) \xrightarrow{\sim}\omega_x(\mathscr{M}_{2,\sigma,\mathfrak{P}})$. 

We again show that $b_{\sigma,\mathfrak{P}}$ preserves the pairings up to homothety: Instead of using Lemma~\ref{lm:homothetyQ}, we use the fact that, possibly enlarging $C$, the representation of the unitary group $G(\mathscr{M}_{i,\sigma,\mathfrak{P}})$ on $\omega_x(\mathscr{M}_{i,\sigma,\mathfrak{P}})$ splits as a direct sum of two mutually dual absolutely irreducible subreps. This induces a split of $\mathscr{M}_{i,\sigma,\mathfrak{P}}$ as the direct sum of two mutually dual coefficient objects $\mathscr{W}_{i,\sigma,\mathfrak{P}}$ and $\mathscr{W}_{i,\sigma,\mathfrak{P}}'$. Up to possible relabeling (of the mutually dual summands), the $G(\mathscr{M}_{\sigma,\mathfrak{P}})$-rep  $b_{\sigma,\mathfrak{P}}:\omega_x(\mathscr{M}_{1,\sigma,\mathfrak{P}}) \xrightarrow{\sim}\omega_x(\mathscr{M}_{2,\sigma,\mathfrak{P}})$ splits as absolutely irreducible $G(\mathscr{M}_{\sigma,\mathfrak{P}})$-subreps $\beta_{\sigma,\mathfrak{P}}:\omega_x(\mathscr{W}_{1,\sigma,\mathfrak{P}}) \xrightarrow{\sim}\omega_x(\mathscr{W}_{2,\sigma,\mathfrak{P}})$ and  $\beta_{\sigma,\mathfrak{P}}':\omega_x(\mathscr{W}'_{1,\sigma,\mathfrak{P}}) \xrightarrow{\sim}\omega_x(\mathscr{W}'_{2,\sigma,\mathfrak{P}})$. By Schur's lemma,  $\beta_{\sigma,\mathfrak{P}}'$ is the inverse of the dual of $\beta_{\sigma,\mathfrak{P}}$ up to a constant. This soon implies (\ref{eq:er34sdf}).
  
The rest of the proof is similar to the orthogonal case: construct $\theta_{\sigma,\mathfrak{P}}$ as in (\ref{eq:eeqeqeclifoord}), then construct $\theta_\mathfrak{P}$ and descend to $\theta_{p}$, and then construct $\Theta_p$, which yields an element $\bm{\delta}_{\bbI}\in\Hom(A_{1,\overline{\eta}},A_{2,\overline{\eta}})\otimes\bQ_p$ that forces the inequality in (\ref{eq:dimcompare}). 
\end{enumerate}
    
\end{proof}
\begin{theorem}\label{thm:MTTllinear2}
    Conjecture~\ref{conj:MT} holds in Setup~\ref{setup2} (hence in Setup~\ref{setup1}).
\end{theorem}
\begin{proof}
Combine Theorem~\ref{thm:MTTllinear} and Corollary~\ref{lm:conjimplicationMTTl}.
\end{proof}
\begin{lemma}\label{lm:homothetyQ}
Let $V$ be a vector space over a characteristic 0 field $E$ such that $\dim_E V\neq 2$. Suppose that $Q_1,Q_2$ are two nondegenerate quadratic forms over $V$. If $\SO(Q_1)=\SO(Q_2)$ as subgroups of $\GL(V)$, then there is some $\lambda\in E^*$ such that $Q_2=\lambda Q_1$.  
\end{lemma}
\begin{proof}
    Combine \cite[Proposition C.3.14]{ReductiveGS} and Schur's lemma (note that $\SO(Q_i)\acts V$ is absolutely irreducible).
\end{proof}

\section{Characteristic $p$ analogue of the André--Oort conjecture}\label{sec:AO}
We begin by working with the general framework, and establish some general results (e.g. the geometric squeeze theorem~\ref{lm:thekeylemma}). After that, we specialize to the product of GSpin case and use our results in previous chapters on the Tate-linear and the mod $p$ Mumford--Tate conjectures to prove the mod $p$ André--Oort conjecture.  
\begin{setup}\label{setup4}
From Construction~\ref{constr:sepa} and Proposition~\ref{prop:AOquasilifts}, we know that to prove Conjecture~\ref{conj:AOAO}, it suffices to show Conjecture~\ref{conj:varAO}. So let $\bbI$ be a finite index set $\bbI$. For each $i\in \bbI$, let $\mathcal{A}_{g_i}$ be a Siegel modular variety with integral canonical model $\mathscr{A}_{g_i}$. Let $$f:X\hookrightarrow\mathscr{A}_{g_\bbI,\Fpbar}^{\ord}$$be a locally closed subvariety that contains a Zariski dense collection $\Xi$ of positive dimensional special subvarieties, such that $f$ \textbf{separates factors} over $\mathscr{A}_{g_\bbI,\Fpbar}^{\ord}$. Let $f_0:X_0\hookrightarrow \mathscr{A}_{g_\bbI,\Fpbar_q}^{\ord}$ be a model of $f$ over a sufficiently large $\Fpbar_q$. There are two more simplifications that do not affect the result: 
\begin{itemize}
    \item We can arbitrarily shrink $X$ to open dense, and replace $\Xi$ by its restriction. 
    \item   We can assume that each $Z\in \Xi$ is smooth and geometrically connected, by throwing away the singular loci, and re-indexing the components. 
\end{itemize}

\end{setup}



\subsection{Formal schemes}\label{sub:prepairformalscheme}
We will make heavy use of the language of formal schemes. This section is a summery of the background. The readers can consult \cite[\S 10]{EGA1} for a detailed treatment. An affinoid formal scheme is a locally topologically ringed space $$(\Spf R ,\mathcal{O}_{\Spf R}):=(\Spec R/I ,\varprojlim\mathcal{O}_{\Spec R/I^n}),$$ where $R$ is a commutative ring which is complete and separated with respect to an ideal $I$ (called the ideal of definition). An affinoid formal scheme is usually just denoted by $\Spf R$.  A formal scheme is a locally topologically ringed space $\mathfrak{X}:=(|\mathfrak{X}|,\mathcal{O}_{\mathfrak{X}})$ which is locally an affinoid formal scheme. 

An affinoid formal scheme $\Spf R$ is called Noetherian, if $R$ is Noetherian. A formal scheme $\mathfrak{X}$ is called Noetherian, if it is quasi-compact, and is locally  affinoid Noetherian. All formal schemes considered in this paper are Noetherian. 

A coherent ideal of a Noetherian formal scheme $\mathfrak{X}$ is a sheaf of ideals $\mathcal{J}\subseteq\mathcal{O}_{\mathfrak{X}}$. On a local affinoid $\Spf R$, it is the sheaf $J^{\Delta}$ associated to an ideal $J\subseteq R$. Let $R\la f^{-1}\ra:=\varprojlim R_f/I^n$ and let $D_f=\Spf R\la f^{-1}\ra$ be the associated dominant open. The sheaf $J^{\Delta}$ is defined by $$J^\Delta(D_f):= \varprojlim (J R_f)/I^n=JR\la f^{-1}\ra.$$ 
The second equality is true when $I$ is finitely generated, which is guaranteed by the Noetherian property. Following \cite[\S 10.14]{EGA1}, closed formal subschemes of a Noetherian formal scheme $\mathfrak{X}$ are defined by coherent ideals: if $\mathcal{J}$ is a coherent ideal, then the closed formal subscheme defined by $\mathcal{J}$ is $\mathfrak{Y}=(|\mathfrak{Y}|,(\mathcal{O}_\mathfrak{X}/\mathcal{J})|_{|\mathfrak{Y}|})$, where $|\mathfrak{Y}|\subseteq |\mathfrak{X}|$ is the support of $\mathcal{O}_\mathfrak{X}/\mathcal{J}$. An open formal subscheme of $\mathfrak{X}$ is the restriction of $(|\mathfrak{X}|,\mathcal{O}_{\mathfrak{X}})$ to an open subset. A locally closed formal subscheme is an open formal subscheme of a closed formal subscheme. A locally closed formal subscheme is again Noetherian.

\subsubsection{Irreducible components}
Let $\mathfrak{X}=\Spf R$ be an affinoid Noetherian formal scheme. We can talk about the formal subscheme of $\mathfrak{X}$ with reduced induced structure, as well as irreducible components: they are closed formal subschemes defined by the nilradical or minimal primes. Clearly, $\Spf R$ has finitely many irreducible components. Unlike the classical scheme theory, irreducibility or reducedness may not carry over to  open formal subschemes. 

\subsubsection{Base changes} Let $X$ be a scheme over $\Fpbar_p$, regarded as a formal scheme with ideal of definition $(0)$. Let $\mathfrak{X}\rightarrow X$ be a formal scheme over $X$. Let $Y\rightarrow X$ be a morphism of schemes. The base change of $\mathfrak{X}$ to $Y$ is the completed base change $\mathfrak{X}\widehat{\times}_X Y$, denoted by $\mathfrak{X}_Y$.

\subsubsection{Formal tori} Let $X$ be a  scheme over $\Fpbar_p$.  A formal torus $\mathfrak{T}$ of rank $d$ over $X$ is a formal scheme prorepresenting the sheaf over $X_{\mathrm{fppf}}$ associated to a multiplicative $p$-divisible group of rank $d$ over $X$. By \cite[Lemma 3.1.1, 3.1.2]{SW}, $\mathfrak{T}$ is a formal Lie group. More precisely, $\mathfrak{T}$ is isomorphic to $\Spf R[[t_1,t_2,...,t_d]]$ (as a formal scheme) when base changed along an Zariski affine cover $\Spec R\rightarrow X$. It is Noetherian when $X$ is.

The cocharacter lattice of $\mathfrak{T}$ is defined to be $X_*(\mathfrak{T}):=\underline{\mathrm{Hom}}_{\mathrm{group}/X}(\bG_{m}^{\wedge},\mathfrak{T})$, which is a lisse $\bZ_p$-sheaf over $X_{\et}$. Conversely, a lisse $\bZ_p$-sheaf $\mathscr{F}$ gives rise to a formal torus with cocharacter lattice $\mathscr{F}$. We denote the formal torus by $$\mathscr{F}\otimes_{\bZ_p} \bG_m^{\wedge}.
$$

A saturated lisse $\bZ_p$-subsheaf $\mathscr{F}'\subseteq \mathscr{F}$ gives rise to a subtorus $\mathscr{F}'\otimes_{\bZ_p} \bG_m^{\wedge}\subseteq \mathscr{F}\otimes_{\bZ_p} \bG_m^{\wedge}$ which is a closed formal subscheme.

\subsection{Global deformation spaces and Serre--Tate lisse sheaves}\label{subsec:speclisse} 
In this section we work with a single Siegel modular variety. The results generalize to products of Siegel modular varieties with obvious modifications. 

Let $\mathcal{A}_g$ be a Siegel modular variety attached to a self-dual symplectic $\bZ_{(p)}$-lattice $H$, with a sufficiently small hyperspecial level structure and integral canonical model $\mathscr{A}_{g}$. Let $q$ be a sufficiently large power of $p$. Let $\mathscr{G}=A[p^\infty]$ be the universal $p$-divisible group over $\mathscr{A}_{g,\Fpbar_q}^{\ord}$. Let $x\in \mathscr{A}_{g,\Fpbar_q}^{\ord}(\Fpbar)$, and let $\mu$ be the canonical Hodge cocharacter of $x$. We will identify $\omega(\mathscr{G}_x)$ with $H_{\bZ_{p}}$ via \S\ref{eq:comparison111}. The lisse sheaf $X^*(\mathscr{G}^{\loc})\oplus T_p(\mathscr{G}^{\et})^\vee$ corresponds to a representation 
\begin{equation*}
\mathrm{can}_{\mathscr{G}}:\pi_1^{\et}(\mathscr{A}_{g,\Fpbar_q}^{\ord},x)\rightarrow \GL(\omega(\mathscr{G}_x))=\GL(H_{\bZ_p})  
\end{equation*}
whose image lies in the Levi of $\GSp(H_{\bZ_p})$ fixing $\mu$. Let $U_{\GL,\mu^{-1}}, U_{\GSp, \mu^{-1}}$ be the opposite unipotents of $\GL(H_{\bZ_p}),\GSp(H_{\bZ_p})$ with respect to $\mu$. Then $\mathrm{can}_{\mathscr{G}}$ induces adjoint representations
\begin{align*}
  & \pi_1^{\et}(\mathscr{A}_{g,\Fpbar_q}^{\ord},x)\acts \Lie U_{\GL, \mu^{-1}},\\
 & \pi_1^{\et}(\mathscr{A}_{g,\Fpbar_q}^{\ord},x)\acts \Lie U_{\GSp, \mu^{-1}},
\end{align*}
The first one corresponds to the $\bZ_p$-lisse sheaf $\mathscr{E}=X_*({\mathscr{G}}^{\loc})\otimes_{\bZ_p} T_p({\mathscr{G}}^{\et})^\vee$, the second one corresponds to a saturated $\bZ_p$-lisse sheaf $\mathscr{C}\subseteq \mathscr{
E}$. The sheaf $\mathscr{C}$ can be interpreted via global canonical coordinates (\S\ref{sub:globalST}): Let $q_{\mathscr{G}}$ be the canonical pairing of $\mathscr{G}$. Then $\mathscr{C}$ can be identified with $\ker(q_\mathscr{G})^{\perp}$ (see the proof of Theorem~\ref{thm:etale-glob}). Let $\mathscr{C}_{\Fpbar}$ be the pullback of $\mathscr{C}$ to $\mathscr{A}_{g,\Fpbar}^{\ord}$. We 
have the following canonical isomorphism of formal schemes over $\mathscr{A}_{g,\Fpbar}^{\ord}$ (cf. \cite[\S1.2]{Ch03}): 
\begin{equation}\label{eq:globalserreTtae}
(\mathscr{A}_{g,\Fpbar}^{\ord}\times \mathscr{A}_{g,\Fpbar}^{\ord})^{/\Delta}\simeq \mathscr{C}_{\Fpbar}\otimes_{\mathbb{Z}_p}\mathbb{G}_{m}^{\wedge}.
\end{equation}
The left hand side should be thought of as a formal torus over $\mathscr{A}_{g,\Fpbar}^{\ord}$, whose fiber over each $x\in  \mathscr{A}_{g,\Fpbar}^{\ord}(\Fpbar)$ is the Serre--Tate formal torus $ \mathscr{A}_{g,\Fpbar}^{/x}$. The right hand side captures the cocharacter lattice of this formal torus. A proof of (\ref{eq:globalserreTtae}) can be sketched as follows: first, there is a tautological deformation of the polarized $p$-divisible group $A[p^{\infty}]$ to  $\mathscr{C}_{\Fpbar}\otimes_{\mathbb{Z}_p}\mathbb{G}_{m}^{\wedge}$, and by Serre--Tate theory, this gives a deformation of the polarized abelain scheme $A$. From moduli interpretation, we get a canonical morphism $\alpha:\mathscr{C}_{\Fpbar}\otimes_{\mathbb{Z}_p}\mathbb{G}_{m}^{\wedge}\rightarrow (\mathscr{A}_{g,\Fpbar}^{\ord}\times \mathscr{A}_{g,\Fpbar}^{\ord})^{/\Delta}$ which is an isomorphism when restricted to any $\Fpbar$-point of $\mathscr{A}_{g,\Fpbar}^{\ord}$. It is then easy to check that $\alpha$ is an isomorphism.

\subsubsection{Serre--Tate lisse sheaves}\label{subsub:gds} Let $X_0$ be smooth and geometrically irreducible and let $f_0:X_0\rightarrow \mathscr{A}_{g,\Fpbar}^{\ord}$ be any morphism. Let $f:X\rightarrow \mathscr{A}_{g,\Fpbar}^{\ord}$ be its base change to the algebraic closure. Let $x\in X_0(\Fpbar)$ and let $G=G(\bH_{p,X_0},x)$, where $\bH_p$ is the Dieudonné isocrystal of $\mathscr{G}$. Then the image of $$\can_{\mathscr{G}_{X_0}}:\pi_1^{\et}(X_0,x)\rightarrow \pi_1^{\et}(\mathscr{A}_{g,\Fpbar_q}^{\ord},x) \xrightarrow{\can_{\mathscr{G}}} \GL(H_{\bQ_p})=\GL(\omega_x(\bH_{p,X_0}))$$
lies in $\Lv_{G,\mu}$, the Levi of $G$ fixing $\mu$: this is a consequence of  Theorem~\ref{Thm:Tatelocal} and 
\cite[Theorem 2.1]{Crew1987}, which imply that $\Lv_{G,\mu}=G(\gr\bH_{p,X_0},x)$ is generated by $\im \mu$ and the Zariski closure of  ${\im \can_{\mathscr{G}_{X_0}}}$. Let $\mathcal{T}_{[X_0]}\subseteq \mathscr{C}_{X_0}$ be the saturated $\bZ_p$-lisse sheaf whose rationalization is the $\bQ_p$-lisse sheaf given by the adjoint representation
\begin{equation*}
    \pi_1^{\et}(X_0,x)\acts \Lie U_{G,\mu^{-1}}=T_{f,x}.
\end{equation*} 
We call $\mathcal{T}_{[X_0]}$ the \textbf{Serre--Tate lisse sheaf} on $X_0$. Let $\mathcal{T}_{[X]}$ be the pullback of $\mathcal{T}_{[X_0]}$ to $X$.

The space $\mathcal{T}_{[X]}\otimes_{\bZ_p} \bG_m^{\wedge}$ is a formal torus over $X$ whose fiber over each $y\in X(\Fpbar)$ is the smallest formal subtorus of $\mathscr{A}_{g,\Fpbar}^{/y}$ containing the image of $X^{/y}$. In other words, $\mathcal{T}_{[X],y}=\mathcal{T}_{f,y}$, where the later is a $\bZ_p$-lattice defined in (\ref{eq:ttlattices}). There are canonical morphisms of formal schemes over $X$:  
\begin{equation}\label{eq:canmorphismform}
    (X\times X)^{/\Delta} \rightarrow \mathcal{T}_{[X]}\otimes_{\bZ_p} \bG_m^{\wedge}\hookrightarrow \mathscr{C}_{X} \otimes_{\bZ_p} \bG_m^{\wedge} \simeq (X\times \mathscr{A}_{g,\Fpbar}^{\ord})^{/\Gamma_{f}}.
\end{equation}
Here the last isomorphism is from (\ref{eq:globalserreTtae}). If $f$ is a locally closed immersion, then the first morphism in (\ref{eq:canmorphismform}) is a locally closed immersion, as can be checked on fibers. If $X$ is furthermore a smooth special subvariety, then the first morphism in (\ref{eq:canmorphismform}) is an isomorphism. 

\subsubsection{Serre--Tate representations}\label{sub:mongpetale} 
Let $G$ be a reductive group over a characteristic 0 field $F$. Consider a cocharacter $\mu: \bG_m\rightarrow G$ defined over $F$. Denote by $\Lv_{G,\mu}$ and $U_{G,\mu^{-1}}$ the Levi and the opposite unipotent with respect to $\mu$. The adjoint representation
$$\ad_{G,\mu}:\Lv_{G,\mu} \acts \Lie U_{G,\mu^{-1}}$$
will be called the \textbf{Serre--Tate representation}. We have been using constructions of this kind in the previous chapters without giving it a name. Let $X,X_0,f,x$ be as in \S\ref{subsub:gds}. Let $G=G(\bH_{p,X_0},x)$, and let $\mu$ be the canonical Hodge cocharacter of $x$. By Theorem~\ref{Thm:Tatelocal}, $\ad_{G,\mu}$ can be identified with $$\ad_{f,x}:G(\gr\bH_{p,X_0},x)\acts T_{f,x}.$$ 
As we have noted at the beginning of \S\ref{subsub:gds}, $G(\gr\bH_{p,X_0},x)$ is generated by $\im \mu$ and the Zariski closure of  ${\im \can_{\mathscr{G}_{X_0}}}$. The restriction of $\ad_{f,x}$ to $\im \can_{\mathscr{G}_{X_0}}$ corresponds to the Serre--Tate lisse sheaf $\mathcal{T}_{[X_0]}$. The restriction of $\ad_{f,x}$ to $\im \mu$ is just  scalar multiplication. Therefore, a $\bQ_p$-lisse subsheaf $\mathscr{B}\subseteq \mathcal{T}_{[X_0]}\otimes \bQ_p$ corresponds canonically to a subrepresentation $\mathscr{B}_x\subseteq T_{f,x}$ of $\ad_{f,x}$.

\subsection{The geometric squeeze theorem}\label{subsec:keylemma}
Use Setup~\ref{setup4}. For each $i\in \bbI$, one has a Dieudonné isocrystal $\bH_{i,p}$ and lisse sheaves $\mathscr{E}_i,\mathscr{C}_i$ over $\mathscr{A}_{g_i,\Fpbar_q}^{\ord}$ as defined in \S\ref{subsec:speclisse}. For $\bbJ\subseteq\bbI$, let $\mathbb{H}_{\bbJ,p},\mathscr{E}_{\bbJ},\mathscr{C}_\bbJ$ be the direct sum of the pullback of various $\mathbb{H}_{i,p},\mathscr{E}_i,\mathscr{C}_i$ over $\mathscr{A}_{g_\bbJ,\Fpbar_q}^{\ord}=\prod_{i\in \bbJ} \mathscr{A}_{g_i,\Fpbar_q}^{\ord}$.

By \S\ref{subsub:gds}, $X_0$ is equipped with a Serre--Tate lisse sheaf 
$$\mathcal{T}_{[X_0]}\subseteq \mathscr{C}_{\bbI,X_0}.$$
Let $x\in X_0(\Fpbar)$ and let $G=G(\bH_{\bbI,p,X_0},x)^\circ$. Consider the Serre--Tate representation $\ad_{G,\mu}$ as per \S\ref{sub:mongpetale}. 
In this section, we prove the following: 
\begin{theorem}[{The geometric squeeze theorem}]\label{lm:thekeylemma} 
Let the setup be as in \ref{setup4}. Possibly shrinking $X$, there exists a nonzero saturated $\bZ_p$-lisse subsheaf $\mathcal{R}_f$ of $\mathcal{T}_{[X_0]}$ over a Galois cover $X_0'\rightarrow X_0$, such that for any $x\in X(\Fpbar)$, the following are true:
\begin{enumerate}
\item\label{it:squeeze1} 
${R}_{f,x}:=(\mathcal{R}_{f,x})_{\bQ_p}\subseteq T_{f,x}$ is a subspace invariant under $\ad_{G,\mu}$ (see Remark~\ref{rmk:squeeze} for our convention for $\mathcal{R}_{f,x}$).
  \item\label{it:squeeze2} $X^{/x}$ is sandwiched between two formal tori: $$\mathcal{R}_{f,x}\otimes \bG_m^\wedge\subseteq X^{/x}\subseteq \mathcal{T}_{f,x}\otimes \bG_m^\wedge.$$
   \item\label{it:squeeze3}
There is a $Z\in \Xi$, such that for $x\in Z(\Fpbar)$:
$${T}_{f|_Z,x}\subseteq {R}_{f,x}\subseteq T_{f,x}.$$ 
Furthermore, for each $i\in \bbI$, the projection image of $Z$ in $\mathscr{A}_{g_i,\Fpbar}$ is positive dimensional if and only if the projection image of $\mathcal{R}_{f}$ in $\mathscr{E}_{i,X_0'}$ has positive rank.
\end{enumerate}
\end{theorem}\begin{remark}\label{rmk:squeeze}    
Let's briefly explain the theorem: 
\begin{enumerate}[label=(\alph*)]
\item  We mean by $\mathcal{R}_{f,x}$ the fiber of $\mathcal{R}_{f}$ at an $\Fpbar$-point of $X_0'$ that lifts $x$, identified as a sublattice of $\mathcal{T}_{f,x}$. Different lifts of $x$ may yield different $\mathcal{R}_{f,x}$ as sublattices of $\mathcal{T}_{f,x}$, but the theorem works for all of them.
\item\label{rmkit:squeeze2} The two sandwiches  (\ref{it:squeeze2}) and (\ref{it:squeeze3}) make Theorem~\ref{lm:thekeylemma} resemble the classical squeeze theorem in calculus. A typical application is as follows: 

If $\ad_{G,\mu}$ is irreducible, then ${R}_{f,x}={T}_{f,x}$ by (\ref{it:squeeze1}) and $X^{/x}$ is squeezed to a formal torus by (\ref{it:squeeze2}). This shows that $X$ is Tate-linear. 
    \item Let $Z$ be as in (\ref{it:squeeze3}). Let $Z_0$ be a suitable finite model of $Z$ and let $G':=G(\bH_{\bbI,p,Z_0},x)^\circ$. Then the sandwich in (\ref{it:squeeze3}) is indeed a sandwich of $\Lv_{G',\mu}$-reps, where $\Lv_{G',\mu}$ acts on ${T}_{f|_Z,x}$. This indicates that $\mathcal{R}_{f}$ is big and is squeezed between ``representations of geometric origin''.  If  (\ref{it:squeeze2}) fails to be an equality, then it ought to fail for geometric reasons. See Remark~\ref{rmk:exceptionalsplit}.
\end{enumerate}
\end{remark}

We will prove Theorem~\ref{lm:thekeylemma} in the rest of this section.


\subsubsection{Some formal schemes}\label{subsub:someformalsubschemes} From (\ref{eq:globalserreTtae}), we have an ambient formal torus $\mathfrak{T}$ over $X$: 
\begin{equation}\label{eq:ambient30}
    \mathfrak{T}=(X\times \mathscr{A}_{g_\bbI,\Fpbar})^{/\Gamma_f}=\mathscr{C}_{\bbI,X}\otimes_{\mathbb{Z}_p}\mathbb{G}_{m}^{\wedge}.
\end{equation}
Possibly shrinking $X$, and replacing $\Xi$ by its restriction, we may and do assume that $X=\Spec B$ is affine, and   \begin{equation}\label{eq:ambient31}
    \mathfrak{T}\simeq \Spf B[[t_1,...,t_d]]
\end{equation}
for some $d$. We will carry this assumption over for the rest of this chapter.

Let $Z\in \Xi$, then from \S\ref{subsub:gds} there is a finite field model $Z_{0}\subseteq X_{0,\bF_{q^n}}$ of $Z$ over some extension $\bF_{q^n}/\bF_q$, and a Tate linear lisse subsheaf $$\mathcal{T}_{[Z_0]}\subseteq (\mathcal{T}_{[X_0]})_{Z_0}\subseteq\mathscr{C}_{\bbI,Z_{0}},$$ such that $$
(Z\times Z)^{/\Delta}=\mathcal{T}_{[Z]}\otimes_{\mathbb{Z}_p}\mathbb{G}_{m}^{\wedge},$$ where $\mathcal{T}_{[Z]}$ is the pullback of $\mathcal{T}_{[Z_0]}$ to $Z$. When $Z$ runs over the Zariski dense collection, the number $n$ may not be bounded. This is the reason why we work over $\Fpbar$ rather than a fixed $\Fpbar_q$. The following diagram summarizes the relation between spaces:
\begin{center}

\begin{tikzcd}
{\mathcal{T}_{[Z]}\otimes_{\mathbb{Z}_p}\mathbb{G}_{m}^{\wedge}} \arrow[d, Rightarrow, no head] \arrow[r, hook] & (X\times X)^{/\Delta} \arrow[d, Rightarrow,no head] \arrow[r, hook] & {\mathscr{C}_{\bbI,X}\otimes_{\bZ_p} \bG_m^{\wedge}} \arrow[d, Rightarrow,no head] \\
(Z\times Z)^{/\Delta} \arrow[r, hook] \arrow[d]                                                        & (X\times X)^{/\Delta} \arrow[d] \arrow[r, hook]             & (X\times \mathscr{A}_{g_\bbI,\Fpbar})^{/\Gamma_f} \arrow[d]                \\
Z \arrow[r, hook]                                                                                      & X \arrow[r, Rightarrow,no head]                            & X                                                                
\end{tikzcd}
\end{center}

For every $Z\in \Xi$, $(Z\times Z)^{/\Delta}$ is a locally closed formal subscheme of $\mathfrak{T}$. By Noetherian property, there is a smallest closed formal subscheme $\mathfrak{Z}_0\subseteq \mathfrak{T}$ containing all $(Z\times Z)^{/\Delta},Z\in \Xi$ (here by ``contain'' we mean $\mathfrak{Z}_0$ admits $(Z\times Z)^{/\Delta}$ as a locally closed formal subscheme). Equivalently, one can define $\mathfrak{Z}_0$ to be the smallest closed formal subscheme that contains all closed formal subschemes of the form $\mathcal{T}_{[Z],x}\otimes_{\bZ_p} \bG_m^\wedge$, where $Z\in \Xi$ and $x\in Z(\Fpbar)$. 

Since $\mathfrak{Z}_0$ only has finitely many irreducible components, it admits, by pigeon hole, an irreducible component which contains all $(Z\times Z)^{/\Delta},Z\in \Xi'$, where $\Xi'$ is a Zariski dense sub-collection of $\Xi$. We will call this irreducible component $\mathfrak{Z}$.

Note that $(X\times X)^{/\Delta}$ is a closed formal subscheme of $\mathfrak{T}$ containing all $(Z\times Z)^{/\Delta}, Z\in \Xi$. By minimality, we have 
\begin{equation}\label{eq:containZX}
\mathfrak{Z}\subseteq \mathfrak{Z}_0\subseteq  (X\times X)^{/\Delta}.    
\end{equation}

\subsubsection{Rigidity and big lisse sheaves} 
The torus $\mathfrak{T}$ admits a scaling automorphism $[\alpha]$ fo each $\alpha\in \bZ_p^*$. The following theorem is of independent interest, and will be proved in Appendix~\ref{subsub:Starshape}.
\begin{theorem}[Rigidity theorem over a geometric base]\label{lm:rigidity}
Let $X=\Spec B$ be a connected smooth affine variety over $\Fpbar$. Let $\mathfrak{T}$ be a formal torus over $X$, which is isomorphic to $\Spf B[[t_1,...,t_d]]$ as a formal scheme over $B$. Suppose that $\mathfrak{Z}$ is an irreducible closed formal subscheme of $\mathfrak{T}$ such that \begin{enumerate}
    \item $\mathfrak{Z}_x\neq \emptyset$ for Zariski dense $x\in X(\Fpbar)$.
    \item $\mathfrak{Z}$ is invariant under the scaling automorphism $[\alpha]$ for infinitely many $\alpha\in \bZ_p^*$.  
\end{enumerate}
Then there is a Galois cover $X'\rightarrow X$, such that $\mathfrak{Z}_{X'}$ is a union of formal subtori of $\mathfrak{T}_{X'}$ which are permuted transitively by $\Gal(X'/X)$.
\end{theorem}

\begin{corollary}\label{prop:locys}
Possibly base change $X_0$ to a larger $\Fpbar_q$, there is a Galois cover $g_0:X'_0\rightarrow X_0$ between geometrically connected $\Fpbar_q$-varieties, and a saturated $\bZ_p$-lisse subsheaf $\mathscr{H}\subseteq \mathcal{T}_{[X_0']}$ over $X_0'$, such that 
\begin{enumerate}
   \item\label{it:1} Let $g:X'\rightarrow X$ be the base change of $g_0$ to $\Fpbar$. Then as subspaces of $(X'\times \mathscr{A}_{g_\bbI,\Fpbar})^{/\Gamma_{f\comp g}}=\mathscr{C}_{\bbI,X'}\otimes_{\mathbb{Z}_p}\mathbb{G}_{m}^{\wedge}$, we have 
   $$\mathscr{H}_{X'}\otimes_{\bZ_p} \bG_m^{\wedge}\subseteq (X'\times X)^{/\Gamma_g}.$$
    \item\label{it:2} There is a $Z\in\Xi$ and a geometrically connected smooth subvariety $Z'_0\subseteq X'_0$ whose base change to $\bF$ lifts $Z$, such that $$\mathcal{T}_{[Z'_0]}\subseteq \mathscr{H}_{Z'_0}.$$    
    Furthermore, for each $i\in \bbI$, the projection image of $Z$ in $\mathscr{A}_{g_i,\Fpbar}$ is positive dimensional if and only if the projection image of $\mathscr{H}$ in $\mathscr{E}_{i,X_0'}$ has positive rank.
\end{enumerate}
\end{corollary}
\proof  
Notation as in \S\ref{subsub:someformalsubschemes}. For each $\alpha\in \bZ_p^*$ and each $Z\in \Xi$, the scaling automorphism $[\alpha]$ takes $(Z\times Z)^{/\Delta}=\mathcal{T}_{[Z]}\otimes_{\bZ_p}\bG_m^\wedge$ to itself. So $[\alpha]$ takes $\mathfrak{Z}_0$ to itself by minimality of $\mathfrak{Z}_0$. It then follows that $[\alpha]$ takes each irreducible component of $\mathfrak{Z}_0$ to itself. In particular, $[\alpha]$ takes $\mathfrak{Z}$ to itself.

Since there is a Zariski dense sub-collection $\Xi'\subseteq\Xi$, such that $\mathfrak{Z}$ contains all $(Z\times Z)^{/\Delta}$ for $Z\in \Xi'$, we see that $\mathfrak{Z}_x\neq \emptyset$ for Zariski dense $x\in X(\Fpbar)$. Therefore $\mathfrak{Z}$ meets the conditions of Theorem~\ref{lm:rigidity}. As a result, there is a Galois cover $g:X'\rightarrow X$, such that $\mathfrak{Z}_{X'}$ is a union of formal subtori $\mathfrak{T}_{1,X'},\mathfrak{T}_{2,X'},...,\mathfrak{T}_{l,X'}\subseteq \mathfrak{T}_{X'}$ of the same rank.

Each $\mathfrak{T}_{i,X'}$ gives rise to a saturated 
lisse subsheaf $X_*(\mathfrak{T}_{i,X'})\subseteq \mathscr{C}_{\bbI,X'}$. Pick a $Z\in \Xi'$. Take a connected component of $Z\times_X X'$ and denote it by $Z'$. Possibly enlarging $q$ and base changing $X_0$ to this new $\Fpbar_q$, we have geometrically connected $\bF_q$-varieties $Z_0',X_0'$ such that \begin{enumerate}[label=(\alph*)]
   \item $Z'_0\subseteq X'_0$, and there is a Galois covering map $g_0:X'_0\rightarrow X_0$, such that $Z'\subseteq X'$ \textit{resp}. $g:X'\rightarrow X$ is the base change from $\Fpbar_q$ to $\Fpbar$ of the map $Z_0'\subseteq X'_0$ \textit{resp}. $g_0:X_0'\rightarrow X_0$,
    \item $Z_0'(\bF_{q})\neq \emptyset$.
\end{enumerate} 
Let $\mathcal{T}_{[Z_0']}$ be the Serre--Tate lisse sheaf over $Z_0'$. Pick a point $x_0'\in Z_{0}'(\bF_{q})$ and a point $x'\in Z_{0}'(\Fpbar)$ lifting $x_0'$. Then the stalk $\mathcal{T}_{[X_0'],x'}$ is a continuous $\pi_1(X'_{0},x')$-module. Let ${\mathscr{H}}_{x'}$ be the saturated submodule of $\mathcal{T}_{[X_0'],x'}$ generated by the set \begin{equation}\label{eeeeee}
    \{g\cdot \mathcal{T}_{[Z_0'],x'}|g\in \pi_1(X',x')\}.
\end{equation} Since $\mathcal{T}_{[Z_0'],x'}$ is invariant under $\Gal(x'|x_0')$-action, we can use the fact that
$\pi_1(X'_{0},x')=\pi_1(X',x')\rtimes\Gal(x'|x_0')$ to show that $\mathscr{H}_{x'}$ is invariant under $\pi_1(X'_{0},x')$-action. It is then routine to check that $\mathscr{H}_{x'}$ is a continuous $\pi_1(X'_{0},x')$-representation, which thus gives rise to a saturated $\bZ_p$-lisse subsheaf $\mathscr{H}\subseteq\mathcal{T}_{[X_0']}$ that by definition satisfies (\ref{it:2}). 


Since $(Z\times Z)^{\Delta}$ is contained in $\mathfrak{Z}$, $\mathcal{T}_{[Z']}\otimes_{\bZ_p} \bG_m^\wedge$ is contained in $\mathfrak{Z}_{X'}$. Therefore $\mathcal{T}_{[Z']}\otimes_{\bZ_p} \bG_m^\wedge$ is contained in one of the formal tori in $\{\mathfrak{T}_{i,X'}\}_{i=1}^l$. Say, $\mathcal{T}_{[Z']}\otimes_{\bZ_p} \bG_m^\wedge$ is contained in $\mathfrak{T}_{1,X'}$. Then $\mathcal{T}_{[Z']}$ is contained in the restriction of $X_*(\mathfrak{T}_{1,X'})$ to ${Z'}$. It follows from definition (\ref{eeeeee}) that ${\mathscr{H}}_{X'}\subseteq X_*(\mathfrak{T}_{1,X'})$. Since $\mathfrak{Z}\subseteq (X\times X)^{/\Delta}$ by (\ref{eq:containZX}), we deduce that 
\begin{equation}\label{eq:222}
{\mathscr{H}}_{X'}\otimes_{\bZ_p} \bG_{m}^{\wedge} \subseteq X_*(\mathfrak{T}_{1,X'})\otimes_{\bZ_p} \bG_m^{\wedge}= \mathfrak{T}_{1,X'}\subseteq \mathfrak{Z}_{X'}\subseteq (X'\times X)^{/\Gamma_g}.
\end{equation}This proves (\ref{it:1}). $\hfill\square$
\begin{proof}[Proof of Theorem~\ref{lm:thekeylemma}]
Let $g_0:X_0'\rightarrow X_0$, $Z_0'\subseteq X_0'$ and $\mathscr{H}$ be as in Lemma~\ref{prop:locys}. 
Define $\mathcal{R}_{f}:=\mathscr{H}$. 
We need to check that this construction satisfies Theorem~\ref{lm:thekeylemma}(\ref{it:squeeze1})$\sim$(\ref{it:squeeze3}).


Let $x\in X_0(\Fpbar)$ and let $x'\in X_0'(\Fpbar)$ be a lift of $x$. Let $\mathcal{R}_{f,x}$ be the fiber of $\mathcal{R}_f$ at $x'$. Identify $\pi_1^{\et}(X_0',x')$ as a subgroup of $\pi_1^{\et}(X_0,x)$. Identify $G(\bH_{\bbI,p,X'_0},x')$ as a  subgroup of $G(\bH_{\bbI,p,X_0},x)$. 

Since $\mathscr{H}\subseteq \mathcal{T}_{[X_0']}$, we see that $R_{f,x}\subseteq (\mathcal{T}_{[X_0'],x'})_{\bQ_p}=T_{f,x}$ is a subspace invariant under $\pi_1^{\et}(X_0',x')$. As noted in \S\ref{sub:mongpetale}, this means that ${R}_{f,x}\subseteq T_{f,x}$ is $G(\gr\bH_{\bbI,p,X'_0},x')^\circ$-equivariant. Since $G=G(\bH_{\bbI,p,X_0},x)^\circ=G(\bH_{\bbI,p,X'_0},x')^\circ$, we see that ${R}_{f,x}\subseteq T_{f,x}$ is also $\Lv_{G,\mu}$-equivariant, this proves Theorem~\ref{lm:thekeylemma}(\ref{it:squeeze1}). 

Theorem~\ref{lm:thekeylemma}(\ref{it:squeeze2}) is a direct consequence of Lemma~\ref{prop:locys}(\ref{it:1}). Theorem~\ref{lm:thekeylemma}(\ref{it:squeeze3}) follows from  Lemma~\ref{prop:locys}(\ref{it:2}) by a similar argument that we used to show Theorem~\ref{lm:thekeylemma}(\ref{it:squeeze1}). \end{proof}

\subsection{Proof of the conjecture}\label{subsec:AOGSpin}
Now we specialize to the product of GSpin case (Setup~\ref{setup1} or \ref{setup2} with $f$ being an immersion). We begin with the single case.
\begin{theorem}\label{thm:AOsingle}
 Conjecture~\ref{conj:varAO} is true in Setup~\ref{setup2} with $f$ being an immersion and $\#\bbI=1$. In particular, $X$ is already special. 
\end{theorem}
\proof Take $\bbI=\{1\}$, and we will drop it from subscript. We need to show that $X$ is special. The case $\dim X=1$ is trivial. So we can assume that $\dim X\geq 2$. If we can show that $X$ is Tate-linear, then Theorem~\ref{thm:TlinSingle} guarantees that $X$ is special.

In the following, we can assume that our $X$ is sufficiently shrunk. Take $\mathcal{R}_{f}$ as in Theorem~\ref{lm:thekeylemma}. Let $x\in X(\Fpbar)$ be a point. By Remark~\ref{rmk:squeeze}(\ref{rmkit:squeeze2}), the Tate-linearity of $X$ follows from the irreducibility of $\ad_{G,\mu}$. We claim that $\ad_{G,\mu}$ is irreducible when $\dim T_{f,x}\geq 3$. This can be deduced from the characteristic $p$ Mumford--Tate conjecture~\ref{thm:MTinSingle}, which says that we can identify $G$ with $\MTT(f)_{\bQ_p}$. Then the irreducibility of $\ad_{G,\mu}$ is equivalent to the 
irreducibility of $$\ad_{\Hg(f)_{\bQ_p},\mu}:\Lv_{\Hg(f)_{\bQ_p},\mu}\acts T_{f,x}.$$ The structure of $\Hg(f)$ is well understood from Proposition~\ref{Thm:classMTgp}. In fact, $G(\mathbb{M}_{p,f},x)^\circ$ is the unique factor of $\Hg(f)_{\bQ_p}$ through which $\mu$ factors nontrivally. And the irreducibility of $\ad_{\Hg(f)_{\bQ_p},\mu}$ is equivalent to the irreducibility of $$\ad_{G(\mathbb{M}_{p,f},x)^\circ,\mu}:\Lv_{G(\mathbb{M}_{p,f},x)^\circ,\mu}\acts T_{f,x}.$$
Now $\Lv_{G(\mathbb{M}_{p,f},x)^\circ,\mu}= G(\gr\mathbb{M}_{p,f},x)^\circ$. The restriction of $\ad_{G(\mathbb{M}_{p,f},x)^\circ,\mu}$ to $G(\gr_0\mathbb{M}_{p,f},x)$ is the standard representation on $T_{f,x}$ of either $\SO(T_{f,x})$ or $\SU(\varphi_{f,0})$, as observed in Remark~\ref{eq:mtbaby}. If $\dim T_{f,x}\geq 3$, then the action of $\SO(T_{f,x})$ or $\SU(\varphi_{f,0})$ on $T_{f,x}$ is irreducible. This forces $X$ to be Tate-linear.

The remaining case is when $\dim T_{f,x}\leq 2$. But then we have $\dim X\leq \dim T_{f,x}\leq 2$. Since we already assumed that $\dim X\geq 2$ from the beginning, it follows that $\dim X=\dim T_{f,x}=2$. So $X$ is Tate-linear. This finishes the proof of the theorem. $\hfill\square$
\begin{remark}\label{rmk:exceptionalsplit}
  Consider the case $\dim T_{f,x}=2$ in the above proof. We claim that, even though $\ad_{G,\mu}$ may not be irreducible, we still have ${R}_{f,x}=T_{f,x}$. The reason is as follows. As observed in the last paragraph of the above proof, we have $\dim X=2$. Suppose that $\dim{R}_{f,x}=1$. Let $Z\in \Xi$ be the special subvariety in
   Theorem~\ref{lm:thekeylemma}(\ref{it:squeeze3}). Pick $y\in Z(\Fpbar)$. Then $\dim Z=1$ and $T_{f|_Z,y}={R}_{f,y}$. This forces $G(\bM_{p,f},x)^\circ$ to split. We must have  $\Hg(X)=\Res_{F/\bQ}\SO(\mathcal{Q})$, where $\mathcal{Q}$ has signature $(2,2)$ at a real place. Let $G'$ be the image of $\Hg(Z)_F\subseteq\Hg(X)_F\rightarrow \SO(\mathcal{Q})$. Then $G'$ is a proper normal subgroup of $\SO(\mathcal{Q})$ (which can be checked after base change to $F_\mathfrak{p}$). So $\SO(\mathcal{Q})$ splits. This contradicts our assumption that $f$ separates factors (i.e. $\Hg(X)$ is almost simple). 
\end{remark}

\begin{theorem}\label{thm:AOprod} Conjecture~\ref{conj:varAO} is true in Setup~\ref{setup2} with $f$ being an immersion. Therefore Conjecture~\ref{conj:AOAO} is true in Setup~\ref{setup1} with $f$ being an immersion. 
\end{theorem}
\begin{proof}
It suffices to show the first assertion, and the second follows from the first assertion and  Construction~\ref{constr:varsepa} and Proposition~\ref{prop:AOquasilifts}. In the following, we assume that our $X$ is sufficiently shrunk. Let $\mathcal{R}_{f}$ be as in Theorem~\ref{lm:thekeylemma}. Let $\bbI_s\subseteq \bbI$ be the set of indices consisting of those $i$ such that the projection of $\mathcal{R}_{f}$ to $\mathscr{E}_{i,X_0'}$ is not 0. Project $X$ to $\prod_{i\in \bbI_s}\mathscr{X}_{f_{i},\Fpbar}^\ord$ and let $C_s$ be a component of the later space containing the image. Project $X$ to $\prod_{i\in \bbI\setminus \bbI_s}\mathscr{X}_{f_{i},\Fpbar}^\ord$ and let $Y$ be the Zariski closure of the image in the later space. We have\begin{equation}\label{eq:proddd}
\overline{X}\subseteq C_s\times Y.
\end{equation}
It suffices to show that (\ref{eq:proddd}) is an equality. Let $x\in X(\Fpbar)$. By Theorem~\ref{thm:MTTllinear} or \ref{thm:MTTllinear2}, we have \begin{equation}
    T_{f,x}=\bigoplus_{i\in \bbI} T_{f_i,x}.
\end{equation}
 
\begin{claim}\label{cm:dfq}
 ${R}_{f,x}=\bigoplus_{i\in \bbI_s}T_{f_{i},x}$.
\end{claim}
By our assumption, we have 
\begin{equation}\label{eq:eredsd}
    {R}_{f,x}\subseteq\bigoplus_{i\in \bbI_s}T_{f_{i},x}.
\end{equation} 
Let $G=G(\bH_{\bbI,p,X_0},x)^\circ$ and let $G_i=G(\bH_{i,p,X_0},x)^\circ$. By Theorem~\ref{lm:thekeylemma}(\ref{it:squeeze1}), (\ref{eq:eredsd}) is $\Lv_{G,\mu}$-equivariant.  

For $i\in \bbI_s$, let $R_{f_i,x}$ be the projection of $R_{f,x}$ to $T_{f_i,x}$. Then $R_{f_i,x}$ is a nonzero subrep of the adjoint rep $\ad_{G_i,\mu}:\Lv_{G_i,\mu}\acts T_{f_i,x}$. When $\dim T_{f_i,x}=1$, $\ad_{G_i,\mu}$ is trivially irreducible. When $\dim T_{f_i,x}\geq 3$, the argument in the proof Theorem~\ref{thm:AOsingle} shows that $\ad_{G_i,\mu}$ is irreducible. So in these cases, $R_{f_i,x}=T_{f_i,x}$. When $\dim T_{f_i,x}= 2$. We claim that  $R_{f_i,x}=T_{f_i,x}$, following the idea of Remark~\ref{rmk:exceptionalsplit}: suppose for the sake of contradiction that $\dim {R}_{f_i,x}=1$. One can take $Z$ as in Theorem~\ref{lm:thekeylemma}(\ref{it:squeeze3}). The projection of $Z$ to $\mathscr{A}_{g_i,\Fpbar}$ is an one dimensional special subvariety. Arguing as 
Remark~\ref{rmk:exceptionalsplit}, we find that $\Hg(f_i)$ is not almost simple, so $f$ does not separate factor. This is the desired contradiction. 

Therefore $R_{f,x}$ projects surjectively onto each $T_{f_i,x}$, for $i\in \bbI_s$. Let $J_i\subseteq \Lv_{G_i,\mu}$ be the image of $\mu$ in $G_i$, which is a rank 1 torus. The restriction of $\ad_{G_i,\mu}:\Lv_{G_i,\mu}\acts T_{f_i,x}$ to $J_i$ is pure of weight 1 (i.e. it acts by scaling on $T_{f_i,x}$). Let $\overline{J}_i$ be the image of $J_i$ in $G_i^{\ad}$. Since $G^{\ad}=\prod_{i\in \bbI} G_i^{\ad}$ by Theorem~\ref{thm:MTTllinear2}, we see that $\Lv_{G^{\ad},\mu}$ contains $\prod_{i\in \bbI}\overline{J}_i$. So $R_{f,x}$ is invariant under the natural action $$\prod_{i\in \bbI}\ad_{G_i,\mu}:\prod_{i\in \bbI}{J}_i\acts \bigoplus_{i\in \bbI}T_{f_{i},x}.$$This implies \textit{Claim}~\ref{cm:dfq}.\\

Now, let $\overline{x}$ be the projection of $x$ to $Y$. Since $\mathcal{R}_{f,x}\otimes \bG_m^\wedge\subseteq X^{/x}\subseteq \mathcal{T}_{f,x}\otimes \bG_m^\wedge$ by Theorem~\ref{lm:thekeylemma}(\ref{it:squeeze2}), \textit{Claim}~\ref{cm:dfq} shows that 
$$X^{/x}\supseteq \mathcal{R}_{f,x}\otimes \bG_m^\wedge=(C_s\times \{\overline{x}\})^{/x}.$$ 
So $\overline{X}\supseteq C_s\times \{\overline{x}\}$.
Now let $x$ vary over $X(\Fpbar)$. Since the projection of $X$ to $Y$ is Zariski dense, we find that $\overline{X}\supseteq C_s\times Y$. This implies that (\ref{eq:proddd}) is an equality. 
\end{proof}

\begin{proof}[Proof of two special cases Theorem~\ref{Mainthm:AOSingle} and \ref{Mainthm:AOmodularcurve}]
    From the way that these theorems are stated, we need to use Setup~\ref{setup1}. 
    
    For the single GSpin case, if $\MTT(f)$ is almost simple, then after modification~\ref{constr:varsepa}, we are in the $\#\bbI=1$ case of the Setup~\ref{setup2}. The theorem then follows from Theorem~\ref{thm:AOsingle}. If $\MTT(f)$ is not almost simple, then $\mathcal{X}_f$ must be a product of two special curves, so $\dim X\leq 2$. Since $X$ contains a Zariski dense collection of positive dimensional special subvarieties, $X$ can not be one dimensional (otherwise $X$ is a special curve, and $\mathcal{X}_f$ will be one dimensional). It follows $\dim X=2$. But then by dimension reasons, $X$ is an open subset of $\mathscr{X}_{f,\Fpbar}^{\ord}$. Therefore $X$ is special. 
    
     For the case of products of modular curves, we use Theorem~\ref{thm:AOprod} and Remark~\ref{rmk:quantatative}. The details will be left to the readers. 
\end{proof}

\appendix

\section{Rigidity theorem over a geometric base}\label{subsub:Starshape}
We refers the readers to \S\ref{sub:prepairformalscheme} for conventions. The main goal is to prove Theorem~\ref{lm:rigidity}. Let $\mathfrak{Z}$ be as in Theorem~\ref{lm:rigidity}. Then Chai's rigidity theorem (cf. \cite{Chai08}) implies that for every geometric point $\overline{s}\rightarrow X$, each irreducible component of $\mathfrak{Z}_{\overline{s}}$ is a formal subtorus of $\mathfrak{T}_{\overline{s}}$. A formal subscheme of $\mathfrak{T}$ having this property is called \textbf{star-shaped}, and its structure can be well understood. We will treat this in \S\ref{subsub:srshapsub}.
\subsection{Some commutative algebra} Let $B$ be a Noetherian ring, and let $B[[\underline{t}]]:=B[[t_1,...,t_d]]$ be a power series ring over $B$. It is well-known that $B[[\underline{t}]]$ is Noetherian. Let $f:B\rightarrow D$ be a ring homomorphism, then the completed base change $B[[\underline{t}]]\rightarrow B[[\underline{t}]]\widehat{\otimes}_BD=D[[\underline{t}]]$ can be identified with the natural map extending $f$ and sending $\underline{t}$ to $\underline{t}$.

\begin{lemma}\label{lm:fnfinfn}
    Let $B\rightarrow D$ be finite. Then the natural map $B[[\underline{t}]] \rightarrow D[[\underline{t}]]$ is finite.
\end{lemma}
\begin{proof}
    Let $w_1,...,w_N\in D$ be the generators of $D$ as a $B$-module. We claim that they generate $D[[t_1,...,t_d]]$ as a $B[[t_1,...,t_d]]$-module. Let $f=\sum \delta_{\underline{i}}\,\underline{t}^{\underline{i}}\in D[[t_1,...,t_d]],\,\, \delta_{\underline{i}}\in D$. Writing $\delta_{\underline{i}}=\sum_{k=1}^Nw_k \beta_{\underline{i},k}$ for $\beta_{\underline{i},k}\in B$, we find $f=\sum_{k=1}^Nw_kf_k$, where $ f_k= \sum\beta_{\underline{i},k}\,\underline{t}^{\underline{i}}\in B[[t_1,...,t_d]]$.  
\end{proof}
\begin{lemma}\label{lm:basetogeneric}
    Let $B$ be integral with fraction field $K$. Then the natural map $B[[\underline{t}]]\hookrightarrow K[[\underline{t}]]$ is flat.
\end{lemma}
\begin{proof}
   This is a direct application of \cite[Tag 0523]{stacks-project}. In \textit{loc.cit}, let $R=B[[\underline{t}]]$, $S=M=K[[\underline{t}]]$ and $I=(\underline{t})$. Then all conditions in \textit{loc.cit} are met, and the flatness follows. 
\end{proof}

\subsection{Structure of star-shaped formal subschemes}\label{subsub:srshapsub}
Let $X$ be a scheme over $\Fpbar_p$, and let $\mathfrak{T}$ be a formal torus over $X$. A closed formal subscheme $\mathfrak{Z}\subseteq \mathfrak{T}$ is called \textbf{star-shaped}, if for every geometric point $\overline{s}\rightarrow X$, the irreducible components of $\mathfrak{Z}_{\overline{s}}$ are formal subtori of $\mathfrak{T}_{\overline{s}}$. In the following, we will establish a structure theory for star-shaped formal subschemes. To simplify the question, we will make the following 
\begin{assumption}
    $X=\Spec B$ is a smooth connected affine variety over $\Fpbar$ and  $\mathfrak{T}\simeq \Spf B[[t_1,...,t_d]]$.
\end{assumption}

\begin{lemma}\label{lm:fundamentalsystem}
Assumption as above. Let $\eta$ be the generic point of $X$ and let $\overline{\eta}$ be its algebraic closure. There is a Galois cover $X'\rightarrow X$ and a finite collection of formal subtori $\{\mathfrak{T}_{i,{X'}}\}_{i=1}^l$ of $\mathfrak{T}_{X'}$, such that \begin{enumerate}
    \item\label{it:iier1} $\{\mathfrak{T}_{i,{\overline{\eta}}}\}_{i=1}^l$ is the set of irreducible components of $\mathfrak{Z}_{\overline{\eta}}$,
    \item\label{it:iier2} Each $\mathfrak{T}_{i,{X'}}$ is a closed formal subscheme of $\mathfrak{Z}_{X'}$.
\end{enumerate}        
\end{lemma}
\begin{proof}
 
Let $\eta^{\mathrm{sep}}$ be the separable closure of $\eta$. By definition, each irreducible component of $\mathfrak{Z}_{\overline{\eta}}$ is a formal torus of $\mathfrak{T}_{\overline{\eta}}\simeq \bG^{\wedge,d}_{m,\overline{\eta}}$. Let's call them $\mathfrak{T}_{1,\overline{\eta}}, \mathfrak{T}_{2,\overline{\eta}},...,\mathfrak{T}_{l,\overline{\eta}}$. Note that $\mathfrak{T}_{{\eta}^{\mathrm{sep}}}=\mathbb{G}_{m,{\eta}^{\mathrm{sep}}}^{\wedge,d}$ as well. So each $\mathfrak{T}_{i,\overline{\eta}}$ descends to a subtorus $\mathfrak{T}_{i,{\eta}^{\mathrm{sep}}}\subseteq \mathfrak{Z}_{{\eta}^{\mathrm{sep}}}$ with the same cocharacter lattice, and they form the set of irreducible components of $\mathfrak{Z}_{{\eta}^{\mathrm{sep}}}$. Now $\Gal({\eta}^{\mathrm{sep}}/\eta)$ preserves $\mathfrak{Z}_{{\eta}^{\mathrm{sep}}}$. So there is a Galois extension $\eta''/\eta$ such that 
$\Gal({\eta}^{\mathrm{sep}}/\eta'')$ leaves each $\mathfrak{T}_{i,{\eta}^{\mathrm{sep}}}$ invariant. Then each $\mathfrak{T}_{i,{\eta}^{\mathrm{sep}}}$ further descends to a subtorus $\mathfrak{T}_{i,{\eta'}}\subseteq \mathfrak{Z}_{{\eta''}}$ with the same cocharacter lattice (but with non-trivial Galois action), and they form the set of irreducible components of $\mathfrak{Z}_{{\eta''}}$. Let $X_*(\mathfrak{T}_{i,{\eta''}})\subseteq X_*(\mathfrak{T}_{\eta''})$ be the cocharacter lattice of $\mathfrak{T}_{i,{\eta''}}$.

Recall that $X$ is smooth, so there is a surjection $\pi_1^{\et}(\eta,\overline{\eta})\rightarrow \pi_1^{\et}(X,\overline{\eta})$. Therefore the image of $\pi_1^{\et}(\eta'',\overline{\eta})$ in $\pi_1^{\et}(X,\overline{\eta})$ is of finite index. There exists a finite étale cover $X'\rightarrow X$ corresponding to $\im \pi_1^{\et}(\eta'',\overline{\eta})$ so that each $X_*(\mathfrak{T}_{i,{\eta'}})$ spreads out to a saturated lisse subsheaf $\mathcal{H}_{i,X'}\subseteq X_*(\mathfrak{T}_{X'})$. Possibly enlarging $X'$, we can assume that it is Galois. As a result, for each $i$, we obtain a formal subtorus $\mathfrak{T}_{i,{X'}}\subseteq \mathfrak{T}_{X'}$ with $\mathcal{H}_{i,X'}=X_*(\mathfrak{T}_{i,{X'}})$. Clearly, the base change of each $\mathfrak{T}_{i,{X'}}$ to $\overline{\eta}$ is $\mathfrak{T}_{i,\overline{\eta}}$. This proves (\ref{it:iier1}). 

Let $\eta'$ be the generic point of $X'$. To show that $\mathfrak{T}_{i,{X'}}$ is a formal subscheme of $\mathfrak{Z}_{X'}$, it suffices to show that $\mathfrak{T}_{i,{X'}}$ is the smallest closed formal subscheme of $\mathfrak{T}_{{X'}}$ whose base change to $\eta'$ contains $\mathfrak{T}_{i,{\eta'}}$. Suppose that there is a closed formal subscheme $\mathfrak{Q}_{i,X'}\subseteq \mathfrak{T}_{i,{X'}}$ whose base change to $\eta'$ contains $\mathfrak{T}_{i,{\eta'}}$, we will show that $\mathfrak{Q}_{i,X'}=\mathfrak{T}_{i,{X'}}$. It suffices to check that for any $x\in X'(\Fpbar)$, we have $\mathfrak{Q}_{i,x}=\mathfrak{T}_{i,x}$. Let ${X'}^{/x}=\Spf B_x$. Then $\mathfrak{T}_{i,\Spec B_x}$ is a trivial tori, so one can choose coordinates so that  $\mathfrak{T}_{i,\Spec B_x}=\Spf B_x[[\underline{s}]]$. Then $\mathfrak{Q}_{i,\Spec B_x}$ corresponds to an ideal $\mathfrak{a}\subseteq B_x[[\underline{s}]]$. Let $K_x$ be the ring of fractions of $B_x$. Then our condition translates to $\mathfrak{a}K_x[[\underline{s}]]=(0)$. Since $B_x[[\underline{s}]]\subseteq K_x[[\underline{s}]]$, this implies that $\mathfrak{a}=(0)$. A further base change to the point $x$ implies that $\mathfrak{Q}_{i,x}=\mathfrak{T}_{i,x}$, as desired. This proves (\ref{it:iier2}).
\end{proof}

\begin{remark}\label{rmkvrmk}
    Let $X'\rightarrow X$ and $\{\mathfrak{T}_{i,{X'}}\}_{i=1}^l$ be as in Lemma~\ref{lm:fundamentalsystem}. Let $\eta$ \textit{resp}. $\eta'$ be the generic point of $X$ \textit{resp}. $X'$. Let $X'=\Spec B', \eta=\Spec K$ and $\eta'=\Spec K'$. Then each $\mathfrak{T}_{i,{\eta'}}$ is cut out by a prime $\mathcal{P}'_i\subseteq K'[[\underline{t}]]$, and $\mathfrak{T}_{i,{X'}}$ is cut out by $\mathcal{P}'_i\cap B'[[\underline{t}]]$. 

Let $G=\Gal(X'/X)$. Then it acts freely on $\mathfrak{T}_{X'}$ (\textit{resp}. $\mathfrak{T}_{\eta'}$) and $\mathfrak{T}_{X'}/G= \mathfrak{T}$ (\textit{resp}. $\mathfrak{T}_{\eta'}/G= \mathfrak{T}_{\eta}$). Then $G$ permutes elements in $\{\mathfrak{T}_{i,{X'}}\}_{i=1}^l$ (\textit{resp}. $\{\mathfrak{T}_{i,{\eta'}}\}_{i=1}^l$), and each orbit descends to an irreducible star-shaped formal subscheme of $\mathfrak{Z}$ (\textit{resp}.  $\mathfrak{Z}_{\eta}$). This can be seen more clearly via algebra: By Lemma~\ref{lm:fnfinfn}, $B[[\underline{t}]]\rightarrow B'[[\underline{t}]]$ (\textit{resp}. $K[[\underline{t}]]\rightarrow K'[[\underline{t}]]$ ) is a finite extension, and the group $G$ acts on $B'[[\underline{t}]]$ (\textit{resp}. $K'[[\underline{t}]]$) via coefficients. The $G$-action permutes the ideals $\{\mathcal{P}'_i\cap B'[[\underline{t}]]\}_{i=1}^{l}$ (\textit{resp}. $\{\mathcal{P}'_i\}_{i=1}^{l}$). For each $i$, the contraction $\mathcal{P}_i'\cap B[[\underline{t}]]$ (\textit{resp}. $\mathcal{P}'_i\cap K[[\underline{t}]]$) cuts out the irreducible star-shaped formal subscheme of $\mathfrak{Z}$ (\textit{resp}.  $\mathfrak{Z}_{\eta}$) obtained by descending the $G$-orbit containing $\mathfrak{T}_{i,{X'}}$ (\textit{resp}. $\mathfrak{T}_{i,{\eta'}}$). 

\end{remark}

\begin{defn}\label{def:bone}
    An irreducible star-shaped formal subscheme of $\mathfrak{Z}$ (\textit{resp}.  $\mathfrak{Z}_{\eta}$) that arises from descending a $G$-orbit of formal subtori as described in the last paragraph is called a \textbf{bone} of $\mathfrak{Z}$ (\textit{resp}.  $\mathfrak{Z}_{\eta}$). 
\end{defn}
\begin{lemma}\label{lm:rquiconditions}
   Let $I(\mathfrak{Z})\subseteq B[[\underline{t}]]$ be the ideal of $\mathfrak{Z}$. The following are equivalent: \begin{enumerate}
        \item\label{lm:rquiconditions1} $I(\mathfrak{Z})\subseteq (\underline{t})$,
        \item\label{lm:rquiconditions2} $\mathfrak{Z}_x\neq \emptyset$ for all $x\in X(\Fpbar)$, 
        \item\label{lm:rquiconditions3}  $\mathfrak{Z}_x\neq \emptyset$ for Zariski dense $x\in X(\Fpbar)$.
    \end{enumerate}
\end{lemma}
\begin{proof}
Suppose that (\ref{lm:rquiconditions1}) holds. Let $\mathfrak{m}\subseteq B$ be the maximal ideal corresponding to  $x$. Then $I(\mathfrak{Z}_x)\subseteq (\underline{t})\frac{B}{\mathfrak{m}}[[\underline{t}]]$, so $\mathfrak{Z}_x$ is not empty. So (\ref{lm:rquiconditions2}) holds. (\ref{lm:rquiconditions2}) $\Rightarrow$ (\ref{lm:rquiconditions3})  is obvious. Suppose that (\ref{lm:rquiconditions3}) holds, then $X$ is a closed formal subscheme of $\mathfrak{Z}$. The ideal cutting out $X$ is $(\underline{t})$, so (\ref{lm:rquiconditions1}) holds.
\end{proof}
\begin{theorem}\label{thm:starshapedmaintheorem}
  Suppose that $\mathfrak{Z}$ is irreducible, and satisfies the equivalent conditions of Lemma~\ref{lm:rquiconditions}, then $\mathfrak{Z}$ has a unique bone and it equals this bone. In particular, there is a finite Galois cover $X'\rightarrow X$ such that $\mathfrak{Z}_{X'}$ is a union of formal subtori of $\mathfrak{T}_{X'}$ which are permuted transitively by $\Gal(X'/X)$.
\end{theorem}
\begin{proof}
Let the notation be the same as Remark~\ref{rmkvrmk}.  By Lemma~\ref{lm:basetogeneric}, $B[[\underline{t}]]\hookrightarrow K[[\underline{t}]]$ is flat. The condition of Lemma~\ref{lm:rquiconditions} implies that $I(\mathfrak{Z})\subseteq (\underline{t})$. So $I(\mathfrak{Z})K[[\underline{t}]]$ is contained in the maximal ideal. By flat going down, there is a prime $\mathcal{P}\subseteq K[[\underline{t}]]$ such that $I(\mathfrak{Z})=\mathcal{P}\cap B[[\underline{t}]]$. It follows that  $\mathcal{P}\supseteq I(\mathfrak{Z})K[[\underline{t}]]$, so $\mathcal{P}$ contains a minimal prime of 
$I(\mathfrak{Z})K[[\underline{t}]]$. 

Such a minimal prime is of form $\mathcal{P}'_i\cap K[[t]]$ for some $1\leq i\leq l$ (i.e., it is the minimal prime that cuts out a bone of $\mathfrak{Z}_{\eta}$). The contraction $\mathcal{P}'_i\cap B[[t]]$ is a bone of $\mathfrak{Z}$. On the other hand, $$I(\mathfrak{Z})\subseteq \mathcal{P}'_i\cap B[[t]]\subseteq \mathcal{P}\cap B[[\underline{t}]]= I(\mathfrak{Z}).$$
So $\mathcal{P}'_i\cap B[[t]]=I(\mathfrak{Z})$. Therefore $\mathfrak{Z}$ contains a unique bone which is cut out by $\mathcal{P}'_i\cap B[[t]]$, and $\mathfrak{Z}$ equals this bone. The last assertion follows from definition of the bones. 
\end{proof}
\begin{proof}[Proof of Theorem~\ref{lm:rigidity}]
Since $\mathfrak{Z}$ is invariant under infinitely many scaling automorphism $[\alpha]$, for each geometric point $\overline{s}\rightarrow X$, $\mathfrak{Z}_{\overline{s}}$ is invariant under infinitely many scaling automorphism $[\alpha]$. Chai's rigidity theorem (cf. \cite{Chai08}) implies that each irreducible component of $\mathfrak{Z}_{\overline{s}}$ is a formal subtorus of $\mathfrak{T}_{\overline{s}}$. So $\mathfrak{Z}$ is star-shaped. Now apply Theorem~\ref{thm:starshapedmaintheorem} to conclude. 
\end{proof}
\bibliographystyle{alpha}
\bibliography{ref}

\end{document}